\RequirePackage{fix-cm}
\documentclass[smallextended]{svjour3}       
\smartqed  
\usepackage{array}
\usepackage{color}
\usepackage{xcolor,colortbl}
\usepackage{tabularx}
\usepackage{graphicx}
\usepackage{amsmath}
\usepackage{amssymb}
\usepackage{empheq}
\usepackage{multirow}
\usepackage{subfig}
\usepackage{moreverb}
\usepackage{nicefrac} 
\usepackage{dsfont}
\usepackage{bm}
\usepackage{multirow}
\usepackage{soul}

\usepackage{multirow} 
\usepackage{longtable,tabu} 
\usepackage{hhline}
\usepackage[linesnumbered,ruled,vlined]{algorithm2e}

\newcommand\BibTeX{{\rmfamily B\kern-.05em \textsc{i\kern-.025em b}\kern-.08em
		T\kern-.1667em\lower.7ex\hbox{E}\kern-.125emX}}

\DeclareMathOperator\erf{erf}

\newcommand{\RR}{p}

\newcommand{\Q}{\mathbf{Q}}

\newcommand{\LL}{\mathcal{L}}

\newcommand{\diff}[2]{\frac{\partial {#1} }{\partial {#2} } }

\usepackage{pgf,tikz}
\usetikzlibrary{arrows}

\definecolor{ttzzqq}{rgb}{0.2,0.6,0}

\definecolor{qqttcc}{rgb}{0,0.2,0.8}
\definecolor{ffxfqq}{rgb}{1,0.5,0}
\definecolor{qqttzz}{rgb}{0,0.2,0.6}
\definecolor{ffqqqq}{rgb}{1,0,0}
\definecolor{qqwuqq}{rgb}{0,0.39,0}
\definecolor{zzttqq}{rgb}{0.6,0.2,0}
\definecolor{qqqqff}{rgb}{0,0,1}
\definecolor{ttttqq}{rgb}{0.2,0.2,0}

\definecolor{qqwwtt}{rgb}{0,0.4,0.2}
\definecolor{ubqqys}{rgb}{0.29,0,0.51}

\definecolor{wwttqq}{rgb}{0.4,0.2,0}
\definecolor{uuuuuu}{rgb}{0.27,0.27,0.27}
\definecolor{qqzzff}{rgb}{0,0.6,1}
\definecolor{xdxdff}{rgb}{0.49,0.49,1}
\definecolor{ccwwqq}{rgb}{0.8,0.4,0}
\definecolor{ttttff}{rgb}{0.2,0.2,1.}
\definecolor{ttqqqq}{rgb}{0.2,0,0}
\definecolor{yqyqyq}{rgb}{0.5019607843137255,0.5019607843137255,0.5019607843137255}
\definecolor{qqzzcc}{rgb}{0,0.6,0.8}
\definecolor{ffqqtt}{rgb}{1.,0.,0.2}
\definecolor{qqccqq}{rgb}{0.,0.8,0.}

\newcommand{\figureAdv}{
	\begin{tikzpicture}[line cap=round,line join=round,>=triangle 45,x=1.0cm,y=1.0cm,scale=0.8]
		\clip(4.738508035803117,4.6178450856057305) rectangle (15.2557003259263,12.126145276832727);
		\draw [->,line width=2.pt] (14.,9.) -- (11.,10.);
		\draw [->,line width=2.pt] (11.,10.) -- (8.,9.);
		\draw [->,line width=2.pt] (8.,9.) -- (6.,7.);
		\draw (14.043785452759789,9.117253177936565) node[anchor=north west] {$x$};
		\draw (10.630921729474782,9.883406258674015) node[anchor=north west] {$x_1^E$};
		\draw (11.95427705074856,10.705279563465096) node[anchor=north west] {$u_1=\mathcal{H}^L_1(x)$};
		\draw (7.984211086927226,8.950092505775666) node[anchor=north west] {$x_2^E$};
		\draw (6.047933301063487,7.125255168019199) node[anchor=north west] {$x_3^E$};
		\draw (8.011063607530594,10.635629283398055) node[anchor=north west] {$u_2=\mathcal{H}^L_2(x)$};
		\draw (4.862870948500343,8.936162449762259) node[anchor=north west] {$u_3=\mathcal{H}^L_3(x)$};
		\begin{scriptsize}
			\draw [fill=black] (14.,9.) circle (3.5pt);
			\draw [fill=black] (11.,10.) circle (3.5pt);
			\draw [fill=black] (8.,9.) circle (3.5pt);
			\draw [fill=black] (6.,7.) circle (3.5pt);
		\end{scriptsize}
	\end{tikzpicture}
}

\newfont{\numerikEleven}{ecrm1000}
\newfont{\numerikTen}{cmss10}
\newfont{\numerikNine}{cmss9}
\newfont{\numerikEight}{cmss8}

\journalname{Journal of Scientific Computing}
\begin{document}
	
	\title{On the construction of conservative semi-Lagrangian IMEX advection schemes for multiscale time dependent PDEs}
	
	\titlerunning{Conservative semi-Lagrangian IMEX advection schemes for multiscale time dependent PDEs}        
	
	\author{Walter Boscheri \and
		Maurizio Tavelli \and
		Lorenzo Pareschi
	}
	
	
	\institute{W. Boscheri \at
		Department of Mathematics and Computer Science, University of
		Ferrara, Ferrara, Italy \\
		\email{walter.boscheri@unife.it}
		\and
		M. Tavelli \at
		Faculty of Science and Technology, Free University of Bozen, Bolzano, Italy \\
		\email{mtavelli@unibz.it}
		\and
		L. Pareschi \at
		Department of Mathematics and Computer Science, University of
		Ferrara, Ferrara, Italy \\
		\email{lorenzo.pareschi@unife.it}
	}
	
	\date{Received: date / Accepted: date}

	\maketitle
	
	\begin{abstract}
		This article is devoted to the construction of a new class of semi-Lagrangian (SL) schemes with implicit-explicit (IMEX) Runge-Kutta (RK) time stepping for PDEs involving multiple space-time scales. The semi-Lagrangian (SL) approach fully couples the space and time discretization, thus making the use of RK strategies particularly difficult to be combined with. First, a simple scalar advection-diffusion equation is considered as a prototype PDE for the development of a high order formulation of the semi-Lagrangian IMEX algorithms. The advection part of the PDE is discretized explicitly at the aid of a SL technique, while an implicit discretization is employed for the diffusion terms. In this way, an unconditionally stable numerical scheme is obtained, that does not suffer any CFL-type stability restriction on the maximum admissible time step. Second, the SL-IMEX approach is extended to deal with hyperbolic systems with multiple scales, including balance laws, that involve shock waves and other discontinuities. A conservative scheme is then crucial to properly capture the wave propagation speed and thus to locate the discontinuity and the plateau exhibited by the solution. A novel SL technique is proposed, which is based on the integration of the governing equations over the space-time control volume which arises from the motion of each grid point. High order of accuracy is ensured by the usage of IMEX RK schemes combined with a Cauchy-Kowalevskaya procedure that provides a predictor solution within each space-time element. The one-dimensional shallow water equations (SWE) are chosen to validate the new conservative SL-IMEX schemes, where convection and pressure fluxes are treated explicitly and implicitly, respectively. The asymptotic-preserving (AP) property of the novel schemes is also studied considering a relaxation PDE system for the SWE. A large suite of convergence studies for both the non-conservative and the conservative version of the novel class of methods demonstrates that the formal order of accuracy is achieved and numerical evidences about the conservation property are shown. The AP property for the corresponding relaxation system is also investigated.
		\subclass{MSC 65 \and MSC 68}
	\end{abstract}
	
	
	\section{Introduction} \label{sec:introduction}
	Dynamic processes in continuum physics are modeled using time-dependent partial differential equations (PDE), which are based on the conservation of some physical quantities, such as mass, momentum and energy. 
	The governing equations may involve different physical processes, like advection, diffusion, pressure gradients and drag forces among many others. The time scale associated to each process is not the same, for instance diffusion processes are defined at a much smaller time scale than advection phenomena, or pressure waves travel much faster than material interfaces or contact discontinuities.  
	Specifically, numerical schemes that are able to handle multiple time scales simultaneously are extremely important for real world applications, as well as algorithms which are numerically stable for a wide range of admissible time steps while advancing the solution in time.
	
	These observations led to the idea of splitting the processes on fast and slow time scales and treating them with different numerical techniques. The main strategy consists in treating implicitly only one part of the system to be solved while keeping the remaining explicit \cite{DegTan,Bog,DLDV2018,ToroVazquez}, hence allowing space and time discretizations to be designed, in which the implicit part of the system is relatively easy to be inverted, typically avoiding nonlinear systems, while keeping robustness and shock-capturing properties in the explicit part. 
	The most popular approach is based on implicit-explicit (IMEX) methods \cite{AscRuuSpi,PR_IMEX,BP2021} that have proven to be very successful in designing asymptotic-preserving (AP) schemes capable to deal with stiff source terms. Such IMEX schemes typically may achieve high order of accuracy under a time step stability constraint independent of the values of the fast scale. Alternatively, other semi-implicit methods have been proposed \cite{ParkMunz2005,Klein,Casulli1999,CasulliCattani,BosRus2019} where a linearly implicit scheme is derived for the stiff terms in the governing equations, thus avoiding any need of iterative solvers. 
	Let us mention that, semi-implicit hybrid finite volume/finite element schemes have been recently proposed in \cite{Hybrid2,Hybrid1}, while semi-implicit methods coupled with discontinuous Galerkin (DG) space discretizations on unstructured staggered meshes have been forwarded for compressible flows \cite{TavelliCNS}, on dynamic adaptive meshes \cite{Fambri2017b} and for axially symmetric flows \cite{Ioriatti2018}. 
	In most of the aforementioned works, the convective terms of the governing equations are discretized explicitly, because they typically involve a nonlinearity which is difficult to be implicitly solved, requiring the usage of computationally time consuming and numerically less stable nonlinear solvers for the resulting system that need to be inverted. Explicit upwind finite difference and Godunov-type finite volume methods are very popular \cite{Godunov1959, LeVequeBook, Munz1994, OS1982}, although these schemes are restricted to a Courant-Friedrichs-Lewy (CFL) stability condition based on the maximum eigenvalues of the Jacobian matrix associated to the hyperbolic system. Another option is given by the so-called semi-Lagrangian (SL) methods, which have recently achieved visibility due to their excellent resolution and stability properties. In the SL framework, the advection term is written in a Lagrangian formulation, which is then discretized accordingly. In particular, semi-Lagrangian schemes require the integration of the material trajectories backward in
	time to find the foot of the characteristic, where the numerical solution is interpolated. These methods have been originally developed for the numerical weather prediction \cite{SL_29,SL_30}. Nowadays, they can be found in environmental engineering applications, such as free surface flows in rivers and oceans \cite{SL_39,SL_40,ADERFSE} as well as in plasma physics  \cite{Vlasov2010} and kinetic equations \cite{Carrillo2007,SLBoltz12,SLkinI_Bos2021}, in applications to image processing \cite{Carlini2013} or for solving the Hamilton-Jacobi equations \cite{Falcone2002}. The semi-Lagrangian approach is the only explicit method for the discretization of convective terms that allows for large time steps without imposing a CFL-type stability condition, therefore it constitutes a very interesting alternative to standard explicit upwind solvers. However, the extension of SL methods to deal with shock waves that occur in hyperbolic PDE is not trivial because conservation must be strictly ensured. Conservative semi-Lagrangian methods are described for instance in \cite{LentineEtAl2011,EulLagWENO2012}. In \cite{Shu_ConsSL2011} a conservative WENO finite difference scheme with SL treatment of advection is proposed for incompressible flows, while recently in \cite{ConsSL2021} a novel approach for the development of conservative semi-Lagrangian schemes has been introduced, which is based on the backward integration of the scalar advection equation onto a space-time control volume defined for each grid point. Let us also recall that implicit-explicit SL discretizations have been considered  in \cite{IMEXSL08} for convection dominated problems, and in \cite{CaiBos21} Runge-Kutta exponential integrators have been coupled with SL discretizations for nonlinear Vlasov equations.
	
	State-of-the-art semi-Lagrangian methods are mostly used in non-conservative form, for instance in the discretization of the advection contribution for incompressible flows. On the other hand, when kinetic equations with linear transport are considered, the use of conservative SL schemes is adopted, where the characteristics are straight lines and therefore can be easily followed for designing high order methods. If IMEX discretizations based on Runge-Kutta time integrators are used, semi-Lagrangian schemes exhibit an intrinsic difficulty which arises from the fact that space and time cannot be decoupled, since SL methods directly solve the Lagrangian form of advection and, as such, no purely spatial advection fluxes can be retrieved as in standard RK time stepping techniques. This makes the application of semi-Lagrangian schemes to RK integrators not very common in the literature.  A recent attempt for the development of conservative SL methods for hyperbolic PDE has been forwarded in \cite{ConsSL2021}, which is restricted to the case of the scalar nonlinear advection PDE.
	
	In this work we aim at constructing a new class of methods that deal with both semi-Lagrangian schemes for advection and IMEX-RK time stepping at the same time. In the first part of the paper, we derive two different algorithms for the non-conservative case and we demonstrate the accuracy and the robustness of the proposed approach by considering the scalar advection-diffusion equation as a prototype, where advection and diffusion are treated explicitly and implicitly, respectively. The second part of the article is devoted to the design of a conservative semi-Lagrangian IMEX scheme for hyperbolic systems of balance laws, namely the one-dimensional shallow water equations (SWE). The 1D SWE represent a simple, but not trivial, prototype example of a hyperbolic system, since they involve both advection and pressure contributions that are suitable for an IMEX scheme with SL treatment of the convective terms. Indeed, a fast scale related to the pressure terms and a slow scale associated with convection coexist. Furthermore, shock waves are part of the eigenstructure of the system, thus requiring a conservative method to properly capture the wave speed and eventually its location. The novel class of methods is addressed with SL-IMEX schemes and is shown to achieve up to third order of accuracy in space and time, while ensuring fully conservation of the state variables. Finally, an asymptotic-preserving SL-IMEX scheme is developed by considering a relaxation system for the SWE, and numerical evidences are proposed.
	
	The outline of this article is as follows. In Section \ref{sec:model} we present the scalar advection-diffusion equation, while Section \ref{sec.EulIMEX} revisits the classical IMEX schemes with Eulerian discretization of the advection terms. Two different non-conservative versions of the novel SL-IMEX algorithms are derived in Section \ref{sec.SL-IMEX} and numerical results for this class of methods are shown in Section \ref{sec.AdvDiff_test}. The second part of the article starts with the introduction of the shallow water model in Section \ref{sec.SWE}, which is followed by a detailed description of the new conservative SL-IMEX schemes. Next, in Section \ref{sec.AP}, the schemes are extended to deal with the presence of source terms  and their AP properties are analyzed. The novel algorithms are validated through a suite of numerical tests that is then proposed in Section \ref{sec.SWE_test}. Finally, conclusions and outlook to future work are given in Section \ref{sec.concl}.
	
	\section{Semi-Lagrangian IMEX schemes for advection-diffusion equations} \label{sec:model}
	
	Let us consider the one-dimensional advection-diffusion equation of a scalar quantity $q=q(x,t)$ over a given velocity field $u(x,t)$ with a constant diffusion coefficient $\alpha \in \mathds{R}^+$: 
	\begin{equation}
		\diff{q}{t} + \diff{u q}{x} = \alpha \frac{\partial^2 q}{\partial x^2},
		\label{sys1.1}
	\end{equation}
	where $t \in \mathds{R}_0^+$ is the time and $x \in \mathds{R}$ represents the spatial coordinate. In the semi-Lagrangian framework, the advection term in \eqref{sys1.1} can be reformulated using the Lagrangian derivative, thus yielding
	\begin{equation}
		\frac{d q}{d t} = \alpha \frac{\partial^2 q}{\partial x^2}.
		\label{sys1.2}
	\end{equation} 
	The solution of both \eqref{sys1.1} and \eqref{sys1.2} propagates along characteristics, that are defined by the trajectory equation
	\begin{equation}
		\frac{d x}{d t} = u(x,t).
		\label{eqn.TrODE}
	\end{equation}
	In the above equation \eqref{sys1.2}, the classical semi-Lagrangian technique is used, with the diffusion terms treated in Eulerian form. The inclusion of the diffusion terms in the semi-Lagrangian approach can be found in \cite{Bonaventura2014,Bonaventura2020}. Despite its simplicity, the governing PDE \eqref{sys1.1} represents a suitable prototype for deriving semi-Lagrangian methods with IMEX time stepping because: i) it contains an advection term which can be discretized at the aid of a semi-Lagrangian scheme; ii) it involves a diffusion contribution that might be conveniently solved implicitly in order to relax the severe stability condition on the time step imposed by parabolic terms. Therefore, the new class of semi-Lagrangian IMEX schemes, denoted as SL-IMEX, will be firstly constructed by studying the advection-diffusion equation \eqref{sys1.1}, or equivalently \eqref{sys1.2}.
	
	\subsection{IMEX schemes with Eulerian advection} \label{sec.EulIMEX}
	The governing PDE \eqref{sys1.1} undergoes an implicit-explicit time discretization, where the nonlinear convective contribution is taken explicitly while the diffusion terms are considered implicitly. Therefore, a first order in time semi-discrete scheme writes
	\begin{equation}
		{\frac{q^{n+1} - q^n}{\Delta t}} + \diff{(u q)^n}{x} = \alpha \frac{\partial^2 q^{n+1}}{\partial x^2},
		\label{sys1.3}
	\end{equation}
	with the time step defined as $\Delta t = t^{n+1}-t^n$. Let $\mathcal{H}(q_E(t),q_I(t))$ denote a spatial approximation of the explicit advection term with $q_E(t)$ and of the implicit diffusion term with $q_I(t)$ , that is
	\begin{equation}
		\mathcal{H}(q_E(t),q_I(t)) = - \diff{(u_E q_E)}{x} + \alpha \frac{\partial^2 q_I}{\partial x^2}.
		\label{eqn.H}
	\end{equation}
	Although the time discretization \eqref{sys1.3} could have been written in a more canonical partitioned form as 
	\[
	\mathcal{H} = \mathcal{H}^E(q_E(t))+\mathcal{H}^I(q_I(t)),
	\]
	 we prefer to write the PDE under the form of an autonomous system following \cite{BosFil2016}. This easily allows for a more general framework, where nonlinearities in the implicit part of the PDE can be properly treated, contrarily to the partitioned form \cite{ODEstiff_book}. Furthermore, the scope of the advection-diffusion equation \eqref{sys1.1} is only to furnish an example of application for the development of the SL-IMEX schemes discussed in this work, which in principle might be extended to other physical models.
	Notice that any Eulerian discretization of the advection term based on flux form fits the formalism \eqref{eqn.H}. The definition \eqref{eqn.H} allows the semi-discrete form \eqref{sys1.3} to be written as
	\begin{equation}
		\begin{aligned}
			\frac{d q_E(t)}{dt} &= \mathcal{H}(q_E(t),q_I(t)), \\
			\frac{d q_I(t)}{dt} &= \mathcal{H}(q_E(t),q_I(t)).
		\end{aligned}
		\label{sys1.4_tmp}
	\end{equation}
	In the above formulation \eqref{sys1.4_tmp}, the number of unknowns has been formally doubled, which is not the case if the PDE is written in partitioned form. However, for a specific choice of time discretizations and for autonomous systems this duplication is indeed only apparent \cite{BosFil2016}, thus yielding
		\begin{equation}
			\frac{d q(t)}{dt} = \mathcal{H}(q_E(t),q_I(t)).
			\label{sys1.4}
		\end{equation}
		
	The above formulation \eqref{sys1.4} can be readily applied to the framework of classical IMEX Runge-Kutta (IMEX-RK) schemes \cite{PR_IMEX,BosFil2016}. These multi-step methods 
	are based on $s$ stages and they are typically described in terms of the double Butcher tableau:
	\begin{equation}
		\begin{array}{c|c}
			\tilde{c} & \tilde{A} \\ \hline \\[-0.3cm] & \tilde{b}^\top
		\end{array} \qquad
		\begin{array}{c|c}
			c & A \\ \hline \\[-0.3cm] & b^\top
		\end{array},
		\label{eqn.butcher}
	\end{equation}
	with the matrices $(\tilde{A},A) \in \mathds{R}^{s \times s}$ and the vectors $(\tilde{c},c,\tilde{b},b) \in \mathds{R}^s$. The tilde symbol refers to the explicit scheme and matrix $\tilde{A}=(\tilde{a}_{ij})$ is a lower triangular matrix with zero elements on the diagonal, while $A=({a}_{ij})$ is a triangular matrix which accounts for the implicit scheme, thus having non-zero elements on the diagonal. From now on we adopt IMEX-RK schemes with $b=\tilde{b}$ and the Butcher tableau of the schemes used in this work are reported in Appendix \ref{app.IMEX}. A general IMEX-RK scheme aims at computing the numerical solution at the next time step $q^{n+1}$ starting from $q^n$, and it can be compactly written as follows:
	\begin{itemize}
		\item Stage values for $i=1 \ldots s$
		\begin{subequations}
			\label{eqn.IMEXRK}
			\begin{align}
				q_E^{(i)}&=q^{n}+\Delta t \sum\limits_{j=1}^{i-1}\tilde{a}_{i,j} \, \mathcal{H}(q_E^{ (j)},q_I^{ (j)}), \label{e1} \\
				q_{*}^{(i)}&=q^n  + \Delta t \sum\limits_{j=1}^{i-1}a_{i,j} \, \mathcal{H}(q_E^{ (j)},q_I^{ (j)}), \label{es} \\
				q_I^{(i)}&=q_{*}^{(i)} +\Delta t \, a_{i,i} \, \mathcal{H}(q_E^{ (i)},q_I^{ (i)}). \label{e2}
			\end{align}
		\end{subequations}
		The solution of \eqref{e2} involves an implicit evaluation which corresponds to the backward Euler scheme \eqref{sys1.3}.
		
		\item Update of the numerical solution in terms of the spatial fluxes evaluated at the previous stages
		\begin{equation}
			q^{n+1}=q^n+ \Delta t \sum\limits_{i=1}^{s} b_i \mathcal{H}(q_E^{ (i)},q_I^{ (i)}),	
			\label{e4}
		\end{equation}	
		where we recall that $b=\tilde{b}$, see Appendix \ref{app.IMEX}.
	\end{itemize}
	
	The above approach is then complemented with a suitable space discretization which can be designed independently of the time discretization by adopting, for example, a hybrid finite volume/finite difference scheme  for the discretization of the spatial flux $\mathcal{H}$. Specifically, a finite volume scheme designed for the hyperbolic advection part can be combined with a centered finite-difference scheme for the diffusion terms (see for example \cite{AscRuuSpi,BPR2017,BP2021}).   

	\subsection{IMEX schemes with semi-Lagrangian advection} \label{sec.SL-IMEX}
	Let us now consider a semi-Lagrangian discretization of the advection term in the governing PDE \eqref{sys1.1}, thus the aim is to solve the transport part of the equation by following the characteristics which move with velocity $u(x,t)$ according to \eqref{eqn.TrODE}. A semi-Lagrangian method is based on two main steps:
	\begin{enumerate}
		\item backward integration of the Lagrangian trajectories, which is nothing but the integration of the ODE \eqref{eqn.TrODE}, in order to find the foot of the characteristic located at $x^L$;
		\item reconstruction, or interpolation, of the transported quantity	$q$ at the point $x^L$, that in principle never lies on a known grid point.
	\end{enumerate} 
	Then, neglecting the diffusion part, a direct discretization of the Lagrangian form of the PDE \eqref{sys1.2} simply yields
	\begin{equation}
		q^{n+1} = q^L,
		\label{eqn.SL}
	\end{equation}
	with $q^L:=q(x^L)$ being the solution interpolated at the foot of the characteristic $x^L$. The semi-Lagrangian method \eqref{eqn.SL} is the only explicit discretization that is not constrained by a CFL-type stability condition, hence allowing for more stability and accuracy compared to classical Eulerian methods.
	The main problem when the semi-Lagrangian approach is considered in combination with an IMEX time discretization is that the solution along the characteristics directly solves the contribution of the Lagrangian derivative $\frac{dq}{dt}$, instead of handling $\left(\diff{q}{t} + \diff{u q}{x}\right)$ separately. As a consequence, it is not possible to split the spatial discretization from the temporal one as required in the IMEX formalism, hence implying that the PDE can no longer be cast into form \eqref{sys1.4}. 
	
	In order to derive a consistent semi-Lagrangian IMEX scheme, labeled with SL-IMEX, we proceed in three steps. First, a pure advection problem is considered. Second, a pure diffusion problem is analyzed. Third, we introduce advection and diffusion simultaneously to solve the governing PDE \eqref{sys1.1}.

	\subsubsection{Advection dominated SL explicit RK schemes}
	For an advection dominated problem, the diffusive terms are neglected ($\alpha=0$), so that there is no need of computing any implicit contribution in \eqref{e2}. Let us define a general semi-Lagrangian operator $\mathcal{L}$ such that
	\begin{equation}
		q^L(x) = q(x^L,t) = \mathcal{L}(q(x,t),\Delta t^L,u(x,t)),
		\label{eqn.Lop}
	\end{equation} 
	which solves the trajectory equation \eqref{eqn.TrODE} with the velocity field $u(x,t)$ over a time interval of $\Delta t^L=\omega \Delta t$ with $\omega \in \mathds{R}$. The result of \eqref{eqn.Lop} is the interpolated quantity $q^L(x)$ at the foot of the characteristic $x^L=x-u(x,t)\Delta t^L$, which is exactly what is required by the semi-Lagrangian scheme \eqref{eqn.SL}. The order of accuracy of the Lagrangian operator depends on both the accuracy of the ODE solver for \eqref{eqn.TrODE} and the reconstruction operator, thus the simple definition $x^L=x-u(x,t)\Delta t^L$ is only first order accurate.
	
\begin{figure}[!htbp]
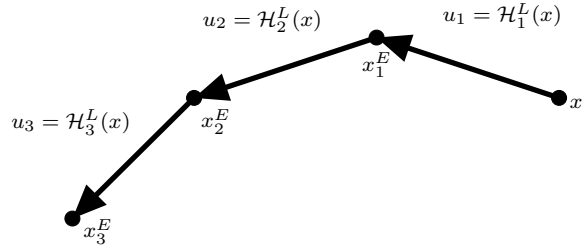

		\begin{center}
			\begin{tabular}{c} 
				\figureAdv      
			\end{tabular} 
			\caption{Example of piecewise trajectory solved with SA-SSP(3,3,2) given by \eqref{eqn.IMEX2}.}
			\label{fig.1}
		\end{center}
	\end{figure}

	Equation \eqref{eqn.Lop} suggests to evolve the coordinates $x$ instead of $q$ itself, meaning that the trajectory ODE \eqref{eqn.TrODE} is the equation to be solved with the RK discretization and not the Lagrangian PDE \eqref{sys1.2}. This is equivalent to consider an augmented PDE system for $\mathbf{q}=(q,x)^\top$, with the fluxes $\mathcal{H}=(0,\mathcal{H}^L)^\top=(0,-u(x,t))^\top$:
	\begin{equation}
			\left\{\begin{aligned}
				\frac{d q}{d t} &= 0, \\
				\frac{d x}{d t} &= \mathcal{H}^L, \qquad \mathcal{H}^L = -u(x,t). \\
			\end{aligned}\right. 
			\label{eqn.advSubSyst}
	\end{equation}
	First, the trajectory equation in \eqref{eqn.advSubSyst} must be solved using the explicit Butcher tableau in \eqref{eqn.butcher}, that is
	\begin{eqnarray}
		x_E^{(i)}=x_E^{n}+\Delta t \sum\limits_{j=1}^{i-1}\tilde{a}_{i,j} \mathcal{H}_j^L, \qquad i=1 \ldots s, 
		\label{eqn.xE}
	\end{eqnarray}
	where $x_E^{n}$ is initialized by $x$, i.e. $x_E^{n}=x$, and the Lagrangian fluxes $\mathcal{H}_j^L$ are simply evaluated at each stage $i$ as
	\begin{equation}
		\mathcal{H}_i^L=-u(x_E^{(i)},t).
		\label{eqn.updateHL}
	\end{equation}
	An example of the trajectory drawn by a grid point for the second order three-step Runge Kutta scheme \eqref{eqn.IMEX2} is shown in Figure \ref{fig.1}.
	The foot of the characteristic is located at $x^L$, which corresponds to the solution of \eqref{eqn.TrODE}. According to \eqref{e4}, this is given by
	\begin{eqnarray}
		x^L=x_E^{n} + \Delta t \sum\limits_{i=1}^{s} b_i \mathcal{H}_i^L.
	\end{eqnarray}
	Finally, the solution at the next time level $t^{n+1}$ for the transported quantity, i.e. $q^L$, is computed by the reconstruction, or interpolation, of the numerical solution at the foot of the characteristic $x^L$, thus leading to
	\begin{equation}
		q^{n+1}(x)=\mathcal{R}\left(q\left(x^L,t^n\right)\right),
		\label{eqn.Rec2}
	\end{equation}
	where $\mathcal{R}(\cdot)$ is a suitable reconstruction operator.
	
	The order of convergence of the scheme directly comes from the order of the Runge-Kutta method given by the explicit tableau in \eqref{eqn.butcher}, which is used to solve the trajectory ODE. Then, the interpolation step \eqref{eqn.Rec2} must also be of the same order of accuracy of the RK scheme in order to provide a consistent semi-Lagrangian discretization of the advection dominated PDE. 
	
	The semi-Lagrangian operator \eqref{eqn.Lop} is then a compact definition which embeds all the equations \eqref{eqn.xE}-\eqref{eqn.Rec2}, that are nothing but the method corresponding to point 1 and 2 described at the beginning of Section \ref{sec.SL-IMEX}.

	\subsubsection{Diffusion dominated implicit RK schemes}
	For diffusion dominated problems the SL discretization does not play any role, thus the scheme must be equivalent to the implicit Runge-Kutta method with $\mathcal{H}=(0,\mathcal{H}^I)$ for $\mathbf{q}=(q,x)$:
	\begin{equation}
			\left\{\begin{aligned}
				\frac{d q}{d t} &= \mathcal{H}^I, \qquad  \mathcal{H}^I = \alpha \frac{\partial^2 q}{\partial x^2}, \\
				\frac{d x}{d t} &= 0. \\
			\end{aligned}\right. 
			\label{eqn.diffSubSyst}
	\end{equation}
	The implicit fluxes are taken into account by the term $\mathcal{H}^I$ and the solution is computed as follows:
	
	\begin{itemize}
		\item Stage values for $i=1 \ldots s$
		\begin{subequations}
			\begin{align}
				q_{*}^{(i)}&=q^n  + \Delta t \sum\limits_{j=1}^{i-1}a_{i,j} \, \mathcal{H}_{j}^I, \\
				q_I^{(i)}&=q_{*}^{(i)} +\Delta t \, a_{i,i} \, \mathcal{H}_i^I { \quad \mbox{where} \quad \mathcal{H}_i^I=\alpha \diff{^2 q_{I}^{(i)}}{x^2}}. \label{e4d_sub_2}
			\end{align}
			\label{e4d_sub}
		\end{subequations}
		\item Update of the solution according to \eqref{e4} as
		\begin{equation}
			q^{n+1}=q^n+ \Delta t \sum\limits_{i=1}^{s} b_i \mathcal{H}_i^I.	
			\label{e4d}
		\end{equation}
	\end{itemize}
	Notice that the determination of the stage value in \eqref{e4d_sub_2} implies the solution of an implicit system. The order of convergence is guaranteed in this case by the implicit part of the IMEX Runge-Kutta scheme \eqref{eqn.butcher}, and no difference with classical IMEX schemes with Eulerian advection arises. For the space discretization we can use, for example, a fourth order central finite difference scheme:
	\begin{equation}
		 \frac{\partial^2 q}{\partial x^2} =  \frac{-q_{i+2} + 16 q_{i+1} -30 q_i + 16 q_{i-1} - q_{i-2}}{12\Delta x^2}+ \mathcal{O}(\Delta x^4),
		\label{eqn.FD4}
	\end{equation}
   with $i$ and $\Delta x$ being the cell index and the mesh spacing, respectively.

	\subsubsection{Advection-diffusion SL-IMEX schemes}
	From the observations illustrated so far, we have seen that, for advection dominated problems, high order of accuracy in time can be achieved by solving the Lagrangian trajectory ODE \eqref{eqn.TrODE} and then by evaluating the solution at the foot of the characteristic relying on a high order spatial reconstruction operator \eqref{eqn.Rec2}. On the other hand, diffusion dominated problems allows classical IMEX discretizations \eqref{e1}-\eqref{e4} to be directly adopted, since a complete splitting between time and space fluxes takes place. In the sequel we will illustrate three different approaches combining the advective SL solver with the finite difference discretization.
	
	\paragraph{Algorithm 0.} To couple advection and diffusion schemes, one is now tempted by i) first moving the solution $q^n$ along the Lagrangian trajectories to obtain $q_E$, ii) and then solving the implicit equation \eqref{e2} for the diffusion terms to evaluate $q_I$. The algorithm in this case would work as follows:
	\begin{itemize}
		\item Stage values for $i=1 \ldots s$
		\begin{subequations}
			\begin{align}
				x_E^{(i)}&=x_E^{n}+\Delta t \sum\limits_{j=1}^{i-1}\tilde{a}_{i,j} \mathcal{H}_j^L \label{ad0} & \quad \\
				q_E^{(j+1)}&=\LL\left(q_E^{(j)},\tilde{a}_{i,j} \Delta t, \mathcal{H}^L_j \right) +\Delta t \, \tilde{a}_{i,j} \mathcal{H}^I_j &\quad j=1 \ldots i-1 \label{ad1} \\
				q_{*}^{(j+1)}&=\LL\left(q_{*}^{(j)},a_{i,j} \Delta t, \mathcal{H}^L_j \right) +\Delta t \, a_{i,j} \mathcal{H}^I_j &\quad j=1 \ldots i-1 \label{ad2}	
			\end{align}
		\end{subequations}
		with the initialization $q_E^{(1)}=q^n$ and $q_{*}^{(1)}=q^n$. The computation of the explicit Lagrangian flux for the trajectory equation \eqref{eqn.TrODE} is given by
		\begin{equation}
			\mathcal{H}^L_i=-u(x_E^{(i)}).
			\label{ad3}
		\end{equation}
		
		The evaluation of the implicit flux $\mathcal{H}^I_i$ for the diffusion terms requires knowledge of the state $q_I^{(i)}$ that is obtained by solving the following equation
		\begin{equation}
			q_I^{(i)}-a_{i,i} \Delta t \, \alpha \frac{\partial^2 q_I^{(i)} }{\partial x^2} = \LL\left(q_{*}^{(i)},a_{i,i}\Delta t, \mathcal{H}^L_i \right), \label{adsys1}\\
		\end{equation}
		which leads to the definition of the implicit fluxes
		\begin{equation}
			\mathcal{H}^I_i=\alpha \frac{\partial^2 q_I^{(i)} }{\partial x^2}.
			\label{eqn.HI}
		\end{equation}
		
		\item Update of the solution $q^{n+1}$ at the next time level using \eqref{e4}.
	\end{itemize}
	\noindent
	Unfortunately this algorithm achieves only first order of accuracy in time, even if high order IMEX schemes are adopted. Indeed, Algorithm 0 is equivalent to a first order splitting, where the governing PDE is split into explicit and implicit part: first, the explicit advection term is solved using \eqref{eqn.Lop}, then the implicit part is computed employing \eqref{e4d_sub}-\eqref{e4d}, and eventually the solution is obtained by summing both contributions. We emphasize that by using higher order splitting approaches the order of accuracy can be increased (see \cite{Strang68} for example). However, in the sequel, we will not explore this direction further even though, as we will see, Algorithm 2 has similar advantages to a splitting approach.	
	
		
	\paragraph{Algorithm 1.} In order to show the idea which yields the construction of high order SL-IMEX schemes, let us consider the simplest combination of advection and diffusion by taking a constant velocity field $u=const$ as well as a non zero diffusion coefficient $\alpha>0$. Let then the initial condition be a step function
	\begin{eqnarray}
		q(x,0)=\left\{
		\begin{array}{ll}
			q_L & x<0 \\
			q_R & x\geq 0
		\end{array}
		\right. .
		\label{ICad}
	\end{eqnarray}
	The exact solution of \eqref{sys1.1} with initial condition \eqref{ICad} reads
	\begin{equation}
		q(x,t)=\frac{1}{2}(q_L+q_R)+\frac{1}{2} \erf \left(\frac{1}{2}\frac{x-ut}{\sqrt{\alpha t}}\right)(q_R-q_L).
	\end{equation}
	Looking at equation \eqref{adsys1}, it is evident that the semi-Lagrangian contribution moves the solution according to the Lagrangian flux $\mathcal{H}^L$ that is derived from \eqref{ad0} and \eqref{ad3}, thus the convective terms are solved at the aid of the explicit Butcher tableau in \eqref{eqn.butcher}. On the other hand, the implicit flux is defined according to \eqref{eqn.HI} which follows from the solution of \eqref{adsys1}, that implies the use of the implicit RK scheme in \eqref{eqn.butcher}. Therefore, the implicit and explicit contributions in \eqref{adsys1} are not compatible, meaning that the Lagrangian part is moving while the implicit flux is defined at the wrong time location within the time step $\Delta t$. For the simple problem given by \eqref{sys1.1}-\eqref{ICad}, this corresponds advecting the solution up to a certain time fractional step of the explicit RK stage, and then applying the diffusion operator at a different time fractional step given by the implicit RK stage. In classical Eulerian IMEX schemes, the fluxes $\mathcal{H}$ are completely independent from time, thus they can be easily used to solve the convective contribution at any desired time fractional step, hence allowing for a perfect compatibility in the solution of the implicit step \eqref{e2}. Contrarily, semi-Lagrangian schemes implies a full coupling of space and time discretization because they directly solve the Lagrangian equation \eqref{sys1.4}.
	
	In order to overcome this problem we underline two aspects. First, if $u=0$ the Lagrangian contribution does nothing and the flux $\mathcal{H}^I$ is automatically correct since it is given by the solution of the pure diffusion equation with the implicit RK scheme \eqref{eqn.butcher}. Second, if the solution is transported with $u \neq 0$ and we assume a Lagrangian reference system which moves with the characteristic equation \eqref{eqn.TrODE}, then the problem becomes again a diffusion dominated one and we reduce to the previous case. In other words, if we apply a change of coordinates from $x$ to $x^\prime$ such as $x^\prime=x-ut$, then we do not see the Lagrangian contribution. Thus, the idea is to move also the spatial flux $\mathcal{H}^I$ using the Lagrangian trajectories, so that it can be shifted with a time fractional step $\omega$ at the same time level of the advection contribution in the RK stages and vice-versa. The way the fluxes are moved has to be consistent with the change of coordinate system and the algorithm reads as follows:
	\begin{itemize}
		\item Stage values for $i=1 \ldots s$
		\begin{subequations}
			\begin{align}
				\omega &= \sum \limits_{r=1}^j a_{j,r} - \sum \limits_{r=1}^{j-1} a_{i,r} & \quad j=1 \ldots i-1 \label{eqn.alg1.1} \\
				\tilde{q}	&= \LL\left(q_{*}^{(j)},\omega \Delta t, -u(x,t) \right) + \Delta t \, a_{i,j} \mathcal{H}^I_j &\quad j=1 \ldots i-1 \label{eqn.alg1.2} \\
				q_{*}^{(j+1)}	&= \LL\left(\tilde{q},(a_{i,j}-\omega)\Delta t, -u(x,t) \right) &\quad j=1 \ldots i-1 \label{eqn.alg1.3}
			\end{align}
		\end{subequations}
		with the initialization $q_{*}^{(1)}=q^n$. The implicit flux $\mathcal{H}^I_j$ is defined at the time fractional step $\omega$ given by \eqref{eqn.alg1.1}. Thus, the solution is firstly advected at the time level of the implicit flux, which can then be added to obtain the intermediate state $\tilde{q}$ with \eqref{eqn.alg1.2}. Finally, the state $q_{*}^{(j+1)}$ can be easily obtained by shifting the entire solution, which now accounts for both advection and diffusion at the same compatible time level, to the time level of the current RK stage, i.e. $a_{i,j}$, see \eqref{eqn.alg1.3}.  
		
		The computation of the implicit flux $\mathcal{H}^I_i$ for the diffusion terms at the time level $a_{i,i}$ implies solving
		\begin{equation}
			q_I^{(i)}-a_{i,i} \Delta t \, \alpha \frac{\partial^2 q_I^{(i)} }{\partial x^2} = \LL\left(q_{*}^{(i)},a_{i,i}\Delta t, -u(x,t) \right), \label{adsysAlg1}\\
		\end{equation}
		and then using $q_I^{(i)}$ to obtain $\mathcal{H}^I_i$ from \eqref{eqn.HI}.
		
		\item The update of the solution $q^{n+1}$ at the next time level makes use of the coefficients $b_i$ for $i=1 \ldots s$, therefore the entire solution must be shifted again to the correct time location in order to sum both explicit and implicit contributions at the same time level. Therefore, for $i=1 \ldots s$ with $q_{*}^{(1)}=q^n$ it holds
		\begin{subequations}
			\begin{align}
				\omega &= \sum \limits_{r=1}^i a_{i,r} - \sum \limits_{r=1}^{i-1} b_{r}  \label{eqn.alg1.s1} \\
				\tilde{q}	&= \LL\left(q_{*}^{(i)},\omega \Delta t, -u(x,t) \right) + \Delta t \, b_i \, \mathcal{H}^I_i \label{eqn.alg1.s2} \\
				q_{*}^{(i+1)}	&= \LL\left(\tilde{q},(b_i-\omega)\Delta t, -u(x,t) \right) \label{eqn.alg1.s3}
			\end{align}
		\end{subequations} 
		and the final solution is given by $q^{n+1}=q_{*}^{s+1}$.
	\end{itemize}	
	
\begin{figure}[!htbp]
		\centering
		\includegraphics[width=0.65\textwidth]{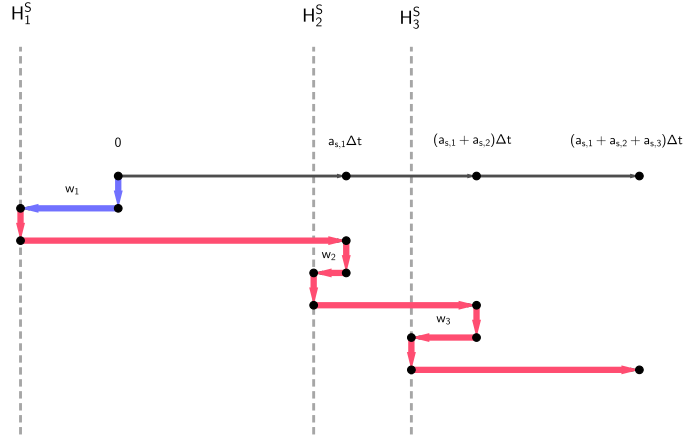}
		\caption{Illustration of the SL-IMEX scheme \eqref{eqn.IMEX2} for Algorithm 1. Blue lines represent the transport of the solution, i.e. the semi-Lagrangian scheme for the advection terms, while red lines are used for the transport of the coupled solution with fluxes.}  
		\label{fig.2}
	\end{figure}

	Algorithm 1 is based on the idea of going backward and forward with the entire solution $q$ in order to properly add the implicit contributions at compatible time levels within the Runge-Kutta framework. To guarantee that the formal order of convergence is achieved, this approach requires the crucial property of the closure of the trajectories: 
	\begin{equation}
		q =\LL(\LL(q,\omega \Delta t, u),-\omega \Delta t, u).
		\label{eqn.LagrClosure}
	\end{equation}
	The above equation implies that the same solution is recovered at the same spatial point after one back and forth round has been done, for a time interval $\omega \Delta t$ over a velocity field $u(x,t)$. This property is strictly exhibited only if the characteristic equation \eqref{eqn.TrODE} is solved exactly, thus providing an analytical expression for the trajectory. Apart from very simple velocity fields, obtaining an exact solution is practically not possible when dealing with nonlinear phenomena that typically occur in real world applications. Therefore, condition \eqref{eqn.LagrClosure} must be satisfied at the discrete level up to the order of accuracy of the method, so that the solution of \eqref{eqn.TrODE} does not affect the convergence of the overall IMEX scheme. This is why all semi-Lagrangian operators $\LL$ present in Algorithm 1 make use of the explicit RK scheme in \eqref{eqn.butcher} to solve the characteristic equation according to \eqref{eqn.xE}-\eqref{eqn.Rec2}. In this way, the numerical solution $q(x,t)$ can be shifted back and forth without spoiling the accuracy of the method, as experimentally proven by the convergence studies reported in Section \ref{sec.AdvDiff_test}.
	A graphical sketch of Algorithm 1 is depicted in Figure \ref{fig.2} for the SL-IMEX scheme \eqref{eqn.IMEX2}.
	
	\paragraph{Algorithm 2.} To improve the efficiency of Algorithm 1, we propose to directly shift the implicit fluxes to the final time level, where the new solution $q^{n+1}$ is assembled by means of the intermediate contributions at the RK stages. The main advantage is that Algorithm 2 does not require a local high order solution of the trajectory equation \eqref{eqn.TrODE}, which is instead mandatory in Algorithm 1. Algorithm 2 works as follows:
	
	\begin{itemize}
		\item Stage values for $i=1 \ldots s$
		\begin{subequations}
			\begin{align}
				x_E^{(j+1)} &= x_E^{(j)} + \tilde{a}_{i,j} \, \Delta t \, \mathcal{H}^L_j  & \quad j=1 \ldots i-1 \label{eqn.alg2.1} \\
				x_{*}^{(j+1)} &= x_{*}^{(j)} + a_{i,j} \, \Delta t \, \mathcal{H}^L_j  & \quad j=1 \ldots i-1 \label{eqn.alg2.2} \\
				\omega &= \sum \limits_{r=1}^j a_{j,r} - \sum \limits_{r=1}^{i-1} a_{i,r} & \quad j=1 \ldots i-1 \label{eqn.alg2.3} \\
				\tilde{\mathcal{H}^I}_{j}	&= \LL\left(\mathcal{H}^I_j,-\omega\Delta t, \mathcal{H}^L_j \right) &\quad j=1 \ldots i-1 \label{eqn.alg2.4}
			\end{align}
		\end{subequations}
		with the initialization $x_E^{(1)}=x_{*}^{(1)}=x$. The quantities $\tilde{\mathcal{H}^I}_j$ denote the transported implicit fluxes at each time fractional step $\omega$ given by \eqref{eqn.alg2.3}, while the Lagrangian flux $\mathcal{H}^L_i$ is updated with \eqref{eqn.updateHL} using $x_E^{(i)}$ coming from \eqref{eqn.alg2.1}. The intermediate state $q_{*}^{(i)}$ is then computed as
		\begin{equation}
			q_{*}^{(i)} = \mathcal{R}\left(q^n\left(x_{*}^{(i)}\right)\right) + \Delta t \sum \limits_{j=1}^{i-1} a_{i,j} \, \tilde{\mathcal{H}^I}_j,
			\label{eqn.qs_alg2}
		\end{equation}
	    where we recall that $\mathcal{R}$ represents a suitable high order spatial reconstruction operator. As for Algorithm 1, the implicit contribution $\mathcal{H}^I_i$ is obtained from \eqref{eqn.HI} where $q_I^{(i)}$ is the solution of the implicit equation
		\begin{equation}
			q_I^{(i)}-a_{i,i} \, \Delta t \, \alpha \frac{\partial^2 q_I^{(i)} }{\partial x^2} = \LL\left(q_{*}^{(i)},a_{i,i}\Delta t, \mathcal{H}^L_i \right), \label{adsysAlg2}\\
		\end{equation}
		in which the semi-Lagrangian operator $\LL$ makes use of a fast first order solver for the trajectory equation, i.e. the foot of the characteristic results to be $x^L=x-a_{i,i}\Delta t \, \mathcal{H}^L_i$.
		
		\item The final solution is updated by adopting the same strategy, hence shifting the implicit fluxes to the time fractional steps given by the coefficients $b_i$ in \eqref{eqn.butcher}. For the stages $i=1 \ldots s$ one computes
		\begin{subequations}
			\begin{align}
				x_E^{(i+1)} &= x_E^{(i)} + b_i \, \Delta t \, \mathcal{H}^L_i  &  \\
				\omega   &= \sum\limits_{r=1}^i a_{i,r}-\sum\limits_{r=1}^{s} b_r \\
				\tilde{\mathcal{H}^I}_{i} &= \LL\left(\mathcal{H}^I_i,-\omega \Delta t, \mathcal{H}^L_i \right)
			\end{align}
		\end{subequations}
		then the solution at the next time level is obtained by combining the semi-Lagrangian scheme for the advection part and the shifted implicit flux contributions, thus
		\begin{equation}	
			q^{n+1} = \mathcal{R}\left(q^n\left(x_{E}^{(s+1)}\right)\right) + \Delta t  \sum\limits_{j=1}^{s}b_j \tilde{\mathcal{H}^I}_{j}.
		\end{equation}
		
	\end{itemize}

\begin{figure}[!htbp]
		\centering
		\includegraphics[width=0.65\textwidth]{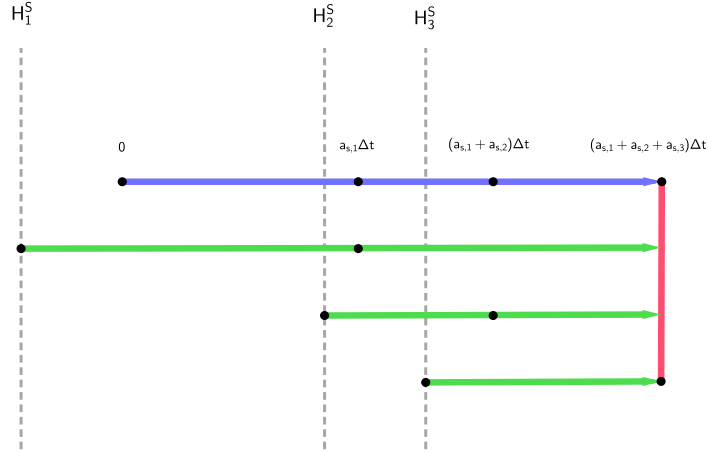}
		\caption{Illustration of the SL-IMEX scheme \eqref{eqn.IMEX2} for Algorithm 2. Blue lines represent the transport of the solution, i.e. the semi-Lagrangian scheme for the advection terms. Green lines refer to the transport of the implicit flux only. Red lines are used for the transport of the coupled solution with fluxes.}   
		\label{fig.3}
	\end{figure}

	Figure \ref{fig.3} shows a schematic of Algorithm 2 for the second order SL-IMEX scheme \eqref{eqn.IMEX2}.
	One advantage, is that this version of the SL-IMEX scheme can be easily embedded in an already existing semi-Lagrangian code for the explicit part. Indeed, Algorithm 2 only requires the transport of the implicit fluxes according to \eqref{eqn.alg2.4}, which have then to be added to the 
	explicit convection contribution in \eqref{eqn.qs_alg2}.
	
	\section{Numerical results for the advection-diffusion equation} \label{sec.AdvDiff_test}
\begin{figure}[!htbp]
		\begin{center}
			\begin{tabular}{ccc} 
				\includegraphics[width=0.33\textwidth]{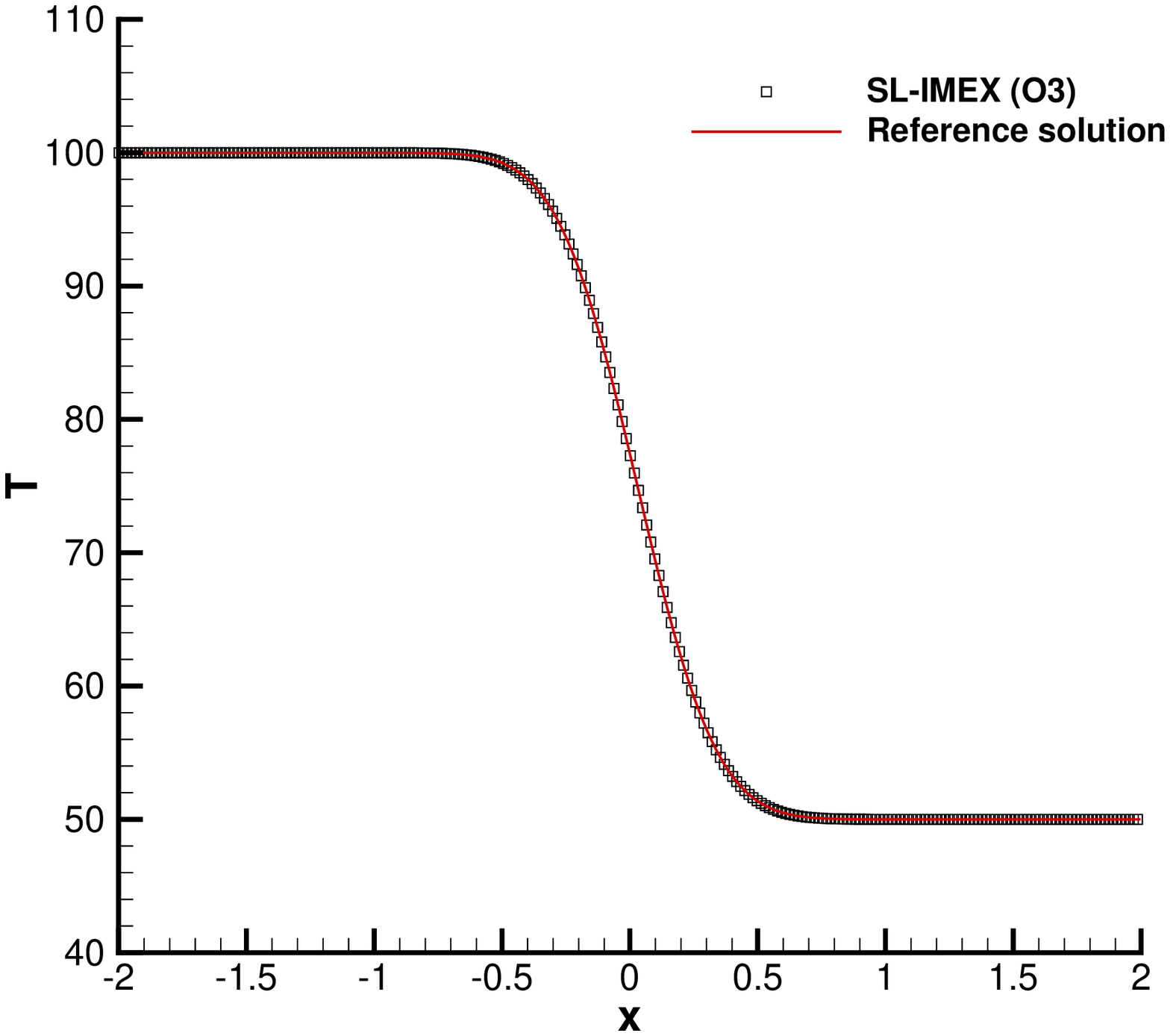}  &                 
				\includegraphics[width=0.33\textwidth]{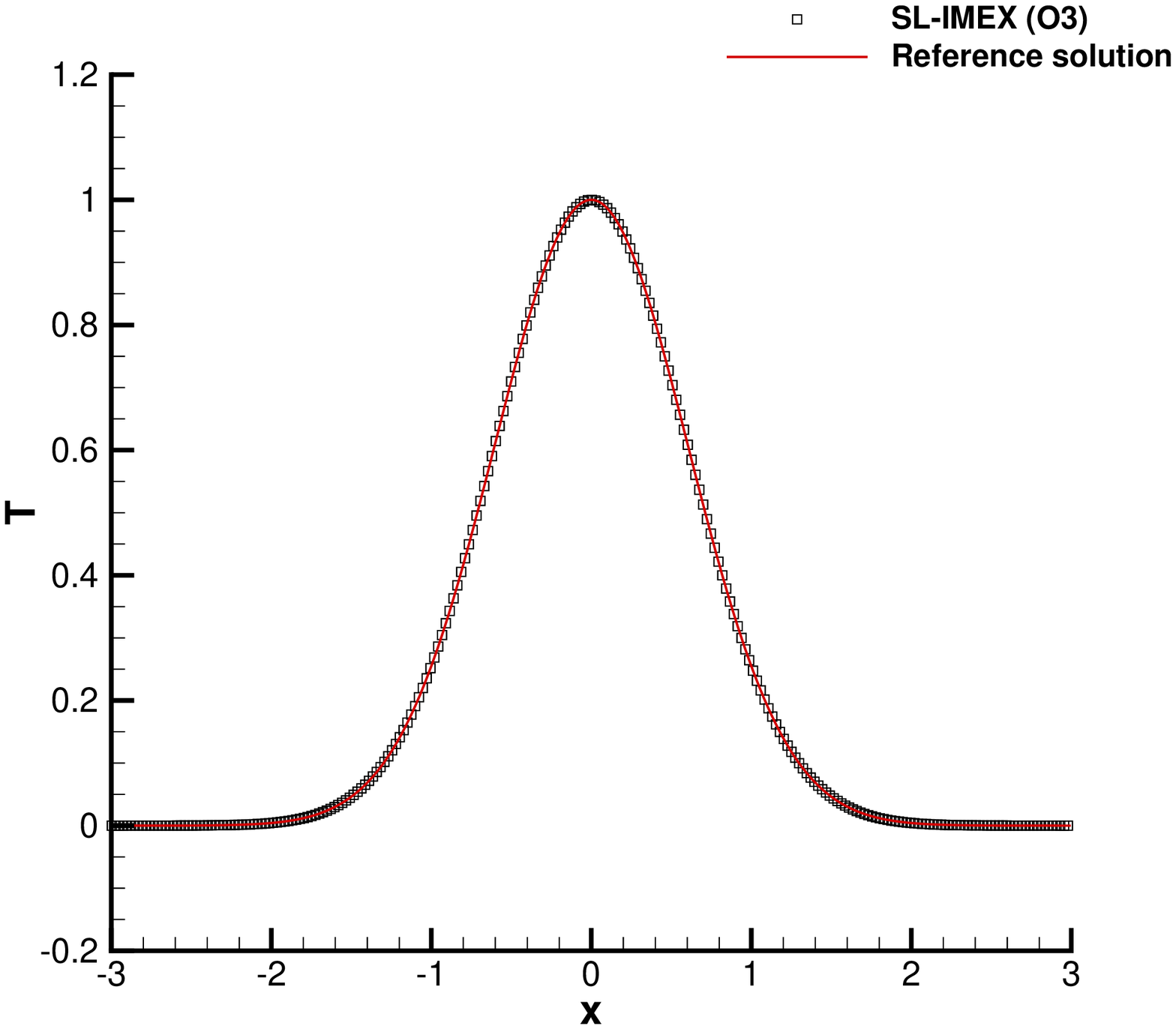}  &                 
				\includegraphics[width=0.33\textwidth]{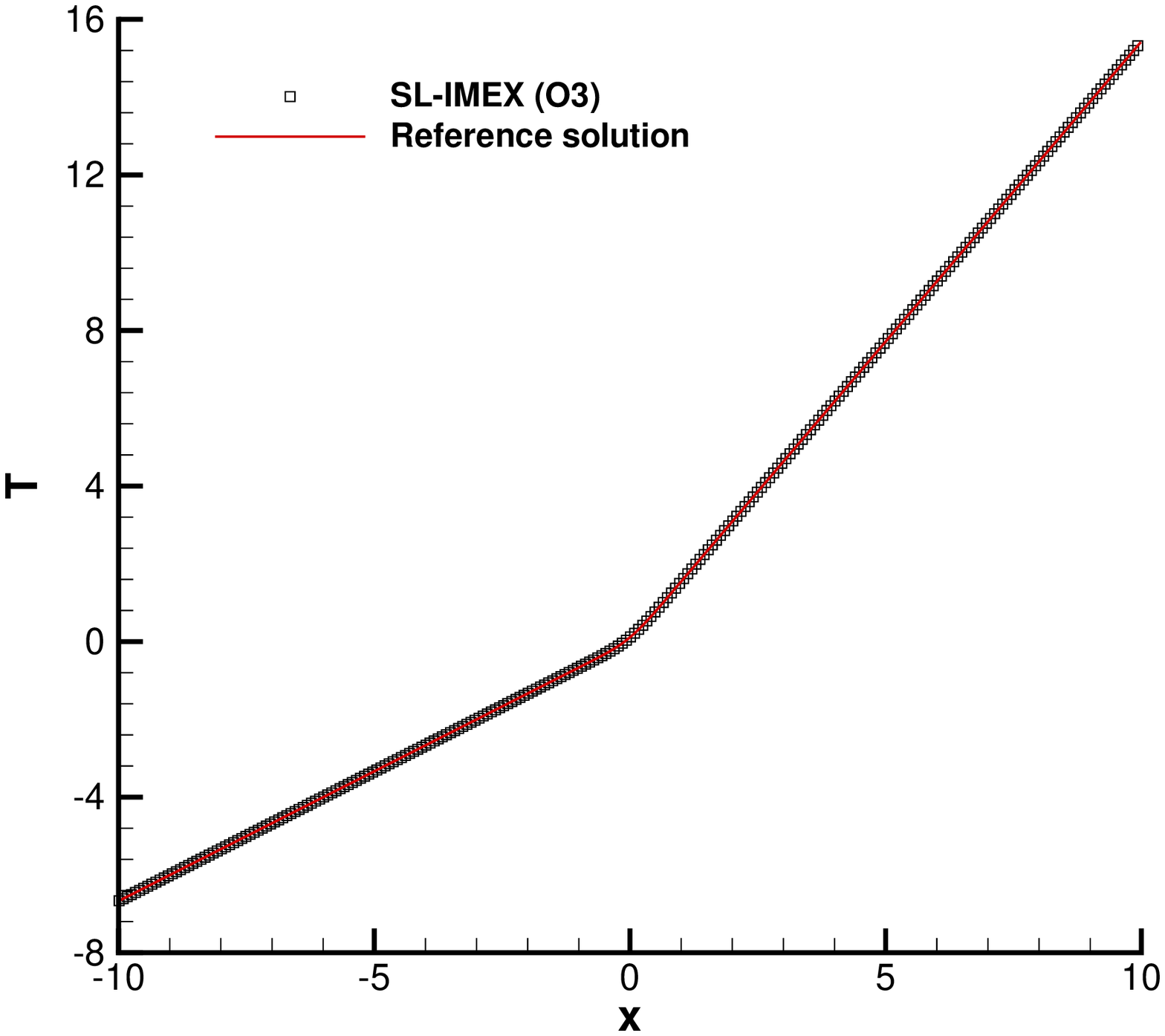} \\
			\end{tabular}
			\caption{Third order numerical results obtained using the SL-IMEX scheme with Algorithm 2 for the test cases involving the advection-diffusion equation: advection-diffusion equation with linear transport at $t=0.3$ (left), advection equation with space-dependent transport at $t=9$ (middle) and advection-diffusion equation with space-dependent transport at $t=0.1$ (right).}
			\label{fig.AdvDiff-test}
		\end{center}
	\end{figure}

	In the following, the new semi-Lagrangian IMEX methods (SL-IMEX) are applied to a set of test cases for the advection-diffusion equation \eqref{sys1.1}. All simulations are run with both Algorithm 1 and Algorithm 2 using as a reconstruction operator $\mathcal{R}$ a simple cubic spline interpolation, which guarantees fourth order of accuracy in space. We study the numerical convergence in time by firstly considering the transport and the diffusion part of the equation separately, then by proposing a non-trivial solution of the PDE that involves advection as well as diffusive terms simultaneously. Specifically, Test 1 is concerned with the solution of \eqref{sys1.1} with linear transport, thus the advective velocity is maintained constant in space and time. Test 2 deals with only the transport part of the PDE with an advection speed that is space dependent, neglecting the diffusive contribution. Finally, the convective terms as well as the diffusion part of the equation are solved in Test 3 for which an analytical solution is derived. Figure \ref{fig.AdvDiff-test} depicts the numerical solution obtained with a third order SL-IMEX scheme for all the three test problems, and a comparison against the reference solution is shown at the final time of each simulation.
	
	The computational domain is defined by $\Omega \in [x_L;x_R]$ and is discretized with a total number of cells $N_x=1000$. In this way, the numerical error related to the spatial discretization is ensured to be smaller than the time error, which is needed in order to properly study the convergence of the IMEX time stepping. 
	In a similar way, $N_t$ will be used to denote the number of time steps so that $\Delta t=(t_{f}-t_{0})/N_t$.
	The spatial discretization is of fourth order with a finite difference scheme for the implicit part and a cubic spline reconstruction for the convective terms, while the time accuracy can achieve up to third order. Therefore, for all test cases, we always guarantee the inequality $(\Delta x)^4<(\Delta t)^3$ to be satisfied and we report both the mesh spacing and the time step used for each simulation.
	The stability of an explicit scheme for the numerical solution of  \eqref{sys1.1} can be studied relying on a Von Neumann analysis, which shows that the time step $\Delta t$ must be chosen respecting the CFL condition for the advection part ($\Delta t^h$), and the parabolic restriction for the diffusion part ($\Delta t^p$), thus yielding
	\begin{equation}
		\Delta t = \min(\Delta t^h,\Delta t^p) \qquad \textnormal{with} \qquad \frac{\Delta t^h}{\Delta x} \max|u| \leq \textnormal{CFL}, \quad \frac{\Delta t^p}{\Delta x^2} \alpha \leq \kappa,
		\label{eqn.timestepAdvDiff} 
	\end{equation}
	with $\textnormal{CFL}\leq 1$ and $\kappa \leq 1/2$. On the contrary, the novel semi-Lagrangian IMEX methods are free from both stability conditions \eqref{eqn.timestepAdvDiff}, and in all the test cases illustrated hereafter we choose the time step as $\Delta t=\Delta t^p$ with $\kappa=400$. Because the parabolic stability condition is typically much more severe than the hyperbolic counterpart, the novel class of SL-IMEX schemes can solve the advection-diffusion equation with a time step that is about 800 times bigger than the time step corresponding to a fully explicit scheme like the well-known finite difference FTCS (Forward in Time Centered in Space) method.

	\subsection{Test 1: advection-diffusion equation with linear transport}
	The semi-Lagrangian discretization of the advection terms implies the coupling between time and space discretization, thus the standard IMEX schemes can not be directly applied. More precisely, the semi-Lagrangian version of IMEX schemes needs to properly transport the implicit fluxes according to the velocity field in order to be compatible with the Runge-Kutta stages of the IMEX scheme in time. To check this consistency we perform a simple but not trivial test concerning advection and diffusion. We first set a constant velocity field $u(x,t)=0.1$ and $\alpha=10^{-3}$ on the computational domain $\Omega \in [-2;2]$. The initial condition is prescribed as a step function centered in $x_0=-0.25$, so that the exact solution is given by the analytical solution of the stationary heat equation which is furthermore transported at a constant advection speed. This explicitly writes
	\begin{equation}
		q(x,t)=\frac{1}{2}-\frac{1}{2} \erf \left( \frac{x-x_0-ut}{\sqrt{4\alpha t}} \right).
	\end{equation}
	In order to avoid oscillations due to the discontinuity at $t=0$, the initial condition is taken at time $t_0$, namely
	\begin{equation}
		q(x,0)=\frac{1}{2}-\frac{1}{2} \erf \left( \frac{x-x_0-ut}{\sqrt{4\alpha (t +t_0)}} \right),
	\end{equation}
	with $t_0=10^{-2}$. SL-IMEX schemes with different accuracy in time are used to compute the solution at the final time $t_f=0.3$ with both Algorithm 1 and Algorithm 2. Table \ref{tab:test1} reports the $L_2$ error norms and the convergence studies for $p=[0,1,2]$, demonstrating that the formal order accuracy is achieved by the novel SL-IMEX methods.
	
	\begin{table}[!h]
		\begin{center}
			\renewcommand{\arraystretch}{1.2}
			\begin{tabular}{c|cccccc}
				\multicolumn{7}{c}{Algorithm 1} \\
				\hline
				\hline
				\multirow{2}{*}{$N_t$}    & \multicolumn{2}{c}{$\RR=0$} & \multicolumn{2}{c}{$\RR=1$} & \multicolumn{2}{c}{$\RR=2$} \\ 
				& $L_2$ & $\mathcal{O}(q)$ & $L_2$ & $\mathcal{O}(q)$ & $L_2$ & $\mathcal{O}(q)$ \\ \hline
				1 &  1.7026E-00 & -    & 1.0879E-00 &  -    & 6.3812E-01 & - \\ 
				2 &  9.4437E-01 & 0.85 & 1.4706E-01 &  2.89 & 5.7912E-02 & 3.46 \\
				4 &  4.9826E-01 & 0.92 & 2.2369E-02 &  2.71 & 6.5048E-03 & 3.15 \\
				8 &  2.5592E-01 & 0.96 & 5.4022E-03 &  2.04 & 9.4678E-04 & 2.78 \\ 
				16 &  1.2967E-01 & 0.98 & 1.3312E-03 &  2.02 & 1.2953E-04 & 2.87 \\ 
			\end{tabular}
		\end{center}
		
		\begin{center}
			\renewcommand{\arraystretch}{1.2}
			\begin{tabular}{c|cccccc}
				\multicolumn{7}{c}{Algorithm 2} \\
				\hline
				\hline
				\multirow{2}{*}{$N_t$}    & \multicolumn{2}{c}{$\RR=0$} & \multicolumn{2}{c}{$\RR=1$} & \multicolumn{2}{c}{$\RR=2$} \\ 
				& $L_2$ & $\mathcal{O}(q)$ & $L_2$ & $\mathcal{O}(q)$ & $L_2$ & $\mathcal{O}(q)$ \\ \hline
				1 &  1.7026E-00 & -    & 1.0879E-00 &  -    & 6.3812E-01 & - \\ 
				2 &  9.4438E-01 & 0.85 & 1.4706E-01 &  2.89 & 5.7912E-02 & 3.46 \\
				4 &  4.9826E-01 & 0.92 & 2.2369E-02 &  2.71 & 6.5048E-03 & 3.15 \\
				8 &  2.5592E-01 & 0.96 & 5.4022E-03 &  2.04 & 9.4666E-04 & 2.78 \\ 
				16 &  1.2967E-01 & 0.98 & 1.3312E-03 &  2.02 & 1.2943E-04 & 2.87 \\ 
			\end{tabular}
		\end{center}	
		\caption{Convergence study for Test 1 involving the advection-diffusion equation with linear transport.}\label{tab:test1}
	\end{table}
	In order to highlight the importance of the time consistency between explicit and implicit contributions, the same test is run with {\itshape Algorithm 0} that does not operate any transport of the fluxes. The resulting errors and convergence rates are reported in Table \ref{tab:test_2} where it is clear that only first order is achieved. Notice that since the velocity is constant, the Lagrangian trajectory ODE \eqref{eqn.TrODE} is solved exactly and the pure diffusion term is of fourth order in space. Therefore the error shown in Table \ref{tab:test_2} is essentially due to the wrong coupling of the transported advection term and the diffusion contribution, which is extremely important in the construction of high order SL-IMEX schemes.
	
	\begin{table}[!h]
		\begin{center}
			\renewcommand{\arraystretch}{1.2}
			\begin{tabular}{c|cccccc}
				\multicolumn{7}{c}{Algorithm 0} \\
				\hline
				\hline
				\multirow{2}{*}{$N_t$}    & \multicolumn{2}{c}{$\RR=0$} & \multicolumn{2}{c}{$\RR=1$} & \multicolumn{2}{c}{$\RR=2$} \\ 
				& $L_2$ & $\mathcal{O}(q)$ & $L_2$ & $\mathcal{O}(q)$ & $L_2$ & $\mathcal{O}(q)$ \\ \hline
				1 &  1.7026E-00 & -    & 2.2826e+00 &  -    & 2.1469e+00 & - \\ 
				2 &  9.4437E-01 & 0.85 & 3.9545e-01 &  2.53 & 3.4371e-01 & 2.64 \\
				4 &  4.9826E-01 & 0.92 & 1.0053e-01 &  1.98 & 9.4635e-02 & 1.86 \\
				8 &  2.5592E-01 & 0.96 & 4.4798e-02 &  1.17 & 4.3513e-02 & 1.12 \\ 
				16 &  1.2967E-01 & 0.98 & 2.1298e-02 &  1.07 & 2.1001e-02 & 1.05 \\ 
			\end{tabular}
		\end{center}	
		\caption{Convergence study for Test 1 involving the advection-diffusion equation with linear transport and Algorithm 0.}\label{tab:test_2}
	\end{table}
	
	\subsection{Test 2: advection equation with space-dependent transport}
	The aim of this test is to verify the high order time accuracy by solving only the advection part of the PDE, while neglecting the diffusion terms. A non-constant velocity field is considered, so that the high order numerical solution of the ODE for the flow trajectories plays a crucial role for obtaining the formal order of accuracy. The computational domain is defined by $\Omega=[-3;3]$ and we set $u(x,t)=u(x,0)=k_0 x$ with $k_0=0.2$. The initial condition for the scalar quantity $q$ reads
	\begin{equation}
		q(x,0)= e^{-50 \, x^2},
	\end{equation}
	and the simulation is run until the final time $t_f=9$. For this test an exact solution can be found by following the characteristic equation, see \cite{TB19}:
	\begin{equation}
		q(x,t)= e^{-50 \, (x \, e^{-k_0 \, t})^2}.
	\end{equation}
	The resulting convergence rates are reported in Table \ref{tab:test2} for both versions of the SL-IMEX algorithms, showing that the order of accuracy is correctly reproduced by all schemes with $p=[0,1,2]$. 
	
	\begin{table}[!h]
		\begin{center}
			\renewcommand{\arraystretch}{1.2}
			\begin{tabular}{c|cccccc}
				\multicolumn{7}{c}{Algorithm 1} \\
				\hline
				\hline
				\multirow{2}{*}{$N_t$}    & \multicolumn{2}{c}{$\RR=0$} & \multicolumn{2}{c}{$\RR=1$} & \multicolumn{2}{c}{$\RR=2$} \\ 
				& $L_2$ & $\mathcal{O}(q)$ & $L_2$ & $\mathcal{O}(q)$ & $L_2$ & $\mathcal{O}(q)$ \\ \hline
				1 &  1.6597E-00 & -    & 7.2176E-02 &  -    & 2.9958E-01 & - \\ 
				2 &  8.2470E-01 & 0.47 & 1.5950E-02 &  2.18 & 2.1377E-02 & 3.81 \\
				4 &  5.9345E-01 & 1.48 & 3.6948E-03 &  2.11 & 1.8254E-03 & 3.54 \\
				8 &  2.2644E-01 & 1.39 & 8.8701E-04 &  2.06 & 1.9175E-04 & 3.25 \\ 
				16 &  1.0085E-01 & 1.17 & 2.1698E-04 &  2.03 & 2.4340E-05 & 2.97 \\ 
			\end{tabular}
		\end{center}
		
		\begin{center}
			\renewcommand{\arraystretch}{1.2}
			\begin{tabular}{c|cccccc}
				\multicolumn{7}{c}{Algorithm 2} \\
				\hline
				\hline
				\multirow{2}{*}{$N_t$}    & \multicolumn{2}{c}{$\RR=0$} & \multicolumn{2}{c}{$\RR=1$} & \multicolumn{2}{c}{$\RR=2$} \\ 
				& $L_2$ & $\mathcal{O}(q)$ & $L_2$ & $\mathcal{O}(q)$ & $L_2$ & $\mathcal{O}(q)$ \\ \hline
				1 &  1.6597E-00 & -    & 5.0765E-01 &  -    & 6.4530E-01 & - \\ 
				2 &  8.2470E-01 & 0.47 & 1.5163E-01 &  1.74 & 6.1088E-02 & 3.40 \\
				4 &  5.9345E-01 & 1.48 & 3.3432E-02 &  2.18 & 4.3401E-03 & 3.82 \\
				8 &  2.2644E-01 & 1.39 & 7.5901E-03 &  2.13 & 3.9663E-04 & 3.45 \\ 
				16 &  1.0085E-01 & 1.17 & 1.7993E-03 &  2.07 & 4.3374E-05 & 3.19 \\ 
			\end{tabular}
		\end{center}	
		\caption{Convergence study for Test 2 involving the advection equation with space-dependent transport.}\label{tab:test2}
	\end{table}

	\subsection{Test 3: advection-diffusion equation with space-dependent transport}
	Since the transport of the velocity field and the fluxes is an important ingredient in order to achieve the formal order of convergence in time, we propose here a test involving advection and diffusion for a non-constant velocity field. Let $u(x,t)=u(x,0)=-x$ and $\alpha>0$. 
	We are looking for an exact solution of the advection-diffusion equation \eqref{sys1.1} written in non-conservative form: 
	\begin{equation}
		\diff{q}{t}  -x\diff{q}{x} = \alpha \frac{\partial^2 q}{\partial x^2}.
		\label{sysNCP}
	\end{equation}
	Let us assume that the solution can be expressed as the product of two independent functions, namely $q(x,t)= G(x)\cdot H(t)$, where $G=G(x)$ and $H=H(t)$ are functions of the space and time only, respectively. Then, $q$ is a solution of \eqref{sys1.1} if it holds
	\begin{equation}
		\diff{H}{t} G  -x H \diff{G}{x} = H \alpha \frac{\partial^2 G}{\partial x^2}.
		\label{sysNCP2}
	\end{equation}
	This requires to find $H$ and $G$ so that
	\begin{equation}
		\left\{ \begin{aligned}
			\frac{dH(t)}{dt}&=H(t) \, c_1 \\
			\alpha \frac{d^2G(x)}{dx^2}&=G(x) \, c_1 - x \frac{dG(x)}{dx}
		\end{aligned} \right. ,
		\label{sysNCP3}
	\end{equation}
	for a generic constant $c_1$. By setting $c_1=1$, a general solution of the previous system writes
	\begin{eqnarray}
		H(t)&=& c_0 e^t, \nonumber \\
		G(x)&=& c_2 x+c_3\left(\sqrt{2 \pi \alpha} \erf\left( \frac{1}{\sqrt{2}x}\right)x -e^{-\frac{x^2}{2\alpha}}\alpha \right),
	\end{eqnarray}
	from which we extract a particular solution by setting $c_0=1$, $c_2=-1$, $c_3=-1/2$, hence obtaining
	\begin{equation}
		q_{ex}(x,t)= e^{t} \cdot \left( -x-\frac{1}{2}\sqrt{2 \pi \alpha} \erf\left( \frac{1}{\sqrt{2}x}\right)x -e^{-\frac{x^2}{2\alpha}}\alpha \right).
	\end{equation}
	We perform the test by setting the initial condition as $q(x,0)=q_{ex}(x,0)$ on a domain $\Omega=[-10;10]$. The diffusion coefficient is $\alpha=0.1$ and the final time is set to $t_f=0.1$. Convergence studies are analyzed and shown in Table \ref{tab:test3}, which confirm the correct behavior of both Algorithm 1 and 2 up to third order in time.
	
	\begin{table}[!h]
		\begin{center}
			\renewcommand{\arraystretch}{1.2}
			\begin{tabular}{c|cccccc}
				\multicolumn{7}{c}{Algorithm 1} \\
				\hline
				\hline
				\multirow{2}{*}{$N_t$}    & \multicolumn{2}{c}{$\RR=0$} & \multicolumn{2}{c}{$\RR=1$} & \multicolumn{2}{c}{$\RR=2$} \\ 
				& $L_2$ & $\mathcal{O}(q)$ & $L_2$ & $\mathcal{O}(q)$ & $L_2$ & $\mathcal{O}(q)$ \\ \hline
				1 &  1.0275E-01 & -    & 2.2454E-04 &  -    & 6.6623E-05 & - \\ 
				2 &  5.3078E-02 & 0.95 & 5.6660E-05 &  1.99 & 8.4538E-06 & 2.99 \\
				4 &  2.6987E-02 & 0.98 & 1.4231E-05 &  1.99 & 1.0654E-08 & 2.99 \\
				8 &  1.3609E-02 & 0.99 & 3.5660E-06 &  2.00 & 1.3573E-07 & 2.97 \\ 
				16 &  6.8338E-03 & 0.99 & 8.9253E-07 &  2.00 & 2.2647E-08 & 2.59 \\ 
			\end{tabular}
		\end{center}
		
		\begin{center}
			\renewcommand{\arraystretch}{1.2}
			\begin{tabular}{c|cccccc}
				\multicolumn{7}{c}{Algorithm 2} \\
				\hline
				\hline
				\multirow{2}{*}{$N_t$}    & \multicolumn{2}{c}{$\RR=0$} & \multicolumn{2}{c}{$\RR=1$} & \multicolumn{2}{c}{$\RR=2$} \\ 
				& $L_2$ & $\mathcal{O}(q)$ & $L_2$ & $\mathcal{O}(q)$ & $L_2$ & $\mathcal{O}(q)$ \\ \hline
				1 &  1.0276E-01 & -    & 1.7405E-03 &  -    & 8.5036E-05 & - \\ 
				2 &  5.3078E-02 & 0.95 & 4.4623E-04 &  1.96 & 1.1189E-05 & 2.93 \\
				4 &  2.6987E-02 & 0.98 & 1.1297E-04 &  1.98 & 1.4854E-06 & 2.91 \\
				8 &  1.3609E-02 & 0.99 & 2.8419E-05 &  1.99 & 2.1375E-07 & 2.80 \\ 
				16 &  6.8338E-03 & 0.99 & 7.1270E-06 &  1.99 & 3.7252E-08 & 2.53 \\ 
			\end{tabular}
		\end{center}	
		\caption{Convergence study for Test 3 involving the advection-diffusion equation with space-dependent transport.}\label{tab:test3}
	\end{table}
	
	\section{Conservative semi-Lagrangian IMEX schemes for the shallow water system} \label{sec.SWE}
	Semi-Lagrangian methods are typically designed for solving the non-conservative form of the advection equation, thus they cannot be directly applied to the solution of hyperbolic systems of conservation laws involving shock waves. In order to extend the SL-IMEX schemes to such systems, with the aim to properly treat the different space-time scales, we propose in the sequel a conservative formulation of the methods. As a prototype problem we will consider the one-dimensional shallow water equations (SWE):
	\begin{eqnarray}
		&&\diff{h}{t}+\diff{V}{x}=0, \label{sw1a}\\
		&&\diff{V}{t}+\diff{uV}{x}+\frac{g}{2}\diff{h^2}{x}=0,
		\label{sw1b}
	\end{eqnarray}
	where $u(x,t)$ is the fluid velocity and $V(x,t)=hu$ denotes the momentum. Let $h(x,t)=b(x)+\eta(x,t)$ be the total water depth that is computed as the sum of the prescribed bottom topography $b(x)$ and the location of the free surface $\eta(x,t)$. For the sake of simplicity, a constant flat bottom is assumed, thus $b(x)=0$ and therefore $h(x,t)=\eta(x,t)$. Finally, $g=9.81$ $m/s^2$ is the constant
	gravity acceleration and the eigenvalues of system \eqref{sw1a}-\eqref{sw1b} are $\lambda_{1,2}=u\pm \sqrt{g h}$. The term containing the gradient of the hydraulic head, i.e. $g/2 \, \partial_x h^2$, might be responsible for a stiffness in the governing equations, especially when the flow velocity is rather low, thus approaching the low Froude regime. This makes the adoption of an implicit-explicit discretization particularly interesting, in order to separate the slow advection scale from the fast scale related to the pressure. 
	
	\subsection{Time discretization}
	A clever discretization based on the IMEX strategy would consider the pressure terms implicitly, which are responsible of the celerity $\sqrt{g h}$ in the eigenvalues, while keeping an explicit scheme for the solution of the nonlinear convective contribution, which is related to the quantity $u$ in the eigenvalues. This choice would improve the efficiency of the scheme in low Froude number flows, since a milder time step condition based only on the fluid velocity would be enough to ensure numerical stability. In the following we propose two different semi-discrete schemes, which will be referred to as SL-IMEX-H and SL-IMEX.
	
	\paragraph{SL-IMEX-H method.} The semi-discrete SL-IMEX-H scheme for the one-dimensional SWE writes
	\begin{eqnarray}
		h^{n+1}&=&h^{n} - \Delta t \diff{V^{n+1}}{x}, \label{sw2a}\\
		V^{n+1}&=&V^{n} - \frac{\partial }{\partial x}\left(\int\limits_{t^n}^{t^{n+1}}u V \, dt \right)-\Delta t \frac{g}{2}\diff{{h}^{n} h^{n+1}}{x},
		\label{sw2b}
	\end{eqnarray}
	where the convective terms retain the time integral notation which will be useful in the sequel. The advection contribution in \eqref{sw2b} is discretized with an explicit conservative SL method, while the flux term in the continuity equation \eqref{sw2a} is approximated by means of an implicit finite difference scheme, thus we adress this method as SL-IMEX-Hybrid scheme. The pressure gradient is discretized relying on a semi-implicit approach, hence obtaining $\partial_x (h^n h^{n+1})$, from which a linear system on the total water depth naturally arises. This strategy has been recently pursued in \cite{BP2021} for the compressible Navier-Stokes equations and makes use of the autonomous form \eqref{eqn.H} for the application of IMEX RK schemes. Along the lines of \cite{ShuOsher1988,Shu_ConsSL2011}, a conservative scheme can then be derived by introducing the following definitions:
	\begin{eqnarray}
		H(x)&:=&\int\limits_{t^n}^{t^{n+1}}u V \, dt, \label{eqn.Hdef} \\
		F(x)&:=&V^n - \diff{H(x)}{x}, \label{eqn.Fdef}
	\end{eqnarray}
	where $H(x)$ is the time integral of the momentum flux and $F$ is an operator that contains the discretization of the explicit convective terms. Formal substitution of the momentum $V^{n+1}$ given by \eqref{sw2b} into the continuity equation \eqref{sw2a} leads to 
	\begin{equation}
		h^{n+1}- \Delta t^2 \, \frac{g}{2} \diff{^2 {h}^{n} h^{n+1}}{x^2}=h^n+\Delta t \diff{F}{x},
		\label{sw4}
	\end{equation}
	where the definition \eqref{eqn.Fdef} has been used and the only unknown is the total water depth at the new time level, i.e. $h^{n+1}$. Equation \eqref{sw4} represents a linear system for the unknown $h^{n+1}$ that can be solved using an efficient matrix-free implementation of the conjugate gradient method. Once the solution is obtained, the momentum can readily be updated according to \eqref{sw2b}.

	\paragraph{Remark on shock discontinuities.} For further stabilization at strong shock waves, an additional numerical flux should be embedded into the continuity equation \eqref{sw2a} in order to introduce some numerical diffusion \cite{TD_SWE2014}. Specifically, the semi-discrete continuity equation takes the form
	\begin{equation}
		h^{n+1}=h^{n} - \Delta t \left[ \diff{V^{n+1}}{x} - \frac{\lambda_{\max}(x)}{2}\diff{\delta^{n+1}}{x} \right],  \qquad \delta^{n+1}=h^{n+1,+} - h^{n+1,-},
		\label{eqn.visc}
	\end{equation}
	with $\lambda_{\max}(x)$ being the maximum eigenvalue of the shallow water system and $\delta^{n+1}$ representing a measure of the jump of the solution across cell interfaces in the spatial discretization between right and left state, i.e. $h^{n+1,+} $ and $h^{n+1,-}$, respectively. This corresponds to the use of a local Lax-Friedrichs flux in the mass conservation equation. Moreover, the introduction of implicit numerical viscosity in the form \eqref{eqn.visc} does not affect neither the stability of the scheme nor the symmetry of the linear system \eqref{sw4}.
	
	\paragraph{SL-IMEX method.} Alternatively, the shallow water system can be discretized in time as 
		\begin{eqnarray}
			h^{n+1}&=&G(x), \qquad G(x):=h^{n} - \Delta t \diff{}{x}\int\limits_{t^n}^{t^{n+1}}V \, dt, \label{sw2aSL}\\
			V^{n+1} &=& F(x) -\Delta t \frac{g}{2}\diff{{h}^{n} h^{n+1}}{x},
			\label{sw2bSL}
		\end{eqnarray}
		where the convective term $G(x)$ in the continuity equation \eqref{sw2aSL} is approximated relying on a SL scheme, following the approach which will be discussed for $F(x)$ hereafter. As a consequence, $h^{n+1}$ is explicitly solvable, thus it is directly substituted into the momentum equation \eqref{sw2bSL} for yielding $V^{n+1}$. The SL algorithm is used twice, namely for the term $G(x)$ and for $F(x)$, hence obtaining a formally semi-implicit method that does not require the solution of an implicit system, differently from the SL-IMEX-H method.
	
	\subsection{Space discretization}
	As done for the advection-diffusion equation, let $\Omega$ be the computational domain defined in the interval $[x_L;x_R]$ and let $N_x$ represent the total number of cells used to discretize $\Omega$. Each cell has a constant spacing of $\Delta x=(x_R-x_L)/N_x$ and a cell-centered space discretization is adopted, thus both conserved quantities are defined at the cell barycenter $x_i$, namely $h_i:=h(x_i)$ and $(hu)_i :=h(x_i)u(x_i)$.
	
	The discrete version of the differential operators $\partial_x$ and $\partial^2_{x^2}$ for a generic quantity $q(x)$ is based on a fourth order finite difference approximation, hence yielding
	\begin{eqnarray}
		\left. \frac{\partial q}{\partial x} \right|_{x=x_i}&=&\frac{1}{12 \Delta x}\left( -q_{i+2} + 8 q_{i+1}-  8 q_{i-1} +q_{i-2} \right), \\
		\left. \frac{\partial^2 q}{\partial x^2} \right|_{x=x_i} &=&\frac{1}{12 \Delta x^2}\left( -q_{i+2} + 16 q_{i+1} + 30 q_{i}  + 16 q_{i-1} -q_{i-2} \right).
	\end{eqnarray}
	
	In order to obtain a high order representation of the numerical solution within each computational cell, a CWENO reconstruction operator $\mathcal{R}$ is adopted (see \cite{LPR:99}). This allows smooth as well as discontinuous solutions to be properly treated ensuring non-oscillatory properties. 
	
	\subsection{Conservative semi-Lagrangian IMEX schemes} \label{sec.consSL}
	We present a conservative formulation of the SL method by considering only the advection term $F(x)$, given by \eqref{eqn.Fdef}, in the momentum equation \eqref{sw2b} for the SL-IMEX-H scheme. The term $G(x)$ in the continuity equation \eqref{sw2aSL} for the SL-IMEX algorithm can then be computed following the same approach outlined for $F(x)$.
	
	For the discretization of the convective contribution $F(x)$ in the momentum equation \eqref{sw2b}, the semi-Lagrangian schemes discussed in Section \ref{sec.SL-IMEX} can be used with no modifications. This would lead to a non-conservative transport of the momentum, which in principle does not represent a problem if incompressible flows are considered \cite{ConsSL2021} or parabolic PDE are likely to be solved. However, if shock waves are part of the eigenstructure of the governing system, as usually occur for hyperbolic PDE like the SWE, then the design of conservative methods is mandatory in order to capture the correct wave speeds and location of the discontinuities and the plateau exhibited by the solution. 
	
\begin{figure}[!ht]
		\begin{center}
			\includegraphics[width=0.9\textwidth]{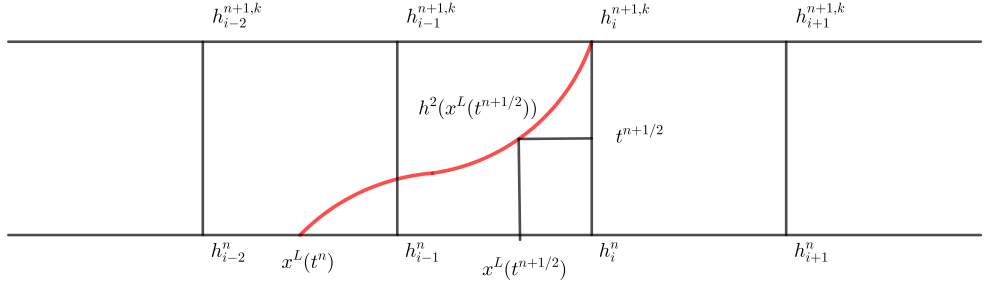}
			\caption{Space-time domain $\Omega_{t,x}$ used to develop the conservative SL-IMEX schemes.}
			\label{fig:HOcons}
		\end{center}
	\end{figure}

	Therefore we aim at constructing a conservative version of the semi-Lagrangian IMEX schemes previously introduced. We propose to use here a new philosophy in order to discretize $F(x)$, which has been very recently introduced in \cite{ConsSL2021} for scalar PDE. The idea is the following: given a starting point $x_i$, the Lagrangian trajectory of this point would travel through the flow lines according to \eqref{eqn.TrODE} up to the point $x_i^L$ at time $t^n$. Looking at Figure \ref{fig:HOcons}, let now $\Omega_{t,x}$ be the space-time domain that is the region bounded by the segments $[x_i^L;x_i]$ and $[t^n;t^{n+1}]$, which lies below the trajectory. 
	
	For the sake of simplicity, let us consider only the advection part of the PDE \eqref{sys1.1}, that is nothing but the simple advection equation
	\begin{equation}
		\diff{q}{t} + \diff{u q}{x} = 0.
		\label{eqn.adv}
	\end{equation}
	By defining the space-time divergence operator as $\nabla_{t,x}=(\partial_t, \partial_x)$, the integration of the advection equation \eqref{eqn.adv} over $\Omega_{t,x}$ and the use of the divergence theorem in the fully space-time framework leads to
	\begin{eqnarray}
		\int_{\Omega_{t,x}} \partial_t q +\partial_x (u q) \, dx dt &=&\int_{\Omega_{t,x}} \nabla_{t,x}\cdot (q,uq)\, dx dt \nonumber \\
		&=& -\int_{\Gamma^{L}}(q,uq)\cdot \vec{n} \, ds - \int_{x_i^L}^{x_i}q \, dx + \int_{t^n}^{t^{n+1}} uq \, dt, \label{adv_approx}
	\end{eqnarray}
	with $\Gamma^{L}$ denoting the boundary defined by the Lagrangian trajectory of the particle traveling from $x_i$ to $x_i^L$ (highlighted in red in Figure \ref{fig:HOcons}). By construction, since the tangent vector to $\Gamma^{L}$ is given by the velocity of the trajectory, the normal vector results to be $\vec{n}=\beta(1,-1/u)$ with $\beta=|u|/\sqrt{(u^2+1)}$ and therefore the integral along $\Gamma^{L}$ in \eqref{adv_approx} vanishes, i.e.
	\begin{equation}
		\int_{\Gamma^{L}}(q,uq)\cdot \vec{n} ds = \int_{\Gamma^{L}} {\beta}\left( q - \frac{1}{u}uq\right) \, ds = 0.
	\end{equation} 
	As a consequence, from \eqref{adv_approx} it follows a simple way to compute $H_i:=H(x_i)$ for the advection equation \eqref{eqn.adv} that only needs spatial information:
	\begin{equation}
		H_i := \int_{t^n}^{t^{n+1}} uq \, dt = \int_{x_i^L}^{x_i}q \, dx.
		\label{eqn.Hi}
	\end{equation} 
	Once $H_i$ is computed for each computational cell, a CWENO reconstruction operator $\mathcal{R}$ can be applied in order to obtain a piecewise high order polynomial representation of $H_i$. At this point, a classical conservative scheme is used to compute $F_i:=F(x_i)$, i.e.
	\begin{equation}
		F_i=q_i^n-\frac{1}{\Delta x}(\hat{H}_{i+1/2}-\hat{H}_{i-1/2}),
		\label{sw7}
	\end{equation}
	where $\hat{H}_{i\pm1/2}$ are numerical fluxes given by a local Lax-Friedrichs scheme.
	
	The boundary extrapolated values of $H(x)$ at the interfaces in \eqref{sw7} are provided by the CWENO reconstruction operator. Notice that in the conservative scheme \eqref{sw7} the time step $\Delta t$ is missing because it is already incorporated in the space-time integral \eqref{eqn.Hi}, thus it is contained in the integral over the interval $x=[x_i^L;x_i]$. Furthermore, this approach allows to reconstruct a purely spatial flux from the semi-Lagrangian discretization that is not possible in a classical formulation of the SL method due to the space-time nature of the Lagrangian trajectories. Therefore, in the construction of a conservative SL-IMEX scheme, for the explicit convective contribution we can use \eqref{sw7} to identify a flux as needed in the Eulerian IMEX scheme \eqref{e1}-\eqref{e4}, even though a semi-Lagrangian approach is adopted. The purely spatial part of the flux simply writes $\hat{H}_{i\pm1/2}/\Delta t$ because the time step is taken into account by the integral \eqref{eqn.Hi}. 
	
	\paragraph{Remark on the computation of the foot of the trajectory.} The conservative SL-IMEX-H scheme can easily be cast into the Eulerian advection form presented in Section \ref{sec.EulIMEX}. However, the velocity field has in any case to be advected according to the previously described Algorithm 1 in Section \ref{sec.SL-IMEX}. Since the trajectory equation associated to the evaluation of $H(x)$ in \eqref{eqn.Hdef} might involve a nonlinear transport, the coordinates $x^L_i$, needed for obtaining $H_i$ in \eqref{eqn.Hi}, have to be computed relying on a nonlinear ODE solver for the trajectory equation \eqref{eqn.TrODE}, such as the class of Runge-Kutta Exponential Integrators \cite{RKEI_Celledoni,CaiBos21,ConsSL2021}, or a predictor-corrector strategy \cite{RussoQiu_PC_2018}, or even a Taylor method \cite{TB19} which we rely on in this work.
	
	\bigskip
	Unfortunately it is easy to show that the conservative semi-Lagrangian approximation \eqref{eqn.Hi}-\eqref{sw7}, if applied to the shallow water equations, would lead to a first order (in time) discretization even for constant solutions. Indeed, let $u(x)$ and $h(x)$ be a steady solution of \eqref{sw1a}-\eqref{sw1b}, which might be non-constant in space. Then, the space-time integral \eqref{adv_approx} is not valid anymore since we also need to introduce the pressure fluxes that are part of the governing PDE system. Therefore, instead of \eqref{adv_approx}, integration of the momentum equation over the same space-time control volume $\Omega_{t,x}$ would lead to
	\begin{eqnarray}
		\int_{\Omega_{t,x}} \partial_t hu +\partial_x (h u^2) \, dx dt &=& \int_{\Omega_{t,x}} \nabla_{t,x}\cdot (hu,hu^2) \, dx dt \nonumber \\
		&=& - \int_{x_i^L}^{x_i}(hu)^n \, dx + \int_{t^n}^{t^{n+1}} hu^2 \, dt +\frac{g}{2}\int_{\Omega_{t,x}}\partial_x h^2 \, dx dt.
		\label{sw8}
	\end{eqnarray}
	Compared to \eqref{adv_approx}, the last term on the right hand side of \eqref{sw8} arises because of the pressure contribution in the momentum equation \eqref{sw1b}.
	Now, this additional multidimensional integral, which will be addressed with $\mathcal{P}$, may be split using the property of normal domains as follows
	\begin{equation}
		\mathcal{P}:=\frac{g}{2}\int_{\Omega_{t,x}}\partial_x h^2 \, dx dt=\frac{g}{2}\int_{t^n}^{t^n+1}{\int_{x^L(t)}^{x_i} \partial_x h^2 \, dx dt}= \frac{g}{2}\int_{t^n}^{t^{n+1}} h^2(x_i)-h^2(x^L(t)) \, dt,
		\label{sw9}
	\end{equation}
	where $x^L(t)$ is the foot of the characteristic located at time $t$, so that  $x^L(t^{n+1})=x_i$ and $x^L(t^n)=x_i^L=\LL(x_i,\Delta t, u)$. Integration of \eqref{sw9} can be done numerically with a quadrature rule of suitable order of accuracy, see Figure \ref{fig:HOcons} where the simple midpoint rule is illustrated. 
	
	\paragraph{Cauchy-Kovalevskaya procedure.} For achieving an order of accuracy greater than one, the quadrature nodes for the evaluation of the additional time integral $\mathcal{P}$ in \eqref{sw9} require the knowledge of $h$ inside the space-time domain. First, we observe that the only available high order information comes from the CWENO reconstruction polynomials defined at time $t^n$. Second, although a time reconstruction would provide the values of $h$ at quadrature nodes, it would demand the knowledge of the solution at $N$ points backward in time, thus needing a sort of initialization which prevents the algorithm to start from time $t_0$ (see for instance the class of BDF or Adams-Moulton schemes for ODE). Therefore, we propose to use a different strategy, that is based on the approximation of $h(t)$ by means of a Taylor series:
	\begin{equation}
		h(t) = h(t^n) + (t-t^n) \diff{h}{t} + \frac{(t-t^n)^2}{2} \diff{^2 h}{t^2} + \mathcal{O}(\Delta t^3),
		\label{eqn.hTaylor}
	\end{equation}
	which holds true for all cell values located at $x_i$ with $i=1 \ldots N_x$. To compute high order time derivatives we rely on the Cauchy-Kovalevskaya procedure which uses repeatedly the governing PDE to convert time partial derivatives into spatial partial derivatives. For the SWE \eqref{sw1a}-\eqref{sw1b}, the terms $\partial_t h$ and $\partial_{t^2}^2 h$ in \eqref{eqn.hTaylor} are evaluated as follows:
	\begin{subequations}
		\begin{align}
			\diff{h}{t} &= - \diff{V}{x}   \label{eqn.dht} \\
			\diff{^2h}{t^2} &= \frac{\partial}{\partial t} \left( - \diff{V}{x} \right) \nonumber \\ 
			&= \frac{\partial}{\partial x} \left( - \diff{V}{t} \right)  \nonumber \\
			& = \frac{\partial^2}{\partial x^2} \left( uV + \frac{g}{2} h^2 \right). \label{eqn.d2ht} 
		\end{align}
	\end{subequations}
	The computation of \eqref{eqn.dht} is readily available from the CWENO reconstruction operator applied to the conserved quantity $V$, i.e. $\mathcal{R}(V(x_i,t^n))$. The evaluation of the term \eqref{eqn.d2ht} requires two steps: i) first, the cell value of the flux term must be defined, that is 
	\begin{equation}
		G(x_i):= u_i V_i + \frac{g}{2} h_i^2 \qquad i=1\ldots N_x,
	\end{equation}
	ii) then, the CWENO reconstruction is used to obtain a high order approximation of $G(x_i)$, i.e. $\mathcal{R}(G(x_i))$, so that the second derivatives can be easily extracted and the term \eqref{eqn.d2ht} explicitly computed. Once the Taylor series \eqref{eqn.hTaylor} is completely determined, the value of $h^2$ needed in \eqref{sw9} can be obtained at any time level within the interval $[t^n;t^{n+1}]$. The integral $\mathcal{P}$ in \eqref{sw9} is numerically computed using quadrature rules of suitable accuracy, that are listed hereafter up to fourth order accurate schemes.
	\begin{itemize}
		\item Midpoint rule for second order of accuracy:
		\begin{equation}
			\mathcal{P} = \frac{g}{2} \Delta t \left[ h^2(x_i,t^{n+1/2})-h^2(x^L(t^{n+1/2})) \right] + \mathcal{O}(\Delta t^2).
			\label{eq:Pmidpoint}
		\end{equation}
		\item Kepler rule for up to fourth order of accuracy:
		\begin{eqnarray}
			\mathcal{P} &= & \frac{g}{2} \Delta t \, \frac{1}{6} \left[ h^2(x_i,t^n)-h^2(x^L(t^{n})) \right]  +\nonumber \\
			&& \frac{g}{2} \Delta t \, \frac{2}{3} \left[ h^2(x_i,t^{n+1/2})-h^2(x^L(t^{n+1/2})) \right]  +\nonumber \\
			&& \frac{g}{2} \Delta t \, \frac{1}{6} \left[ h^2(x_i,t^{n+1})-h^2(x^L(t^{n+1})) \right]  + \mathcal{O}(\Delta t^4),
			\label{eq:Psimpson}
		\end{eqnarray}
		where the solutions $h^2(x_i,t^{n+1/2})$ and $h^2(x_i,t^{n+1})$ are readily computed with the Taylor expansion \eqref{eqn.hTaylor}.
	\end{itemize}

	We recall here that higher order time stepping can be achieved by applying the implicit-explicit discretization presented above to the class of high order IMEX schemes given by \eqref{eqn.IMEXRK}-\eqref{e4} with the tableaux reported in \ref{app.IMEX}.
	
	\paragraph{Remark.} The additional term $\mathcal{P}$ in \eqref{sw9} can be written, using the mean value theorem, as
	\begin{equation}
		\frac{g}{2}\int_{t^n}^{t^{n+1}} h^2(x_i)-h^2(x^L(t)) dt={ \frac{g}{2}}\left[h^2(x_i,t_c)-h^2(x^L(t_c))\right] \, \Delta t \approx { \frac{g}{2}} \delta(t_c) \, \Delta t,
	\end{equation}
	where $t_c \in [t^n;t^{n+1}]$ and $\delta(t_c)$ is the mean value of the integrand function. For a non-constant solution $h=h(x,t)$, or even in the case of steady non-constant solutions like $h=h(x)$, neglecting this contribution as done in \eqref{adv_approx} would introduce a first order local truncation error. 
	
	\section{Asymptotic-preserving semi-Lagrangian IMEX schemes} \label{sec.AP}
	In order to extend the previous scheme \eqref{sw2a}-\eqref{sw2b} with the novel conservative semi-Lagrangian IMEX methods to the case of hyperbolic balance laws, let us consider the shallow water equations with a relaxation term \cite{PR_IMEX}, that is given by
	\begin{eqnarray}
	&&\diff{h}{t}+\diff{V}{x}=0, \label{sw1ra}\\
	&&\diff{V}{t}+\diff{uV}{x}+\frac{g}{2}\diff{h^2}{x}=\frac{h}{\varepsilon}\left( \frac{h}{2} - u \right),
	\label{sw1rb}
	\end{eqnarray}
	with $\varepsilon$ being a relaxation parameter. In the stiff limit, i.e. when $\varepsilon \to 0$, system \eqref{sw1ra}-\eqref{sw1rb} reduces to the inviscid Burgers equation
	\begin{equation}
	\diff{h}{t} + \frac{\partial}{\partial x} \left(\frac{h^2}{2}\right) = 0.
	\label{eqn.Burgers}
	\end{equation}
	On the other hand, if $\varepsilon \to \infty$ the shallow water PDE system \eqref{sw2a}-\eqref{sw2b} is consistently retrieved because the source term simply vanishes.

	Let us now propose two different semi-discrete AP schemes for the relaxation system \eqref{sw1ra}-\eqref{sw1rb}, which arise from the SL-IMEX-H and SL-IMEX methods illustrated in Section \ref{sec.consSL}.
	
	\paragraph{SL-IMEX-H AP method.} The first method directly follows from \eqref{sw2a}-\eqref{sw2b} and yields
	\begin{eqnarray}
	h^{n+1}&=&h^{n} - \Delta t \diff{V^{n+1}}{x}, \label{sw2ra}\\
	V^{n+1}&=&F(x) -\Delta t \frac{g}{2}\diff{{h}^{n} h^{n+1}}{x} + \Delta t \, \left( \frac{{h}^{n} h^{n+1}}{2 \varepsilon} - \frac{V^{n+1}}{\varepsilon} \right),
	\label{sw2rb}
	\end{eqnarray}
	where the relaxation source is discretized implicitly, so that the second term is linear with respect to $V^{n+1}$ and can be readily inverted. Notice that the conservative semi-Lagrangian scheme is already embedded in the definition of the term $F(x)$ according to \eqref{eqn.Fdef}. From \eqref{sw2rb}, by solving for $V^{n+1}$, we obtain
	\begin{eqnarray}
	V^{n+1} \gamma =F(x) -\Delta t \frac{g}{2}\diff{{h}^{n} h^{n+1}}{x} + \Delta t \frac{{h}^{n} h^{n+1}}{2 \varepsilon}, \qquad \gamma = \frac{\varepsilon+\Delta t}{\varepsilon},
	\end{eqnarray}
	that can be substituted into the continuity equation \eqref{sw2ra}, hence yielding
	\begin{equation}
	h^{n+1} = h^{n} - \Delta t \frac{\partial }{\partial x} \left( \frac{F(x)}{\gamma} \right) + \frac{\Delta t^2}{\gamma} \frac{g}{2} \frac{\partial }{\partial x} \left( \diff{{h}^{n} h^{n+1}}{x} \right) - \frac{\Delta t^2}{\gamma \varepsilon} \frac{\partial }{\partial x} \left( \frac{{h}^{n} h^{n+1}}{2} \right).
	\label{eqn.hAP}
	\end{equation}
	Now, the limit of the relaxation system \eqref{sw2a}-\eqref{sw2b} can be analyzed starting from the semi-discrete form \eqref{eqn.hAP}, which represents the novel conservative SL scheme introduced in this work. In the limit for $\varepsilon \to 0$ it holds
		\begin{equation}
		\lim\limits_{\varepsilon \to 0} \frac{1}{\gamma} = 0, \qquad \lim\limits_{\varepsilon \to 0} \frac{1}{\gamma \varepsilon} = \frac{1}{\Delta t},
		\label{eqn.limZero}
		\end{equation}
		therefore the scheme given by \eqref{eqn.hAP} simplifies to 
		\begin{equation}
		h^{n+1} = h^{n} - \Delta t \frac{\partial }{\partial x} \left( \frac{{h}^{n} h^{n+1}}{2} \right),
		\label{eqn.hstiff}
		\end{equation}
		which is a consistent discretization of the limit model \eqref{eqn.Burgers}.
	The proposed SL-IMEX-H method is asymptotic preserving (AP). Furthermore, the semi-Lagrangian discretization does not affect the AP property, which is ensured by the IMEX time stepping. 
	
	\paragraph{SL-IMEX AP method.} The second discretization is based on the scheme \eqref{sw2aSL}-\eqref{sw2bSL}, thus the relaxation system \eqref{sw1ra}-\eqref{sw1rb} can be discretized in time as
	\begin{eqnarray}
	h^{n+1}&=&G(x), \qquad G(x):=h^{n} - \Delta t \diff{}{x}\int\limits_{t^n}^{t^{n+1}}V \, dt, \label{sw2raSL}\\
	V^{n+1} &=& \frac{F(x)}{\gamma} -\frac{\Delta t}{\gamma} \frac{g}{2}\diff{{h}^{n} h^{n+1}}{x} + \frac{\Delta t}{\gamma\varepsilon} \frac{{h}^{n} h^{n+1}}{2 }, \qquad \gamma = \frac{\varepsilon+\Delta t}{\varepsilon},
	\label{sw2rbSL}
	\end{eqnarray}
	Similarly to what discussed above,  we can characterize the limit of \eqref{sw2raSL}-\eqref{sw2rbSL} for $\varepsilon \rightarrow 0$. According to \eqref{eqn.limZero},  from the momentum equation \eqref{sw2rbSL} one has $V^{n+1}\to (h^{n+1})^2/2$. Therefore, $u(=V/h)\to h^{n+1}/2$ and the continuity equation becomes
		\begin{equation}
		h^{n+1}= h^{n} - \Delta t \diff{}{x}\int\limits_{t^n}^{t^{n+1}}u h \, dt = h^{n} - \Delta t \diff{}{x}\int\limits_{t^n}^{t^{n+1}}h^2 \, dt,
		\label{eqn.BurgersSL}
		\end{equation}
		which is a semi-Lagrangian discretization of the inviscid Burgers equation \eqref{eqn.Burgers}. Notice that the advection speed in \eqref{eqn.BurgersSL} is consistently retrieved, hence leading to an AP method.

   The SL-IMEX scheme \eqref{sw2raSL}-\eqref{sw2rbSL} is formally unconditionally stable and completely explicit from the computational viewpoint. Nevertheless, in case of shocks or strong discontinuity exhibited by the solution, it is necessary to supplement the continuity equation \eqref{sw2raSL} with additional numerical viscosity according to \eqref{eqn.visc}. As already done, we adopt an implicit treatment of the numerical dissipation, that leads to the solution of a linear system for the unknown $h^{n+1}$. This algorithm makes use of the conservative SL method for the advection of both mass and momentum, thus it is simply denoted with SL-IMEX.
   
   Higher order SL-IMEX AP methods can be derived by adopting the same implicit-explicit partitioning here described in the high order IMEX setting \eqref{eqn.IMEXRK}-\eqref{e4}. We omit the details for brevity.
   
    \paragraph{Remark on the low Froude number limit.}
    The approach here designed permits to capture also diffusive behaviors of the shallow water system with relaxation terms. Specifically, the stiffness related to the pressure gradient in the SWE can be highlighted by introducing a scaling based on the Froude number $Fr=|u|/\sqrt{g h}$ and assuming dimensionless variables. Following \cite{Kurganov2019}, the shallow water system supplemented with the relaxation source term introduced in \eqref{sw1ra}-\eqref{sw1rb}, can be rewritten in the following form
    \begin{eqnarray}
	&&\diff{h}{t}+\diff{V}{x}=0, \label{sw1ra2}\\
	&&\diff{V}{t}+\diff{uV}{x}+\frac{1}{2Fr^2}\diff{h^2}{x}=\frac{h}{\varepsilon}\left( \frac{h}{2} - u \right),
	\label{sw1rb2}
	\end{eqnarray}
   with eigenvalues $\lambda_{1,2}=u \left(1 \pm 1/Fr\right)$. The derivation of the asymptotic limit of the shallow water equations in the low Froude regime has been studied for instance in \cite{Lukacova2014,Berthon2015} and references therein. Here, we limit us to provide the asymptotic limit of \eqref{sw1ra2}-\eqref{sw1rb2}. Assuming $\varepsilon = O(Fr^2)$, in the limit $\varepsilon\to 0$ system \eqref{sw1ra2}-\eqref{sw1rb2} reduces to the viscous Burgers equation
	\begin{equation}
	\diff{h}{t} + \frac{\partial}{\partial x} \left(\frac{h^2}{2}\right) = \frac12\frac{\partial^2 h^2}{\partial x^2}.
	\label{eqn.Burgersv}
	\end{equation}
    To obtain a numerical scheme which is in principle unconditionally stable, despite the physical scale of the problem under consideration, one might discretize the advection flux term with the novel SL algorithm proposed in this work, while the diffusive terms relying on an implicit finite difference scheme, thus ensuring stability regardless the chosen time step. Under the assumption $\varepsilon = O(Fr^2)$, the semi-discrete scheme reads then as follows:
    \begin{eqnarray}
    	h^{n+1}&=&G(x), \qquad G(x):=h^{n} - \Delta t \int \limits_{t^n}^{t^{n+1}} \diff{V}{x} \, dt, \label{sw2raFr}\\
    	V^{n+1} &=& \frac{F(x)}{\gamma} -\frac{\Delta t}{2\gamma \varepsilon} \diff{{h}^{n} h^{n+1}}{x} + \frac{\Delta t}{\gamma\varepsilon} \frac{{h}^{n} h^{n+1}}{2 }, \qquad \gamma = \frac{\varepsilon+\Delta t}{\varepsilon}.
    	\label{sw2rbFr}
    \end{eqnarray}  
	In the limit $\varepsilon\to 0$, the momentum equation \eqref{sw2rbFr} becomes
	\begin{equation}
		V^{n+1} = -\frac{1}{2} \diff{{h}^{n} h^{n+1}}{x} + \frac{{h}^{n} h^{n+1}}{2 },
	\end{equation}
	which, after substitution into the continuity equation \eqref{sw2raFr}, yields the discrete version of the limit equation \eqref{eqn.Burgersv}:
	\begin{equation}
		h^{n+1} - \hat{G}(x) = \frac{\Delta t}{2} \diff{^2 {h}^{n} h^{n+1}}{x^2}, \qquad \hat{G}(x) := h^n + \Delta t \frac{\partial}{\partial x}\int \limits_{t^n}^{t^{n+1}} \frac{h^2}{2} \, dt. 
	\end{equation}	
	We underline that the advection flux $\hat{G}(x)$ can be computed using the conservative SL schemes, while implicit finite differences might be employed for the diffusive terms.

	\section{Numerical results for the shallow water system} \label{sec.SWE_test}
	In this section, we perform some numerical test cases which aim at demonstrating the accuracy and the robustness of the conservative version of the SL-IMEX schemes. First, numerical convergence studies are presented to show that the formal order of accuracy is achieved. Second, a smooth solution involving a pressure wave propagation is used to highlight the benefits of the high order method in terms of energy dissipation. The total energy $E$ is computed as the sum of kinetic $K$ and potential $U$ energy contribution, which are given by
	\begin{equation}
		E = K + U, \qquad K = \int \limits_{\Omega} \frac{h u^2}{2} \, dx, \qquad U = \int \limits_{\Omega} \frac{g h^2 }{2} \, dx.
		\label{eqn.enDef}
	\end{equation}
	Third, two Riemann problems are solved which deal with shock waves, where a conservative scheme is crucial for obtaining the correct propagation speed. Fourth, the asymptotic preserving property of the SL-IMEX schemes is verified by solving a set of smooth and discontinuous test cases for the relaxation system of the SWE in the stiff limit. 
	
	The CFL stability condition for an explicit solver of the SWE \eqref{sw1a}-\eqref{sw1b} is given by
	\begin{equation}
		\textnormal{CFL} = \frac{\Delta t \, \max \left(|u| + \sqrt{g h}\right) }{\Delta x},
		\label{eqn.CFL}
	\end{equation}
	which leads to the condition $\textnormal{CFL}<1$. Because of both the semi-Lagrangian discretizations of the convection terms (SL-IMEX and SL-IMEX-H) and the implicit treatment of the pressure flux in \eqref{sw2b} and \eqref{sw2bSL}, the proposed method is unconditionally stable for any time step $\Delta t$. As a consequence, the time step can be defined only by the physical time scale of the problem under consideration and not by numerical restrictions. For each test case shown in the sequel, we explicitly report the resulting CFL number, according to the definition \eqref{eqn.CFL}. 
	
	\subsection{Numerical convergence studies}
	Here, a smooth steady state problem is considered to measure the order of accuracy of the SL-IMEX-H schemes. Notice that the semi-Lagrangian approach directly solves the Lagrangian form of the advection term \eqref{sys1.2} and, as such, even a steady solution is not trivial to be maintained since the algorithm formally transports backward and forward the numerical solution. On the other hand, we propose a steady state solution of the SWE \eqref{sw1a}-\eqref{sw1b} because it easily allows an analytical solution to be derived, so that convergence studies can be carried out.
	
	We consider a computational domain $\Omega=[10;10]$ and we prescribe the following initial condition for the fluid velocity at $t=0$:
	\begin{equation}
		u(x)= 1 + a \cos\left( \frac{2\pi}{10}x \right), \qquad a \in \mathds{R}.
		\label{eqn.uex}
	\end{equation} 
	To obtain a steady solution for the momentum, the advection contribution must be exactly balanced by the pressure fluxes, which means solving the following ODE:
	\begin{equation}
		\frac{d hu^2}{d x} = - gh \frac{d h}{x},
	\end{equation}
	that yields the sought water depth 
	\begin{equation}
		\begin{aligned}
			h(x) = &-\frac{a^2 \cos(\alpha x)^2+1+2a\cos(\alpha x)}{g} +  \\
			& -\frac{1}{g} \sqrt{a^4\cos(\alpha x)^4+6 a^2\cos(\alpha x)^2 + 4 a^3 \cos(\alpha x)^3 + 1 + 4 a \cos(\alpha x) - 2g c}, \qquad \alpha = \frac{\pi}{5}.
		\end{aligned}
		\label{eqn.hex}
	\end{equation}
	We set $a=5$ in \eqref{eqn.uex} and $c=-1$ in \eqref{eqn.hex}, then the continuity equation \eqref{sw1a} must be supplemented with a source term $S(h,V)$ that maintains constant over time the water depth given by \eqref{eqn.hex}, thus we solve
	\begin{equation}
		h^{n+1}=h^{n} - \Delta t \diff{V^{n+1}}{x} + \Delta t \, S(h,V)^{n+1}, \qquad S(h,V)^{n+1}=\diff{V^{n+1}}{x},
		\label{sw2a_s}
	\end{equation}
	where the source term is discretized implicitly so that we ensure a perfectly compatible discretization that keeps a constant water depth up to machine precision.
	
\begin{figure}[!htbp]
		\begin{center}
			\begin{tabular}{ccc} 
				\includegraphics[width=0.33\textwidth]{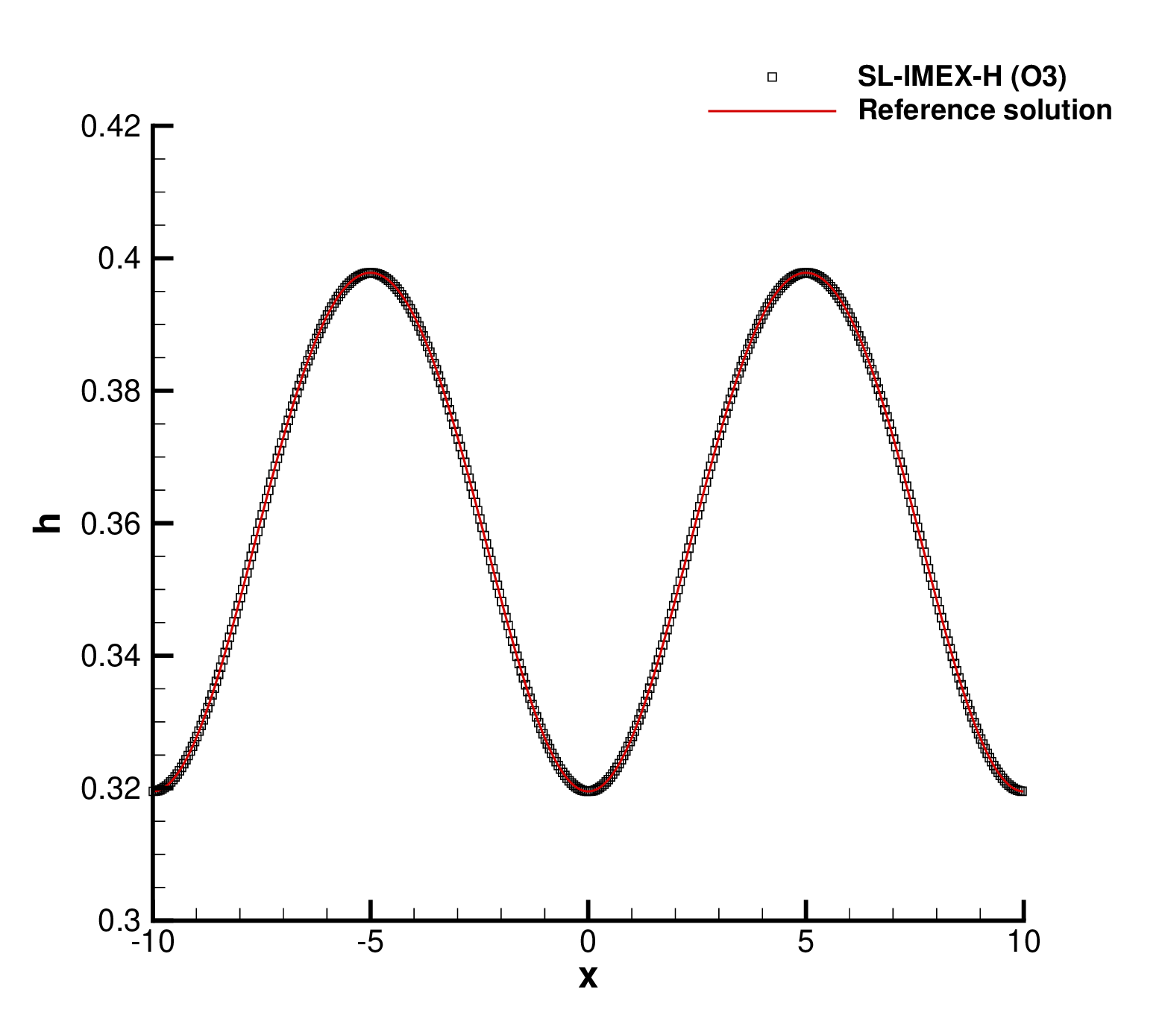}  &                 
				\includegraphics[width=0.33\textwidth]{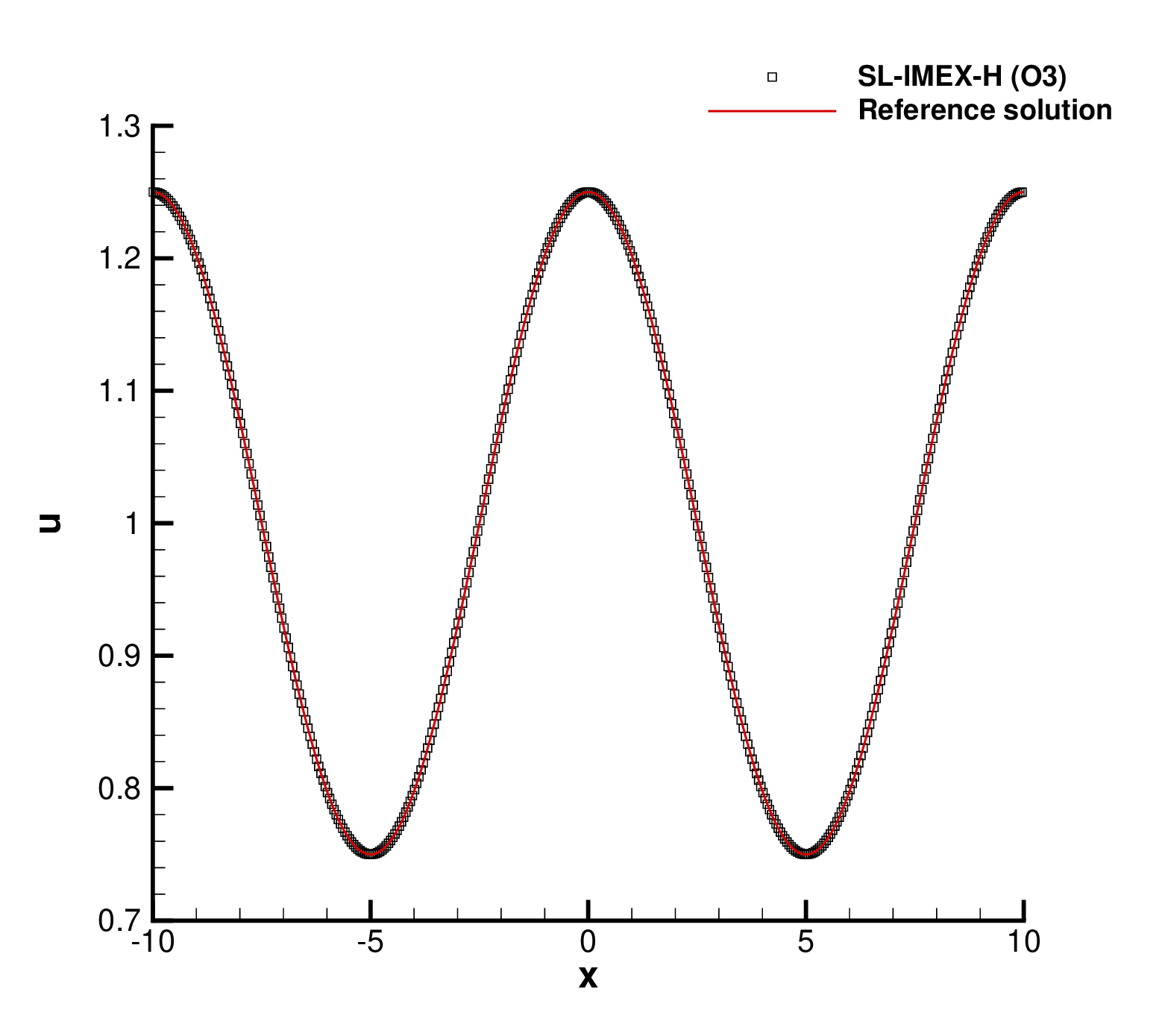} &
				\includegraphics[width=0.33\textwidth]{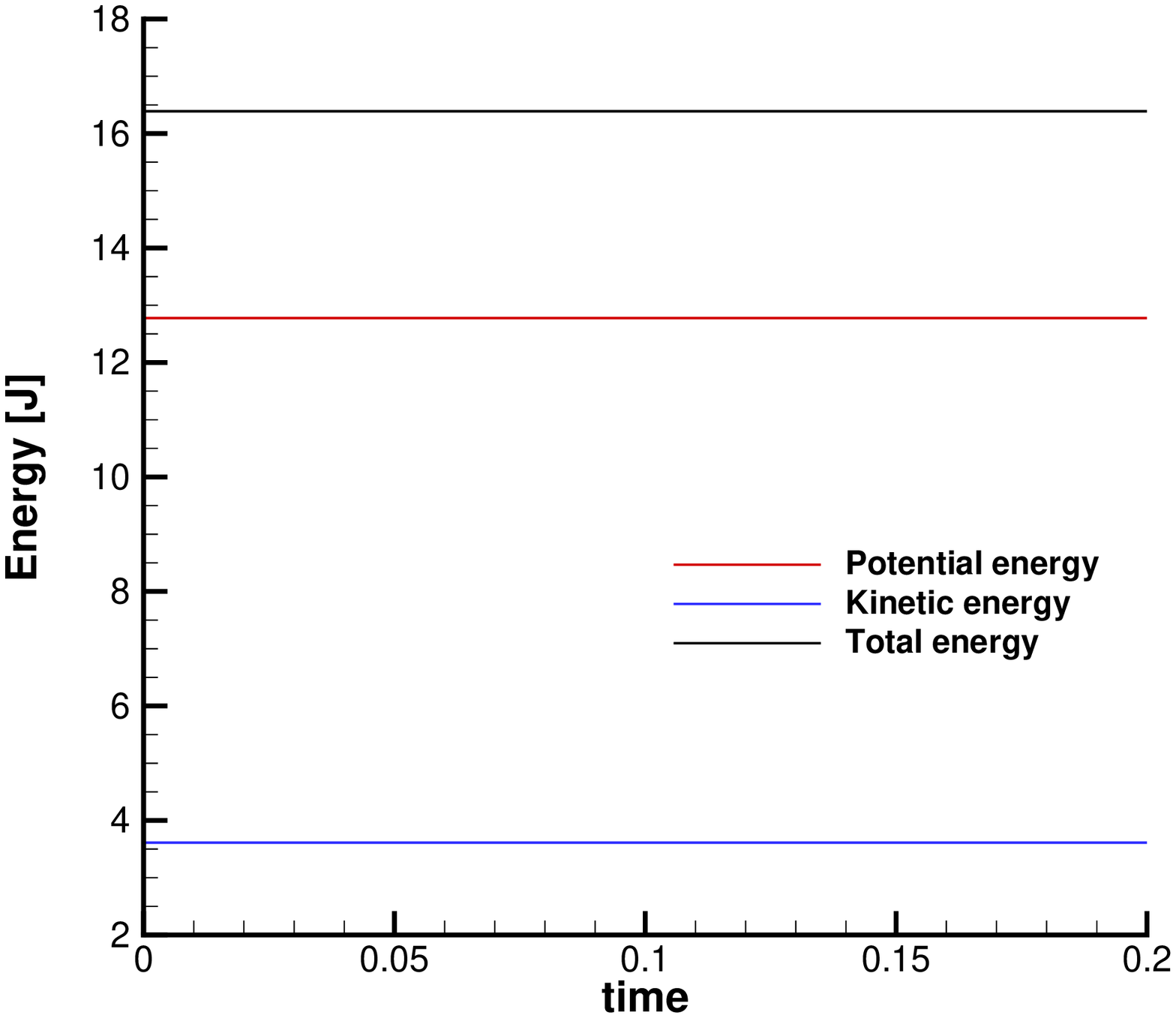} \\
			\end{tabular}
			\caption{Numerical convergence test with third order SL-IMEX-H scheme at time $t=0.2$ with $N_x=400$. Comparison against the reference solution for water height (left) and velocity (middle), and time evolution of potential, kinetic and total energy (right).}
			\label{fig.SWE-convergence}
		\end{center}
	\end{figure}

	The simulation is run until the final time $t_f=0.2$ with CFL=2 on a sequence of five successively refined meshes obtained by a refinement factor of $2$ applied to the number of cells, i.e. $N_x=\{50,100,200,400,800\}$. Figure \ref{fig.SWE-convergence} shows a comparison between a third order numerical solution with $N_x=400$ and the analytical solution, as well as the time evolution of the energy which is perfectly preserved constant over time.
	
	Table \ref{tab.conv} reports the convergence studies for second and third order SL-IMEX-H schemes, where errors are measured in $L_\infty$ norm for the momentum $V$, that is
	\begin{equation}
		L_\infty = \max_i \left|V_i-V_i^{ex}\right|, \qquad i=1 \ldots N_x,
	\end{equation}
	with the exact momentum $V_i^{ex}= h(x_i) u(x_i)$ computed from the analytical solution \eqref{eqn.uex}-\eqref{eqn.hex}. The formal order of accuracy is achieved by the proposed numerical method, which is possible if the additional integral $\mathcal{P}$ in \eqref{sw9} is properly considered. Numerical evidences reported in Table \ref{tab.conv} demonstrate that the SL-IMEX-H schemes, although using a high order CWENO reconstruction and a high order IMEX scheme, can achieve only up to first order of accuracy if the integral $\mathcal{P}$ is neglected, hence solving \eqref{adv_approx} for the convective terms.
	
	\begin{table}[!htbp]  
		\begin{center} 
			\begin{small}
				\renewcommand{\arraystretch}{1.2}
				\begin{tabular}{cccccc} 
					\multicolumn{6}{c}{SL-IMEX-H $\mathcal{O}2$} \\
					\hline
					\hline
					& \multicolumn{2}{c}{with $\mathcal{P}$} & & \multicolumn{2}{c}{without $\mathcal{P}$} \\
					\hline
					$\Delta x$ & $L_\infty$ & $\mathcal{O}(V)$ &  & $L_\infty$ & $\mathcal{O}(V)$ \\ 
					0.4   & 1.3032E-04 & -    & & 3.5714E-04 & -    \\
					0.2   & 2.9772E-05 & 2.13 & & 1.9069E-04 & 0.91 \\
					0.1   & 7.0589E-06 & 2.08 & & 1.0655E-04 & 0.84 \\
					0.05  & 1.7404E-06 & 2.02 & & 5.4242E-05 & 0.97 \\
					0.025 & 4.4074E-07 & 1.98 & & 2.8017E-05 & 0.95 \\
					\hline
					\multicolumn{6}{c}{} \\
					\multicolumn{6}{c}{SL-IMEX-H $\mathcal{O}3$} \\
					\hline
					\hline
					& \multicolumn{2}{c}{with $\mathcal{P}$} & & \multicolumn{2}{c}{without $\mathcal{P}$} \\
					\hline
					$\Delta x$ & $L_\infty$ & $\mathcal{O}(V)$ &  & $L_\infty$ & $\mathcal{O}(V)$ \\ 
					0.4   & 1.2011E-04 & -    & & 5.6089E-04 & -    \\
					0.2   & 1.5396E-05 & 2.96 & & 3.0007E-04 & 0.91 \\
					0.1   & 1.9865E-06 & 2.95 & & 1.6725E-04 & 0.84 \\
					0.05  & 2.7466E-07 & 2.85 & & 8.5115E-05 & 0.97 \\
					0.025 & 4.8296E-08 & 2.51 & & 4.3960E-05 & 0.95 \\
					\hline
				\end{tabular}
			\end{small}
		\end{center}
		\caption{Numerical convergence results for the one-dimensional shallow water equations with CFL=2 using second and third order SL-IMEX-H schemes with and without the time integral term $\mathcal{P}$ in the semi-Lagrangian scheme given by \eqref{sw9}. The errors are measured in $L_\infty$ norm and refer to the momentum $V$ at time $t=0.2$.}  
		\label{tab.conv}
	\end{table}
	
	\subsection{Pressure wave propagation}
\begin{figure}[!htbp]
		\begin{center}
			\begin{tabular}{cc} 
				\includegraphics[width=0.47\textwidth]{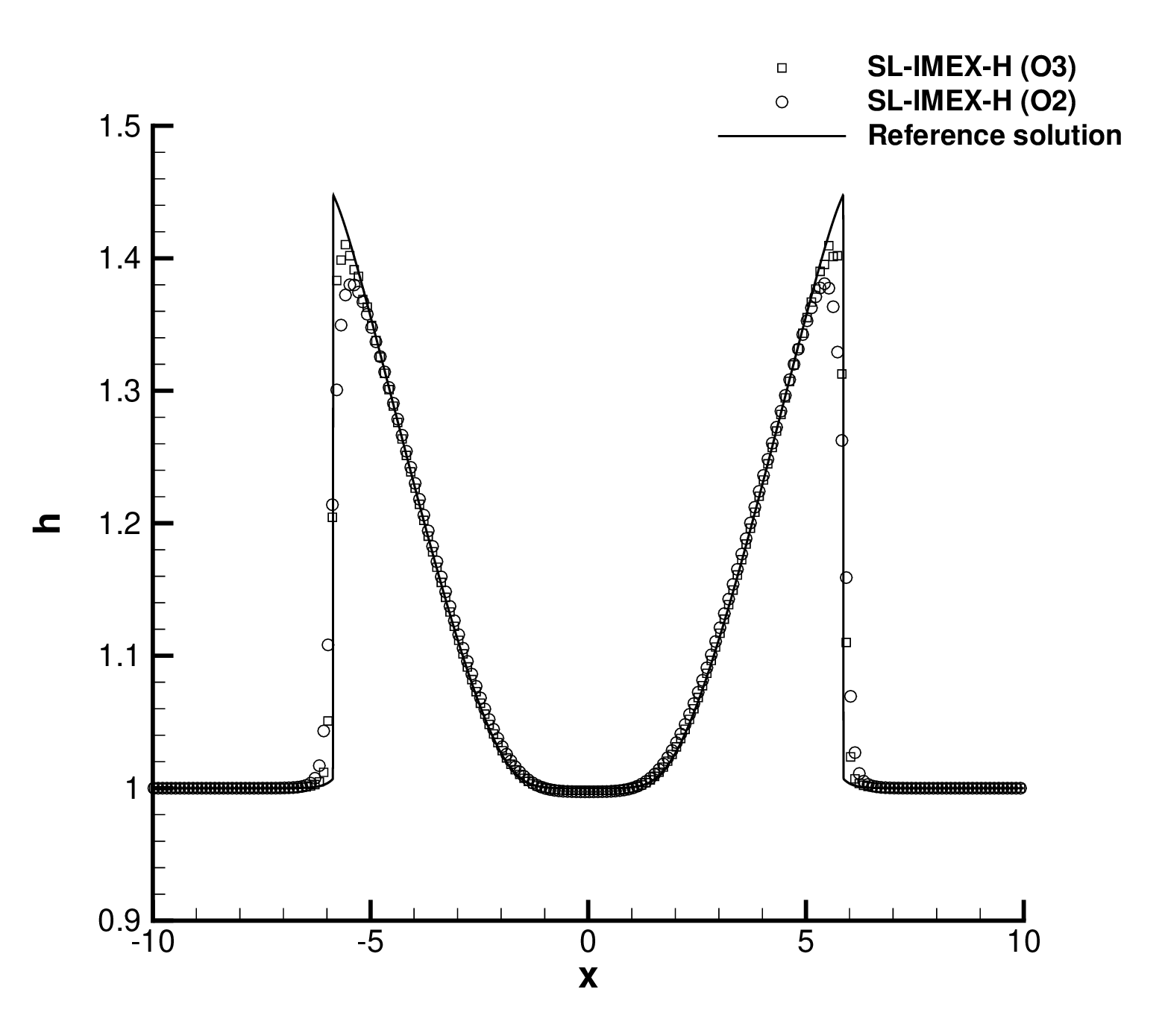}  &                 
				\includegraphics[width=0.47\textwidth]{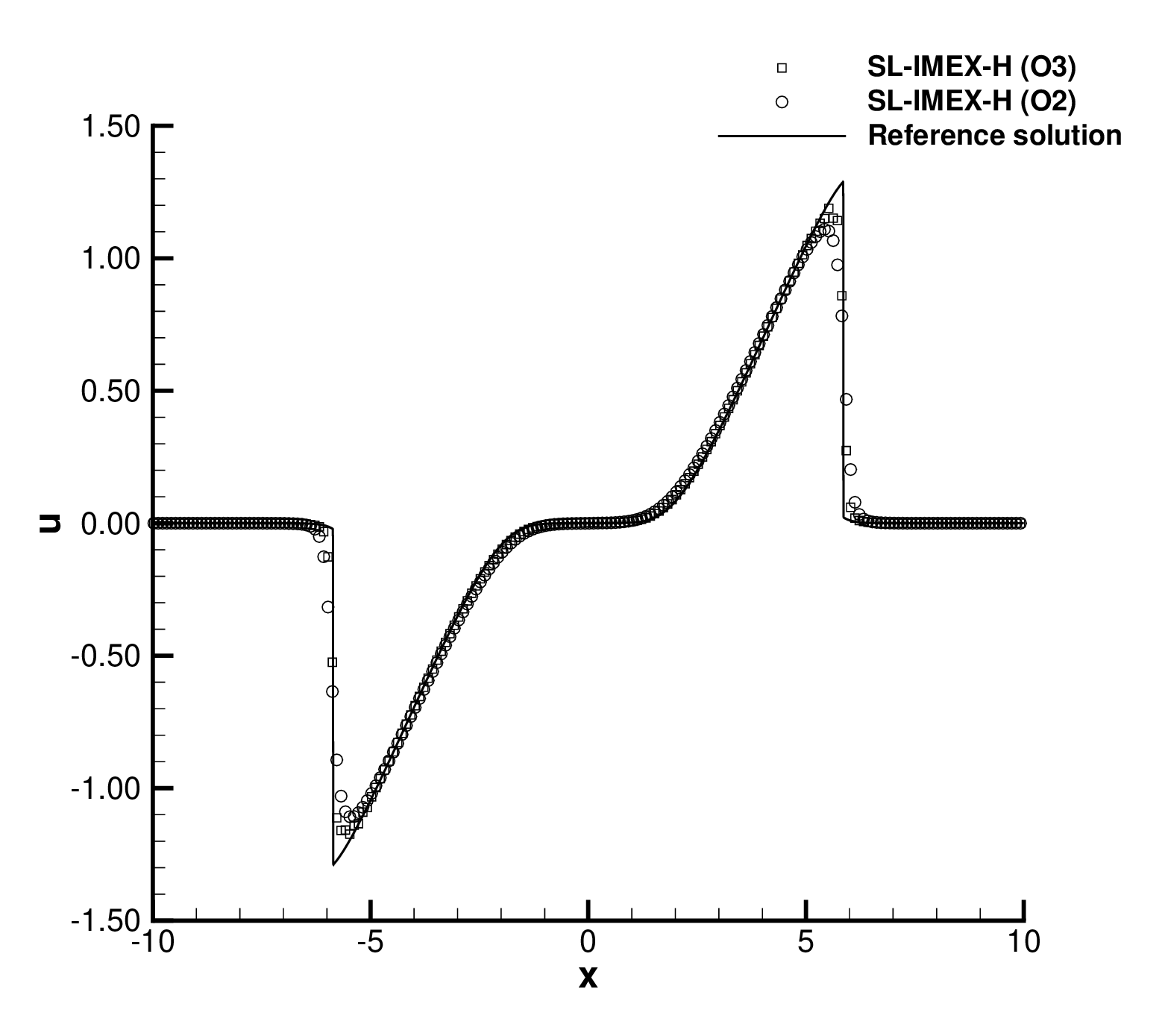} \\
				\includegraphics[width=0.47\textwidth]{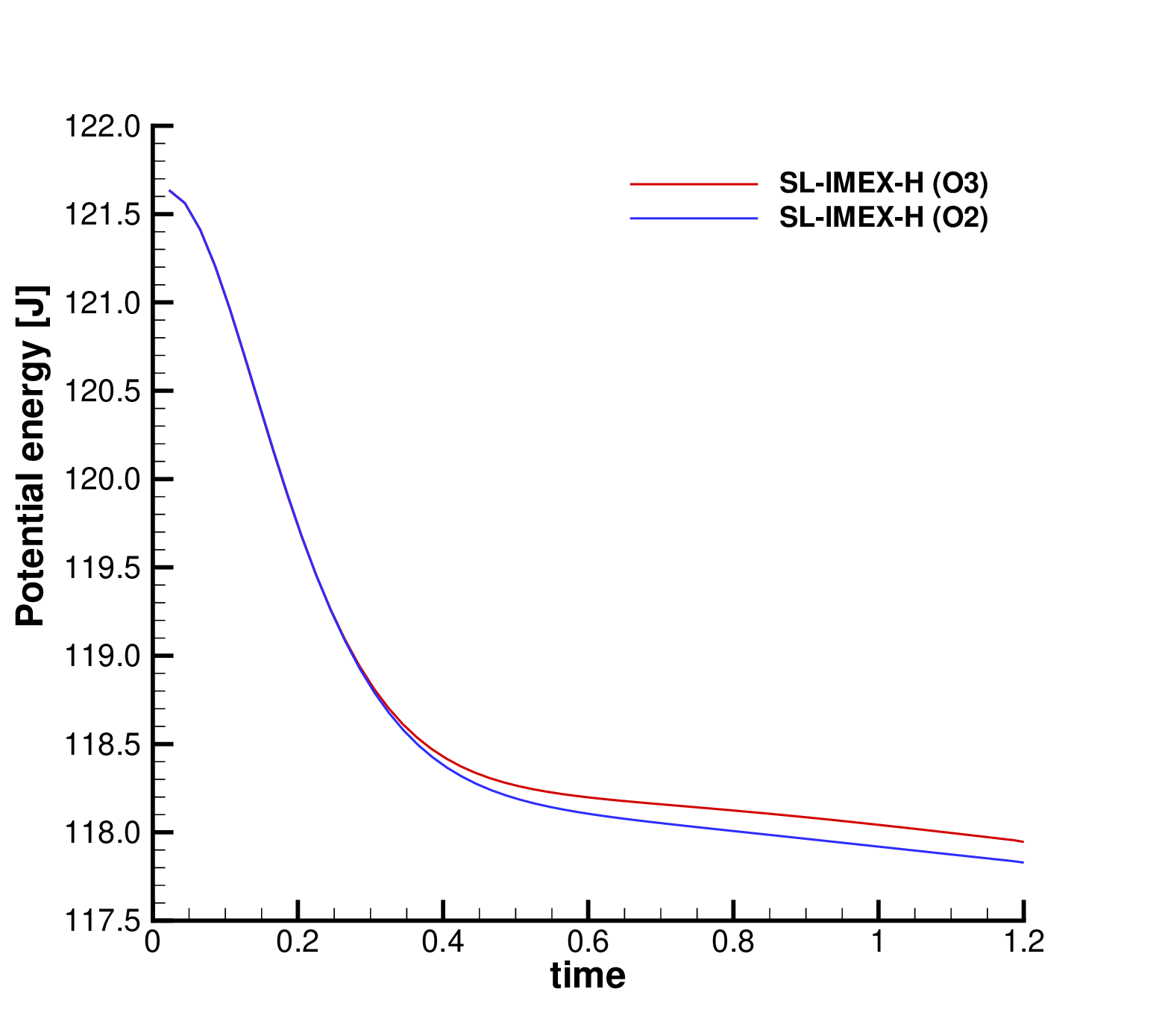}  &                 
				\includegraphics[width=0.47\textwidth]{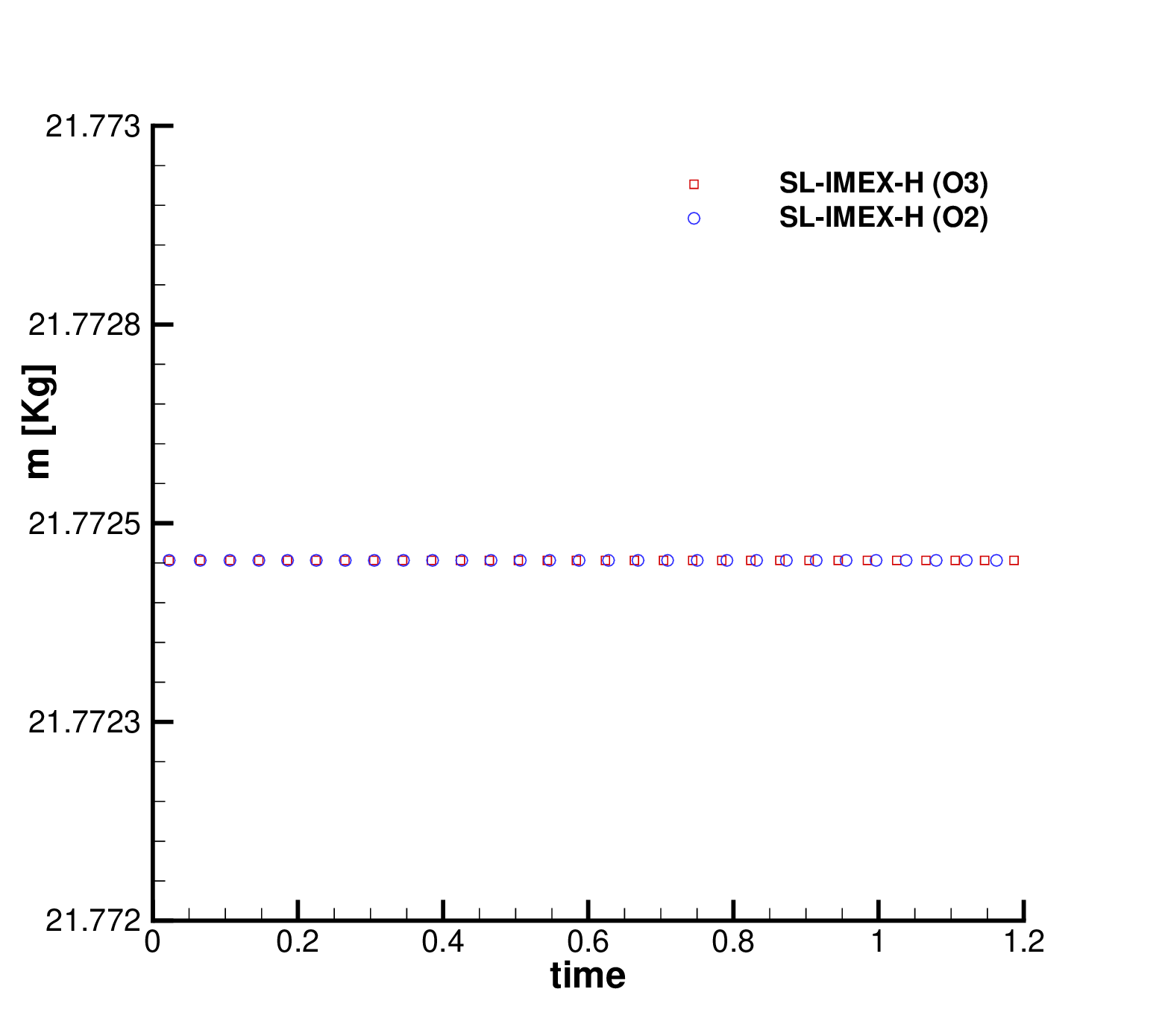} \\
			\end{tabular}
			\caption{Pressure wave propagation at time $t=1.2$ with CFL=2. Top: comparison of second and third order numerical results against the reference solution for water height (left) and velocity (right). Bottom: time evolution of the potential energy (left) and total mass (right) for both second and third order SL-IMEX-H schemes.}
			\label{fig.SWE-Gaussian}
		\end{center}
	\end{figure}

	We consider an initially smooth solution that involves the propagation of a pressure wave over the computational domain $\Omega=[-10;10]$ which is discretized with a total number of cells $N_x=400$. The water is initially at rest and the water depth is assigned a Gaussian profile
	\begin{equation}
		h(x,0) = 1 + e^{-x^2}.
	\end{equation} 
	The setup of the test is symmetric, thus two waves are departing from the center of the domain at $x=0$ and traveling towards the boundaries. The initial profile of the water front is smooth, but a shock wave arises after $t \approx 0.9$. The reference solution is obtained with a MUSCL-TVD scheme on a very fine mesh composed of $20000$ cells, which is used to compare the numerical solution of the SL-IMEX-H methods in Figure \ref{fig.SWE-Gaussian}. The second order scheme is more dissipative compared to the third order version, that is highlighted in the dissipation of the potential energy $U$ according to \eqref{eqn.enDef}. Finally, we also monitor the conservation of the total mass given by
	\begin{equation}
		m=\int \limits_{\Omega} h \, dx,
		\label{eqn.m}
	\end{equation}
	which is satisfied up to machine precision thanks to the conservative version of the SL schemes.  
	
	\subsection{Riemann problems}
	Here, we propose some computational results for two different Riemann problems. The Riemann problem is a special initial value problem for system \eqref{sw1a}-\eqref{sw1b}, where a discontinuity separates two piecewise constant states that represent the initial data. Since these problems involve shock waves, they are traditionally solved by explicit conservative methods, while here we make use of the novel class of implicit-explicit schemes with a conservative semi-Lagrangian treatment of the advection terms. The initial condition is given in terms of the state vector $\Q(x,0)=(h(x,0),V(x,0))$ assigned to the left $L$ and the right $R$ part of the domain, that is
	\begin{equation}
		\Q(x,0) = \left\{ \begin{array}{lll}
			\Q_L & \textnormal{if} & x \leq x_d \\
			\Q_R & \textnormal{if} & x > x_d 
		\end{array} \right.,
		\label{eqn.RPini}
	\end{equation}
	with $x_d$ denoting the location of the initial discontinuity. The initial condition for the two Riemann problems is reported in Table \ref{tab:init} and the computational domain is discretized using $N_x=400$ control volumes. The exact solution can be computed relying on the Riemann solver given in \cite{toro-swe-book}.
	
	\begin{table}[!htbp]  
		\caption{Initialization of Riemann problems. Initial states left (L) and right (R) are reported as well as the final time of the simulation $t_f$, the computational domain $[x_L;x_R]$ and the position of the initial discontinuity $x_d$.}  
		\begin{center} 
			\begin{small}
				\renewcommand{\arraystretch}{1.0}
				\numerikNine
				\begin{tabular}{l|c|ccc|cc|cc} 
					\hline
					Name & $t_{f}$ & $x_L$ & $x_R$ & $x_d$ & $h_L$ & $u_L$ & $h_R$ & $u_R$ \\
					\hline
					RP1 & 1.0 & -10.0 & 10.0 & 0.0 & 1.5 & -1.0  & 1.0 & 2.0 \\
					RP2 & 1.5 & -10.0 & 10.0 & 0.0 & 1.0 &  0.0  & 0.5 & 0.0 \\
					\hline
				\end{tabular}
			\end{small}
		\end{center}
		\label{tab:init}
	\end{table}
	
\begin{figure}[!htbp]
		\begin{center}
			\begin{tabular}{cc} 
				\includegraphics[width=0.47\textwidth]{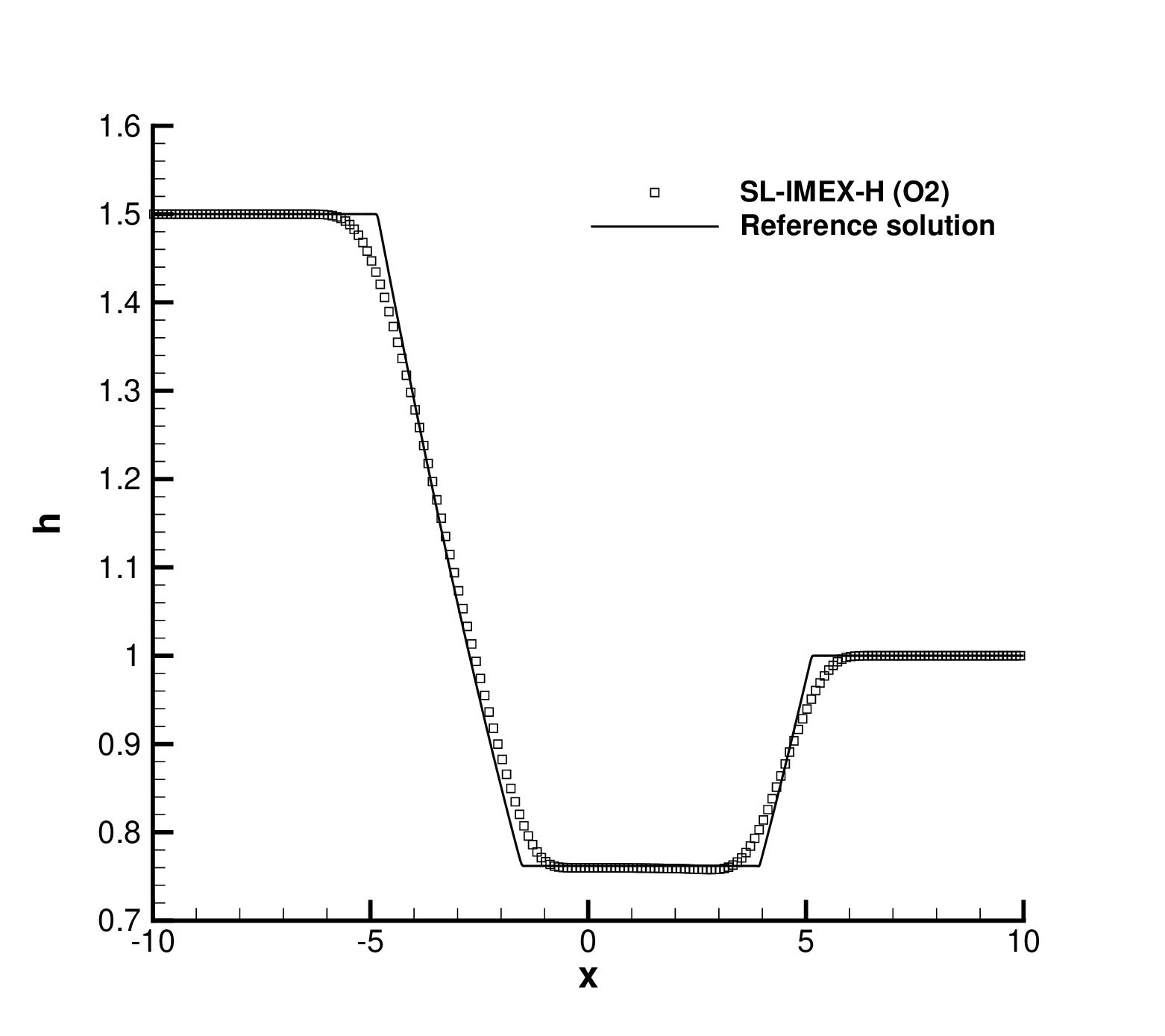}  &                 
				\includegraphics[width=0.47\textwidth]{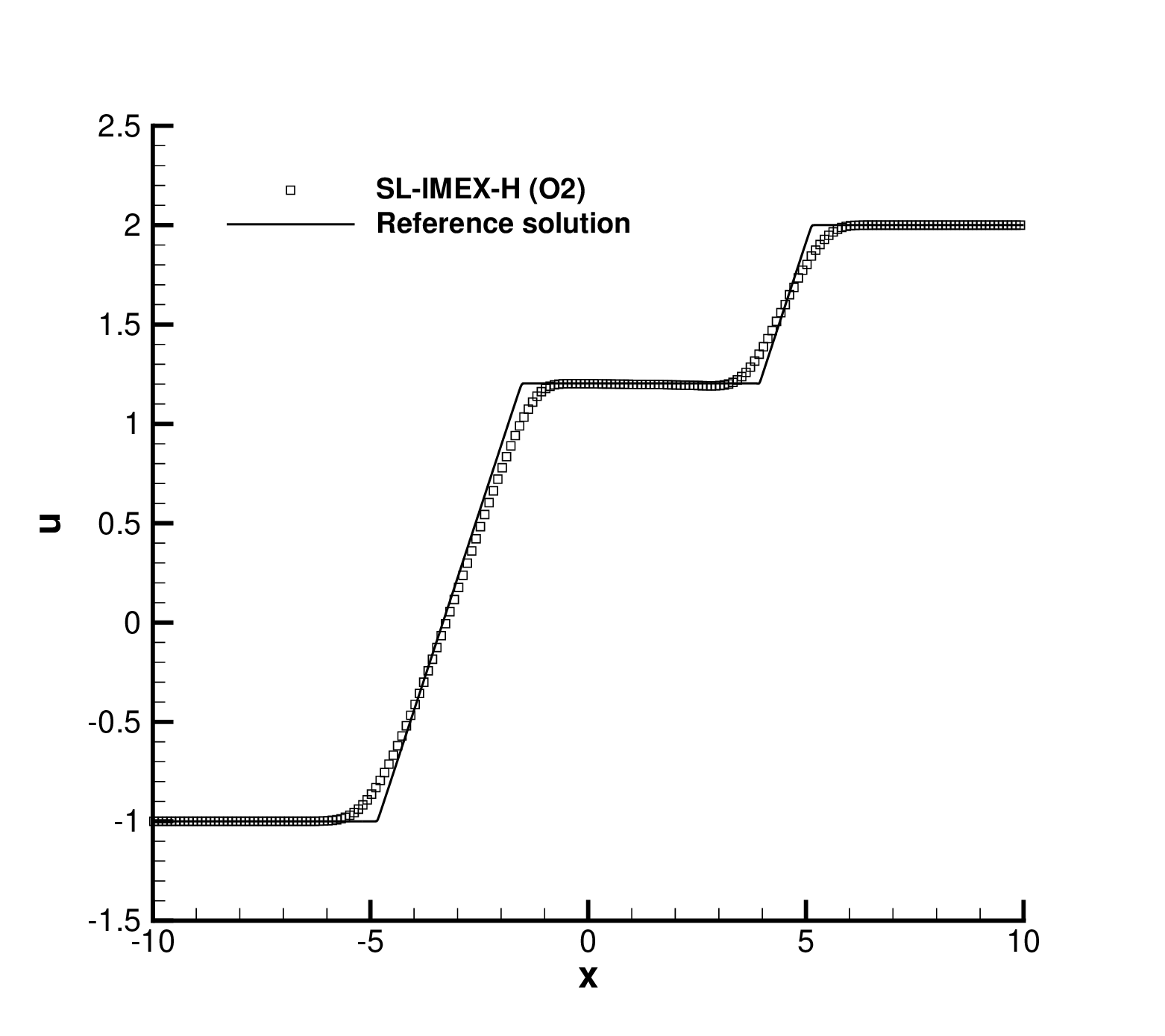} \\
				\includegraphics[width=0.47\textwidth]{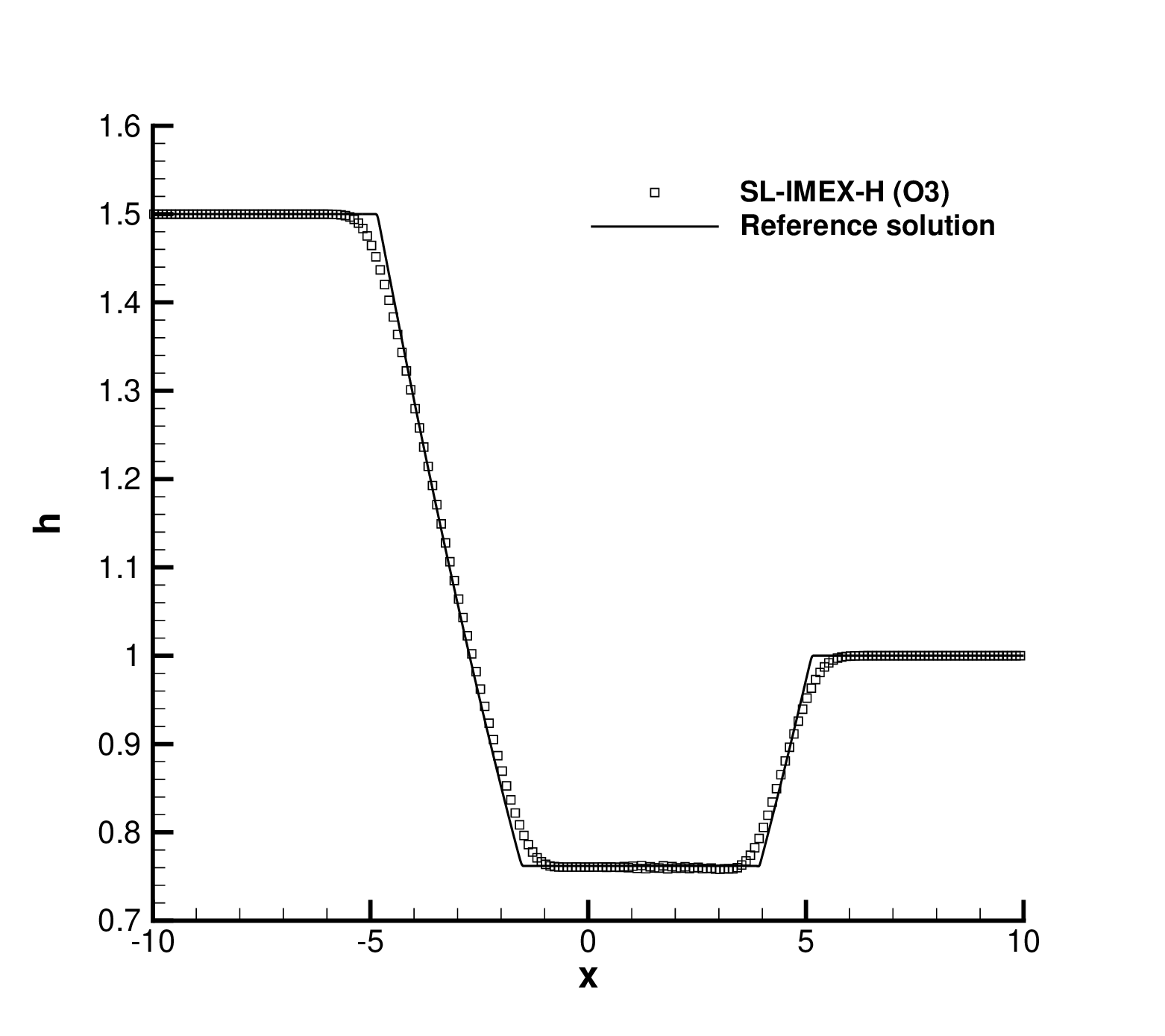}  &                 
				\includegraphics[width=0.47\textwidth]{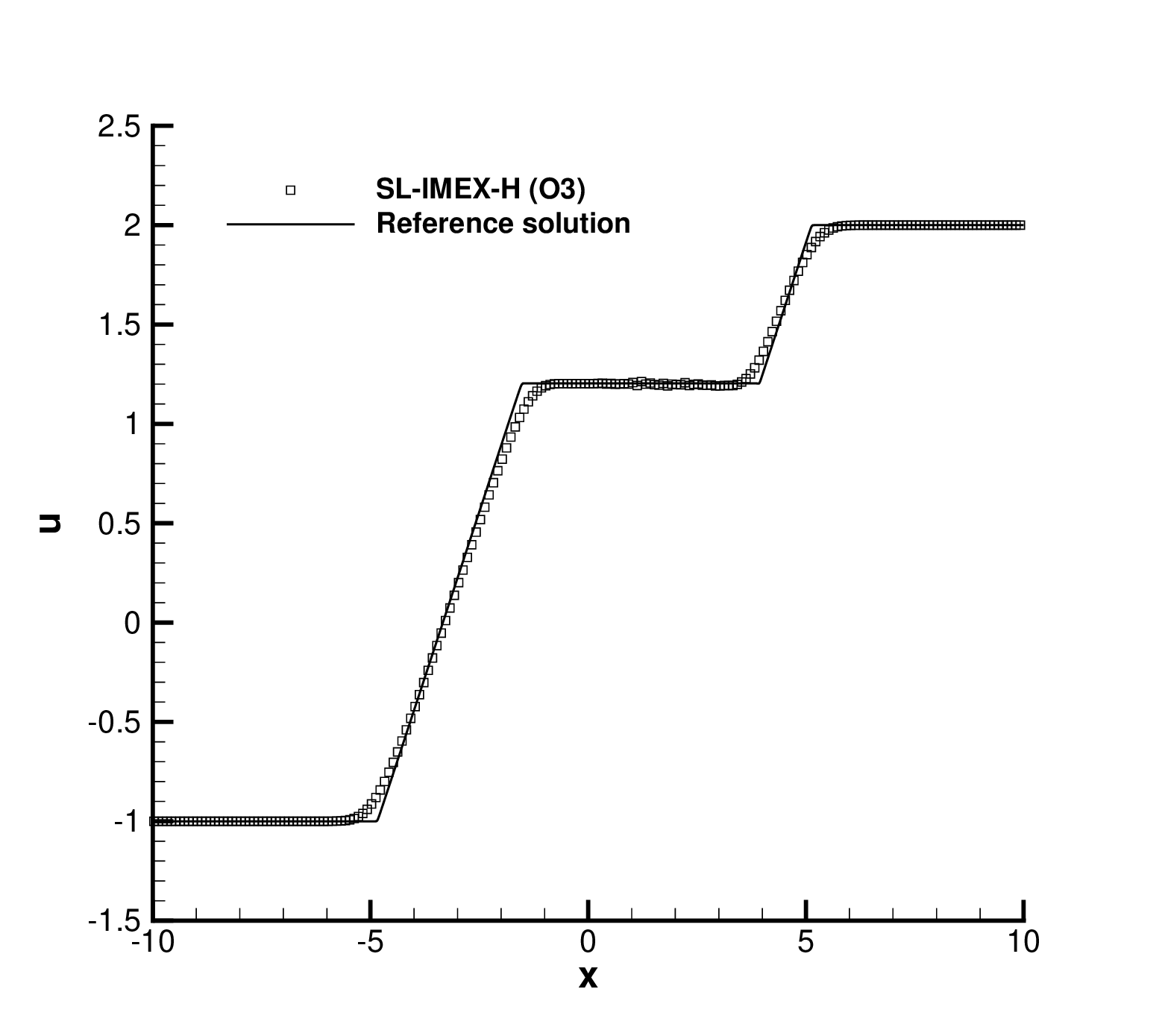} \\
			\end{tabular}
			\caption{Riemann problem RP1 at time $t=1$ with CFL=3. Second (top row) and third (bottom row) order numerical results compared against the reference solution for water height (left) and velocity (right).}
			\label{fig.RP10}
		\end{center}
	\end{figure}

	The first Riemann problem RP1 generates two non-symmetric rarefaction waves, hence this solution is continuous, but with discontinuous derivatives. The results obtained with the conservative SL-IMEX-H scheme run with CFL=3 are compared against the reference solution in Figure \ref{fig.RP10}, where second and third order accurate methods are used. Both numerical schemes generate profiles for the solution that are in excellent agreement with the exact solution and the third order method is slightly less dissipative, which can be appreciated looking at the heads and tails of the rarefaction waves. 
	
\begin{figure}[!htbp]
		\begin{center}
			\begin{tabular}{cc} 
				\includegraphics[width=0.47\textwidth]{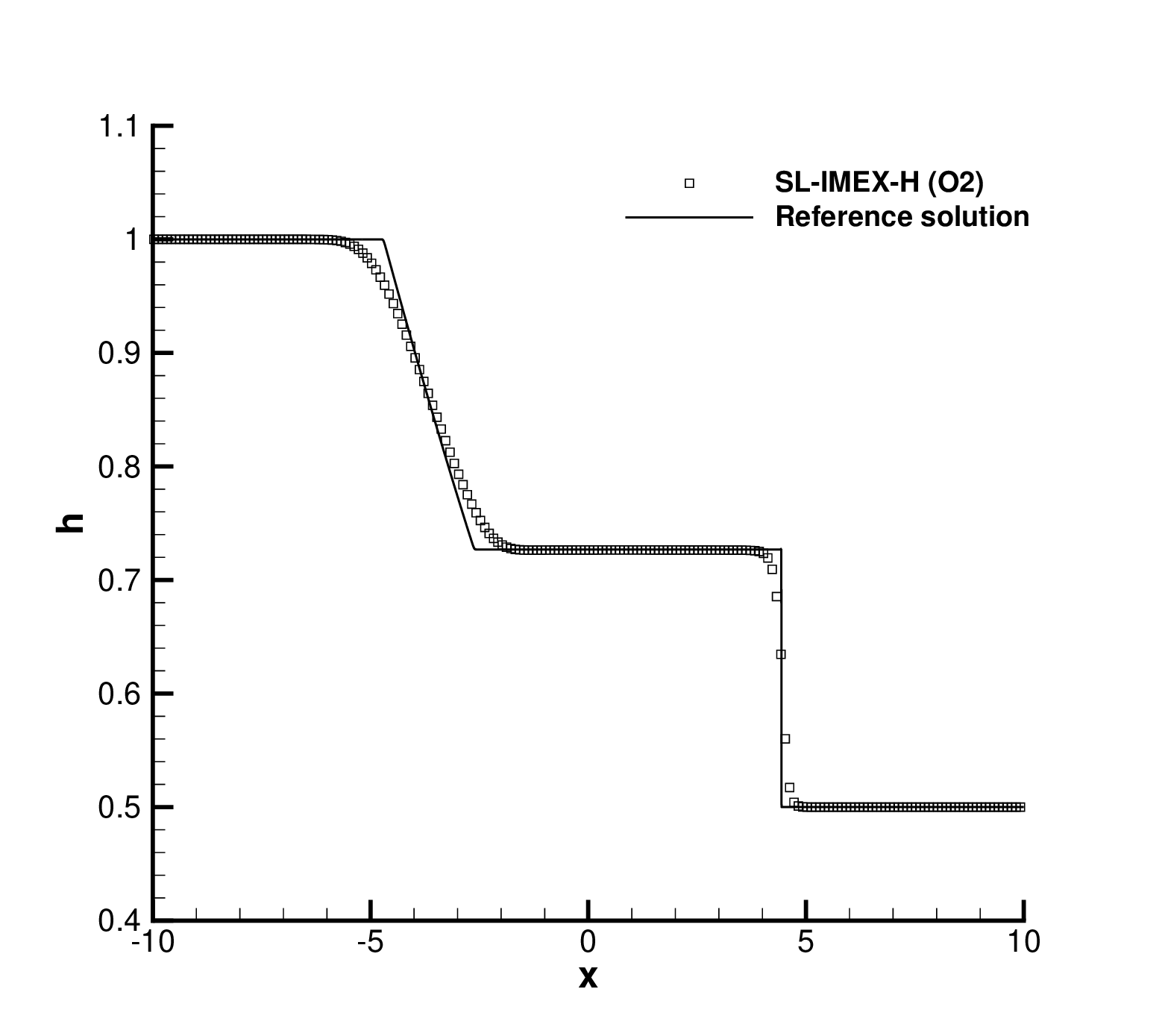}  &                 
				\includegraphics[width=0.47\textwidth]{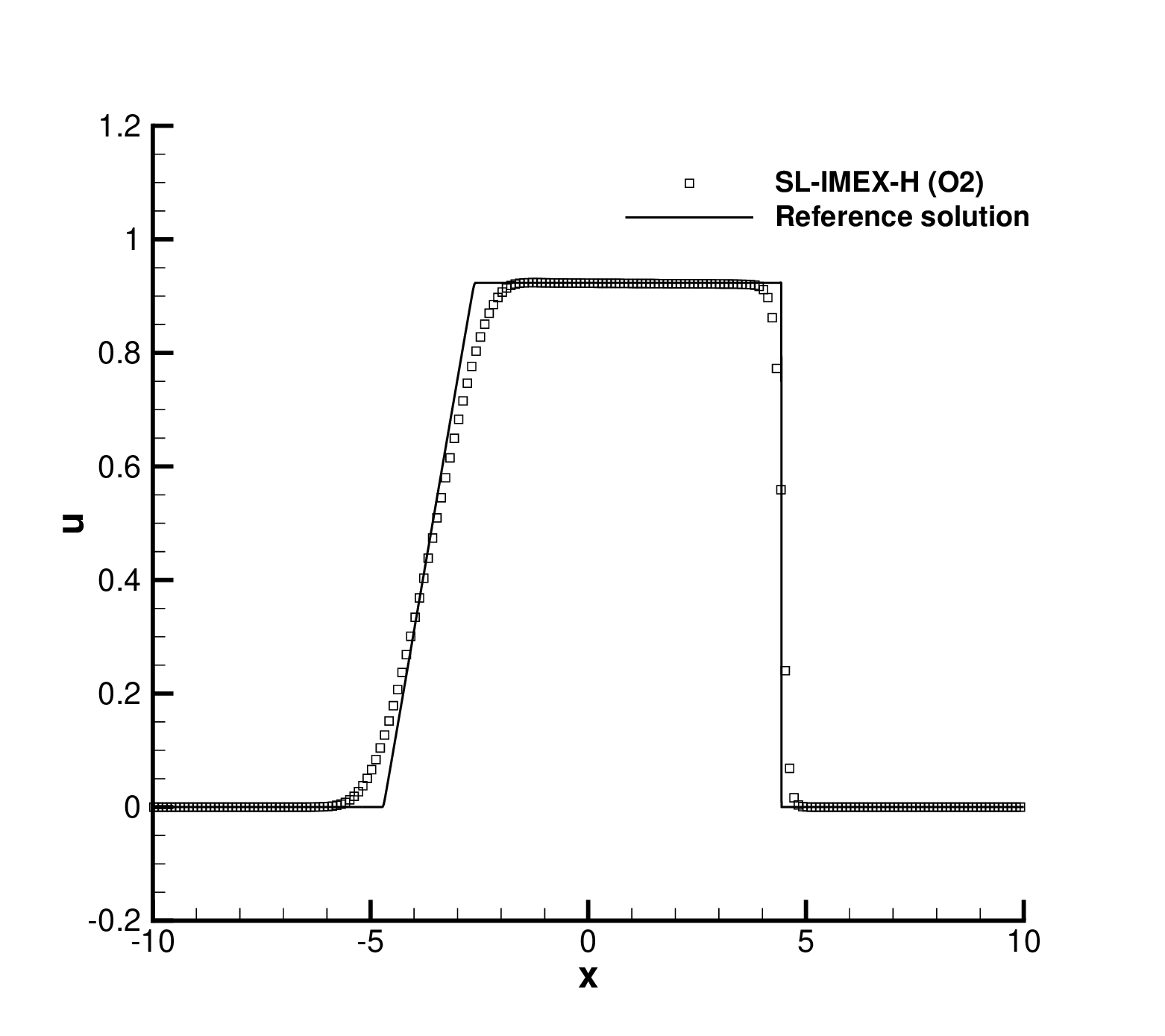} \\
				\includegraphics[width=0.47\textwidth]{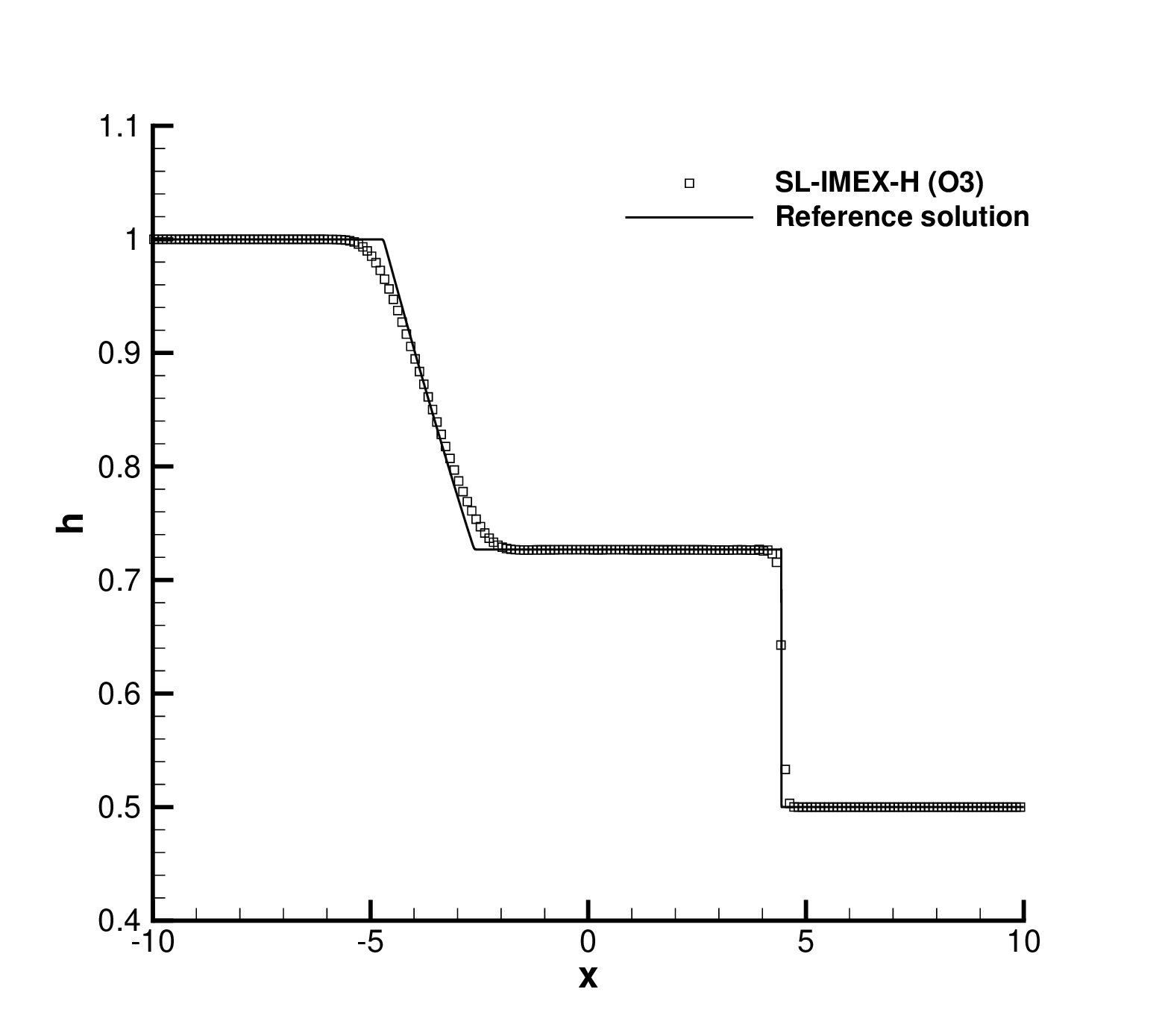}  &                 
				\includegraphics[width=0.47\textwidth]{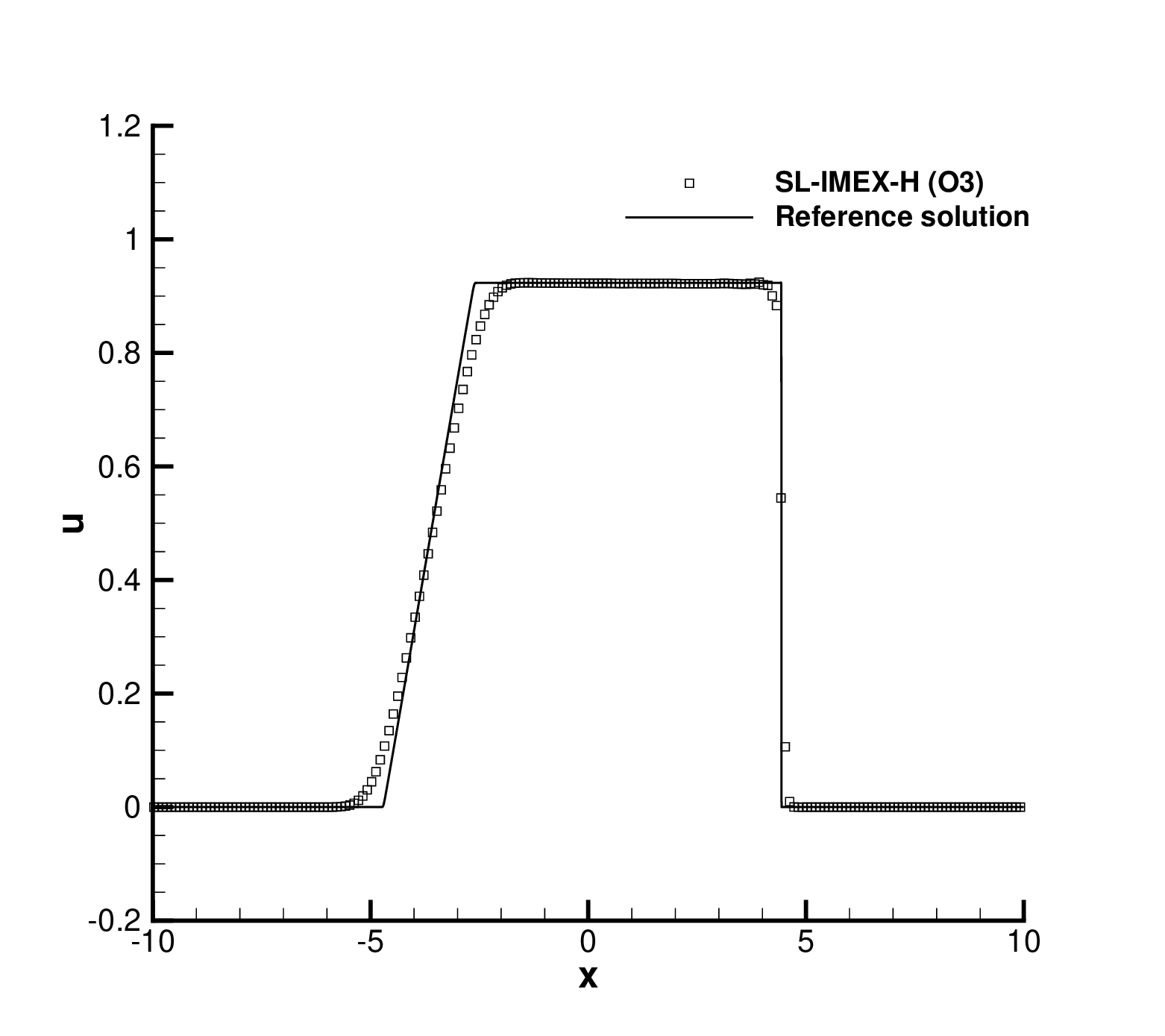} \\
			\end{tabular}
			\caption{Riemann problem RP2 at time $t=1.5$ with CFL=2. Second (top row) and third (bottom row) order numerical results compared against the reference solution for water height (left) and velocity (right).}
			\label{fig.RP12}
		\end{center}
	\end{figure}

\begin{figure}[!htbp]
		\begin{center}
			\begin{tabular}{cc} 
				\includegraphics[width=0.47\textwidth]{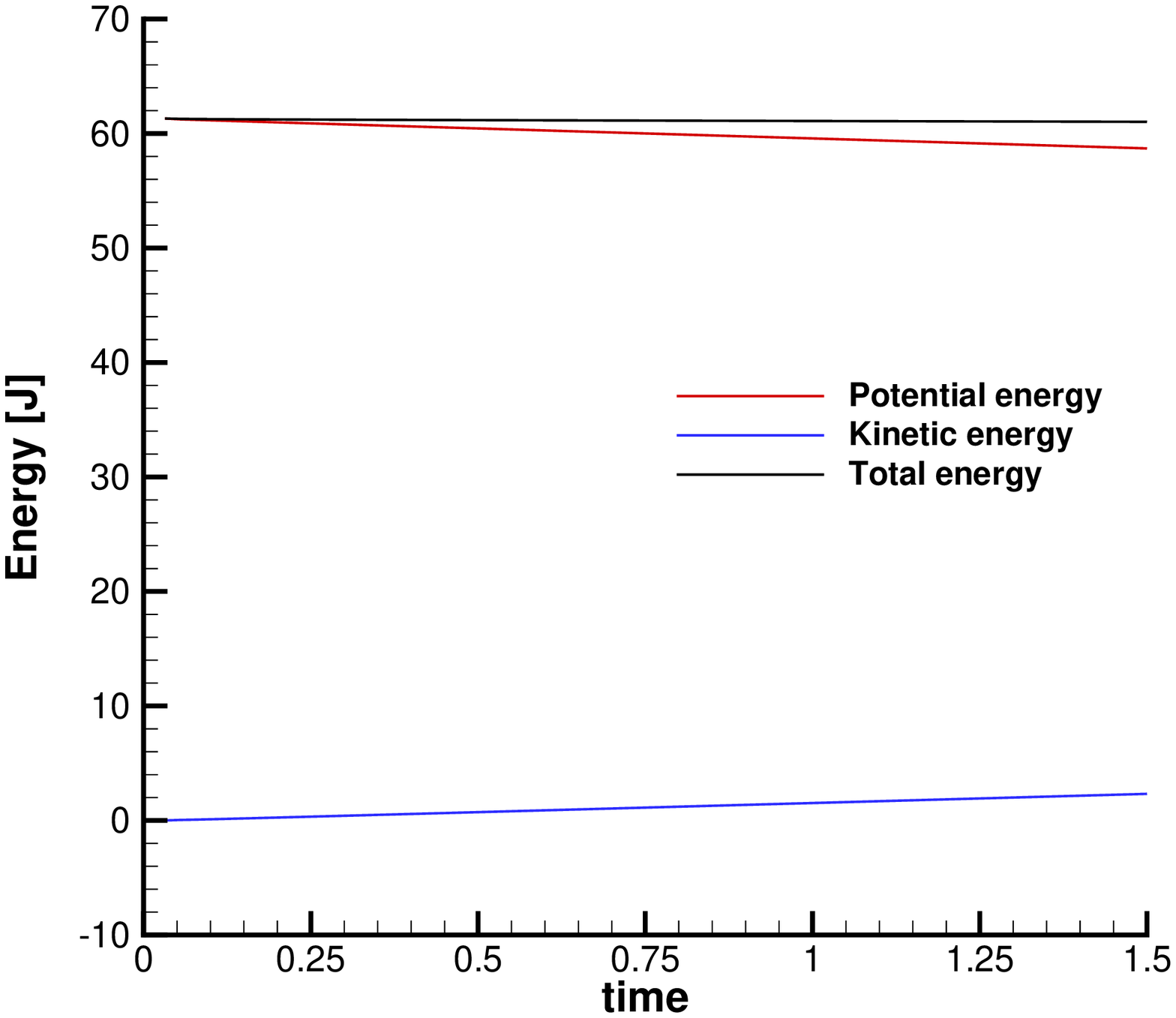}  &                 
				\includegraphics[width=0.47\textwidth]{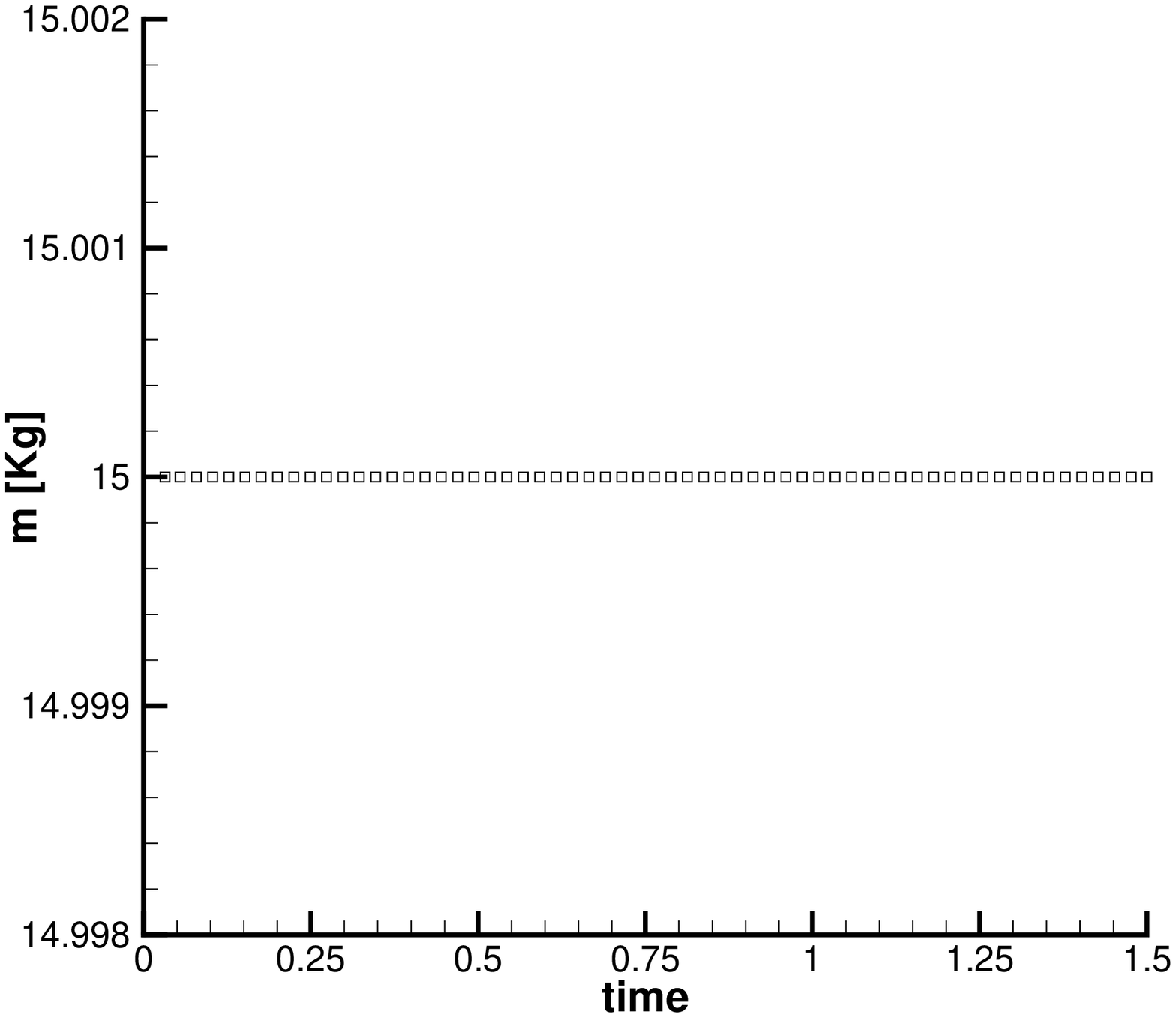} \\
			\end{tabular}
			\caption{Riemann problem RP2. Time evolution of potential, kinetic and total energy (left) and total mass (right) for the third order SL-IMEX-H scheme.}
			\label{fig.RP12_energy}
		\end{center}
	\end{figure}

	The second Riemann problem RP2 is directly inspired from the well-known Sod shock tube problem in gas dynamics \cite{ToroBook} and consists in one rarefaction wave and one shock wave traveling towards opposite directions. We use again both second and third order SL-IMEX-H schemes, and the results depicted in Figure \ref{fig.RP12} allows us to confirm that the conservative schemes correctly capture the location of the discontinuities as well as the values along the plateau exhibited by the exact solution. We underline that this is only possible if a conservative discretization of the momentum equation is adopted. Finally, Figure \ref{fig.RP12_energy} shows the time evolution of the total energy, namely potential and kinetic contributions, as well as the evolution of the total mass computed with \eqref{eqn.m}, that clearly remains constant for the entire simulation.
	
	\subsection{Relaxation system of SWE: stiff limit tests}
	This section is devoted to the validation of the novel SL-IMEX schemes applied to the relaxation system of SWE given by \eqref{sw1ra}-\eqref{sw1rb}. Three different test cases are presented in the stiff limit of the model, thus setting $\varepsilon=10^{-14}$ as relaxation parameter, so that the numerical results can be directly compared with the solution of the inviscid Burgers equation \eqref{eqn.Burgers}. As discussed in Section \ref{sec.AP}, the new methods are Asymptotic Preserving, therefore the proposed numerical schemes must retrieve a consistent discretization of the limit model, according to \eqref{eqn.hstiff} and \eqref{eqn.BurgersSL}, when the relaxation parameter tends to zero. For the sake of brevity, since the SL method is applied to the advection part of the governing PDE, we want to verify that the advection velocity of the limit system is correctly captured by the scheme, thus we only consider the stiffness in the parameter $\varepsilon$ for the relaxation system \eqref{sw1ra}-\eqref{sw1rb}, neglecting the stiffness related to the low Froude number. The SL-IMEX scheme \eqref{sw2raSL}-\eqref{sw2rbSL} is also tested in the limit $\varepsilon\to \infty$ in order to demonstrate the capability of correctly capturing the behavior of the original SWE system.
	
\begin{figure}[!htbp]
		\begin{center}
			\begin{tabular}{cc} 
				\includegraphics[width=0.47\textwidth]{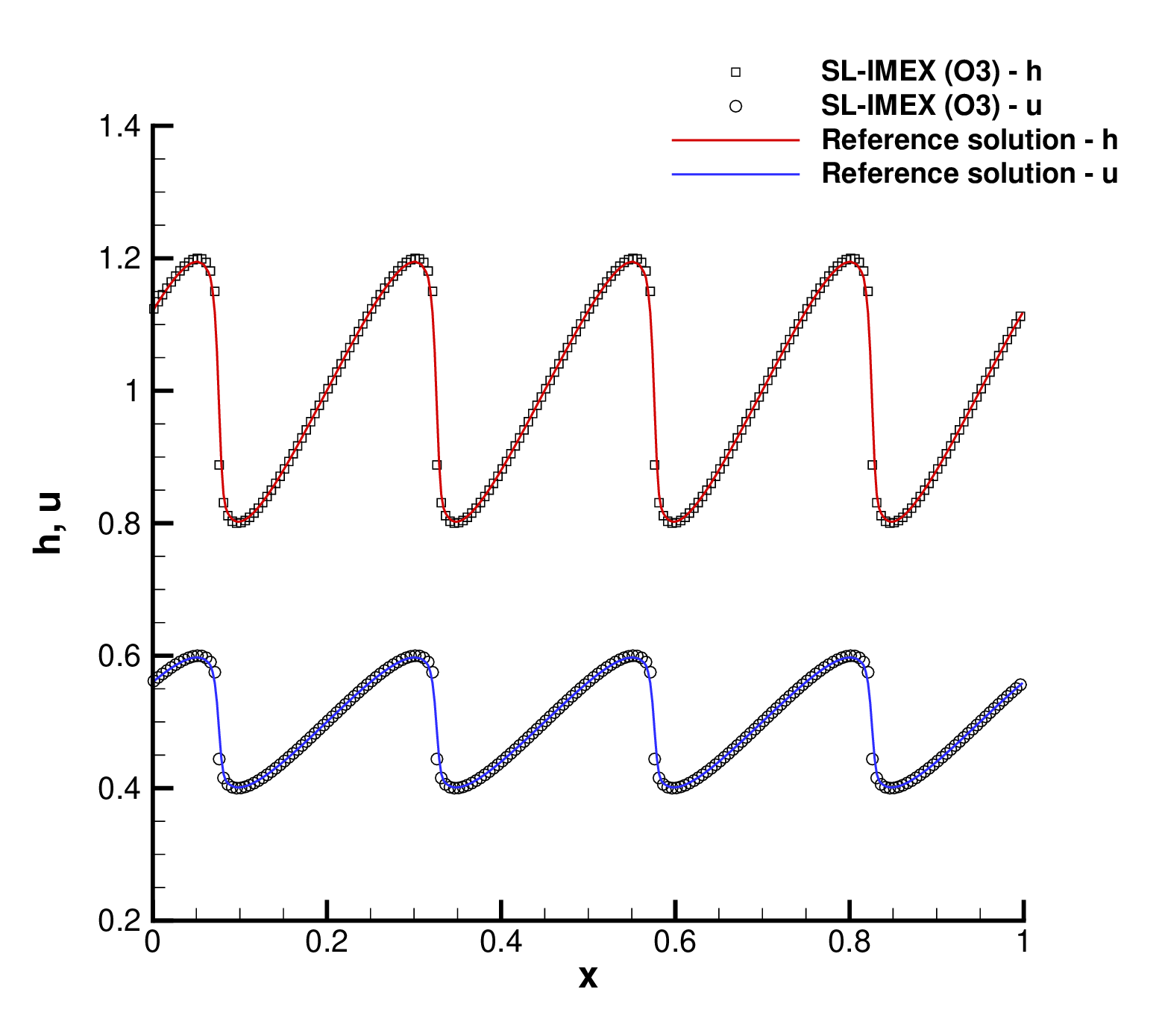}  &
				\includegraphics[width=0.47\textwidth]{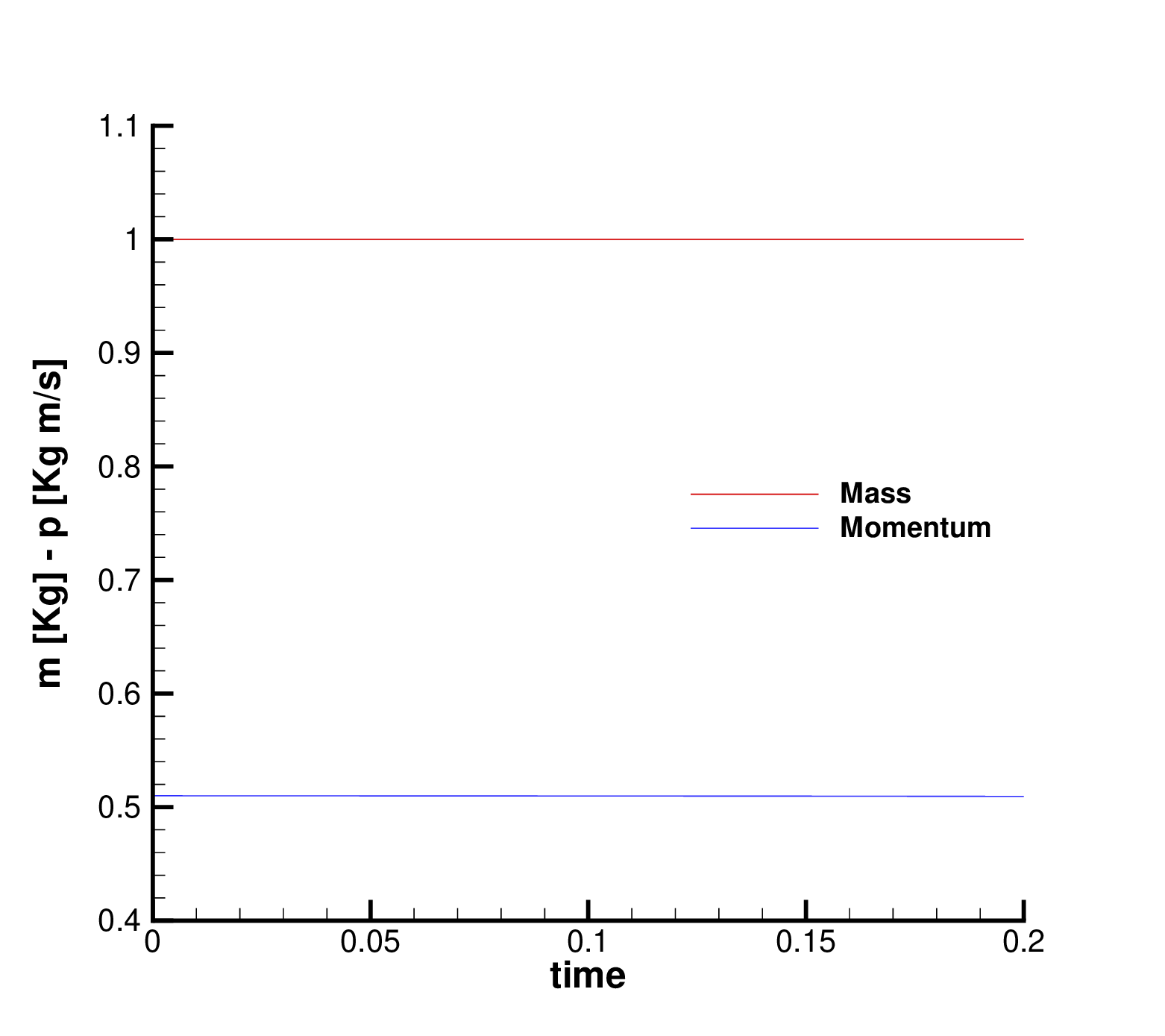} \\
			\end{tabular}
			\caption{Numerical convergence test with third order SL-IMEX scheme at time $t=0.2$ with $N_x=400$ for the relaxation system of SWE with $\varepsilon=10^{-14}$. Comparison against the reference solution for water height and velocity (left), and time evolution of total mass $m$ and momentum $p$ (right).}
			\label{fig.Burgers-convergence}
		\end{center}
	\end{figure}

	First, the numerical convergence in time is studied. To this aim, let us consider a computational domain $\Omega=[0;1]$ with periodic boundary conditions, which is discretized with a total number of cells $N_x=400$. The initial condition reads \cite{PR_IMEX}:
	\begin{equation}
	h(x,0) = 1 + 0.2 \sin(8\pi x), \qquad u(x,0)=\frac{h}{2},
	\end{equation}
	and the final time is set to $t_f=0.05$. Table \ref{tab.convB} reports the convergence rates for six successively refined values of the CFL number in \eqref{eqn.CFL}, showing that the formal order of accuracy is obtained for both second and third order SL-IMEX-H schemes. Next, the SL-IMEX method is used to run the simulation until the final time $t=0.2$, so that no shocks are present in the solution. Figure \ref{fig.Burgers-convergence} plots the numerical solution at the final time for both the water depth and the velocity, as well as the time evolution of the total mass and momentum which are fully preserved by the novel conservative SL methods.

	\begin{table}[!htbp]  
		\begin{center} 
			\begin{small}
				\renewcommand{\arraystretch}{1.2}
				\begin{tabular}{ccc|cc|cc} 
					\multicolumn{3}{c|}{} & \multicolumn{2}{c|}{SL-IMEX-H $\mathcal{O}2$} & \multicolumn{2}{c}{SL-IMEX-H $\mathcal{O}3$} \\				
					CFL & time steps & $\Delta t$ & $L_\infty$ & $\mathcal{O}(h)$ & $L_\infty$ & $\mathcal{O}(h)$ \\
					\hline
					8 &	11 &	4.5455E-03 &	9.7536E-05 &	-&	2.6662E-05 &	- \\
					7 &	13 &	3.8462E-03 &	7.2151E-05 &	1.80 &	1.7027E-05 &	2.68 \\
					6 &	15 &	3.3333E-03 &	5.3516E-05 &	2.09 &	1.0790E-05 &	3.19 \\
					5 &	18 &	2.7778E-03 &	3.8036E-05 &	1.87 &	6.3335E-06 &	2.92 \\
					4 &	22 &	2.2727E-03 &	2.4068E-05 &	2.28 &	3.1208E-06 &	3.53 \\
					3 &	30 &	1.6667E-03 &	1.3134E-05 &	1.95 &	1.1780E-06 &	3.14 \\
				\end{tabular}
			\end{small}
		\end{center}
		\caption{Numerical convergence results in time for the relaxation system of SWE with different CFL numbers using second and third order SL-IMEX-H schemes. The errors are measured in $L_\infty$ norm and refer to the water depth $h$ at time $t=0.05$.}  
		\label{tab.convB}
	\end{table}
	
	Second, two different Riemann problems are solved with the relaxation system of SWE. The initial setting is given in Table \ref{tab:initB} and the initial condition is assigned according to \eqref{eqn.RPini}. 
	
	\begin{table}[!htbp]  
		\caption{Initialization of Riemann problems for the relaxation system of SWE. Initial states left (L) and right (R) are reported as well as the final time of the simulation $t_f$, the computational domain $[x_L;x_R]$ and the position of the initial discontinuity $x_d$.}  
		\begin{center} 
			\begin{small}
				\renewcommand{\arraystretch}{1.0}
				\numerikNine
				\begin{tabular}{l|c|ccc|cc|cc} 
					\hline
					Name & $t_{f}$ & $x_L$ & $x_R$ & $x_d$ & $h_L$ & $u_L$ & $h_R$ & $u_R$ \\
					\hline
					B1 & 0.3 & -1.0 & 1.0 & 0.0 & 1.0 &  0.0  & 2.0 & 0.0 \\
					B2 & 0.4 & -1.0 & 1.0 & 0.0 & 2.0 &  0.0  & 1.0 & 0.0 \\
					\hline
				\end{tabular}
			\end{small}
		\end{center}
		\label{tab:initB}
	\end{table}
	
\begin{figure}[!htbp]
		\begin{center}
			\begin{tabular}{cc} 
				\includegraphics[width=0.47\textwidth]{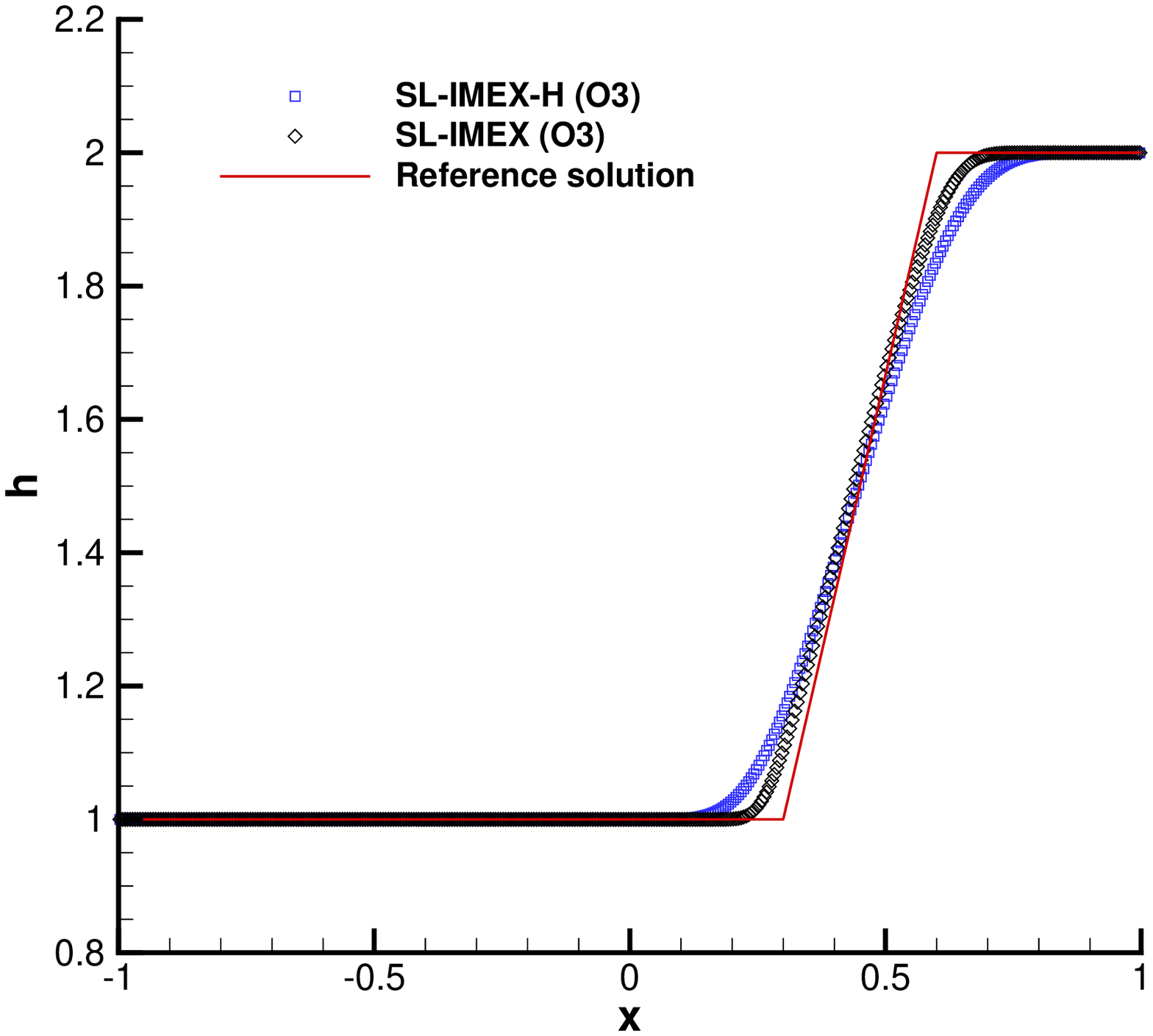}  &
				\includegraphics[width=0.47\textwidth]{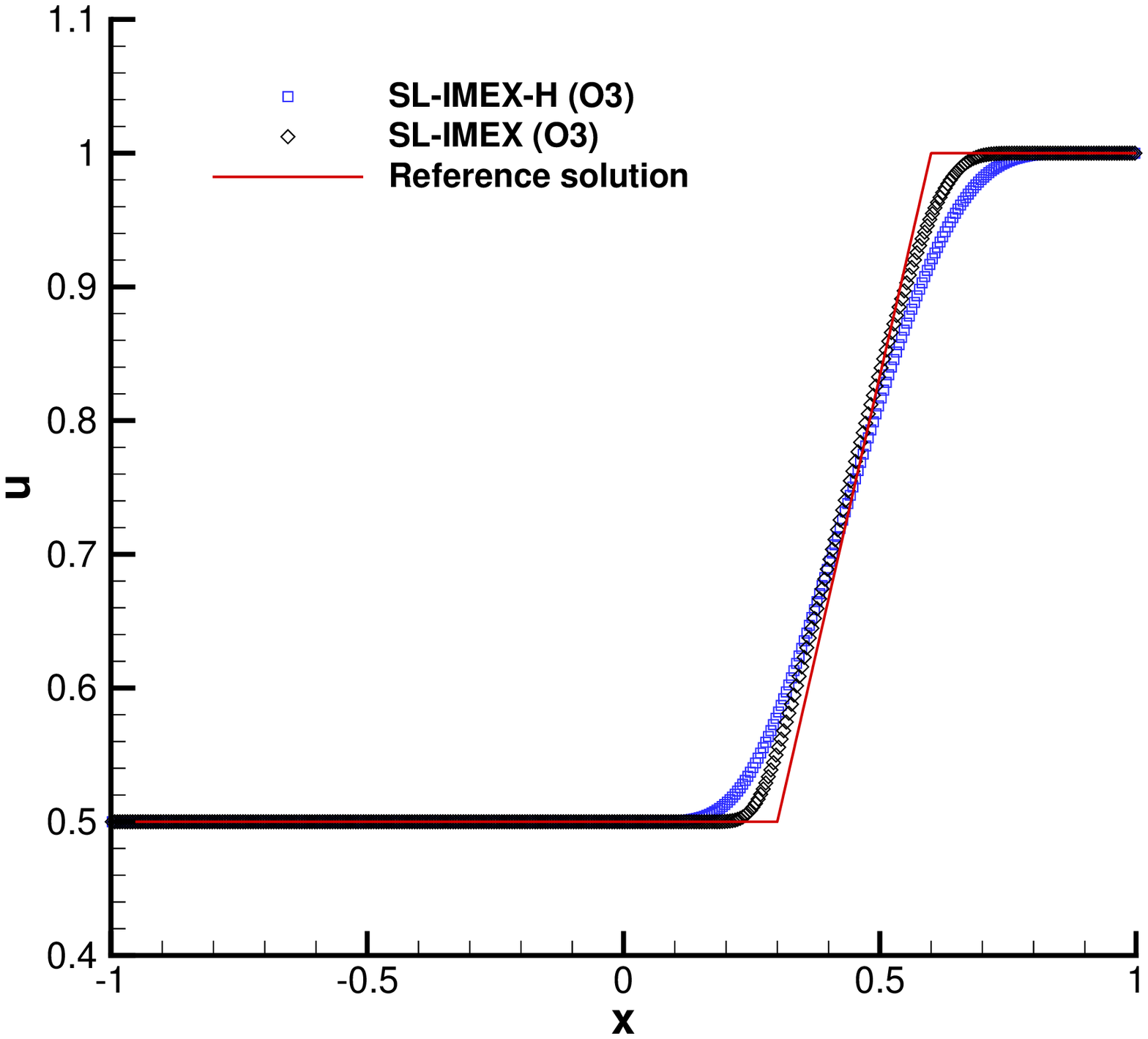} \\
				\includegraphics[width=0.47\textwidth]{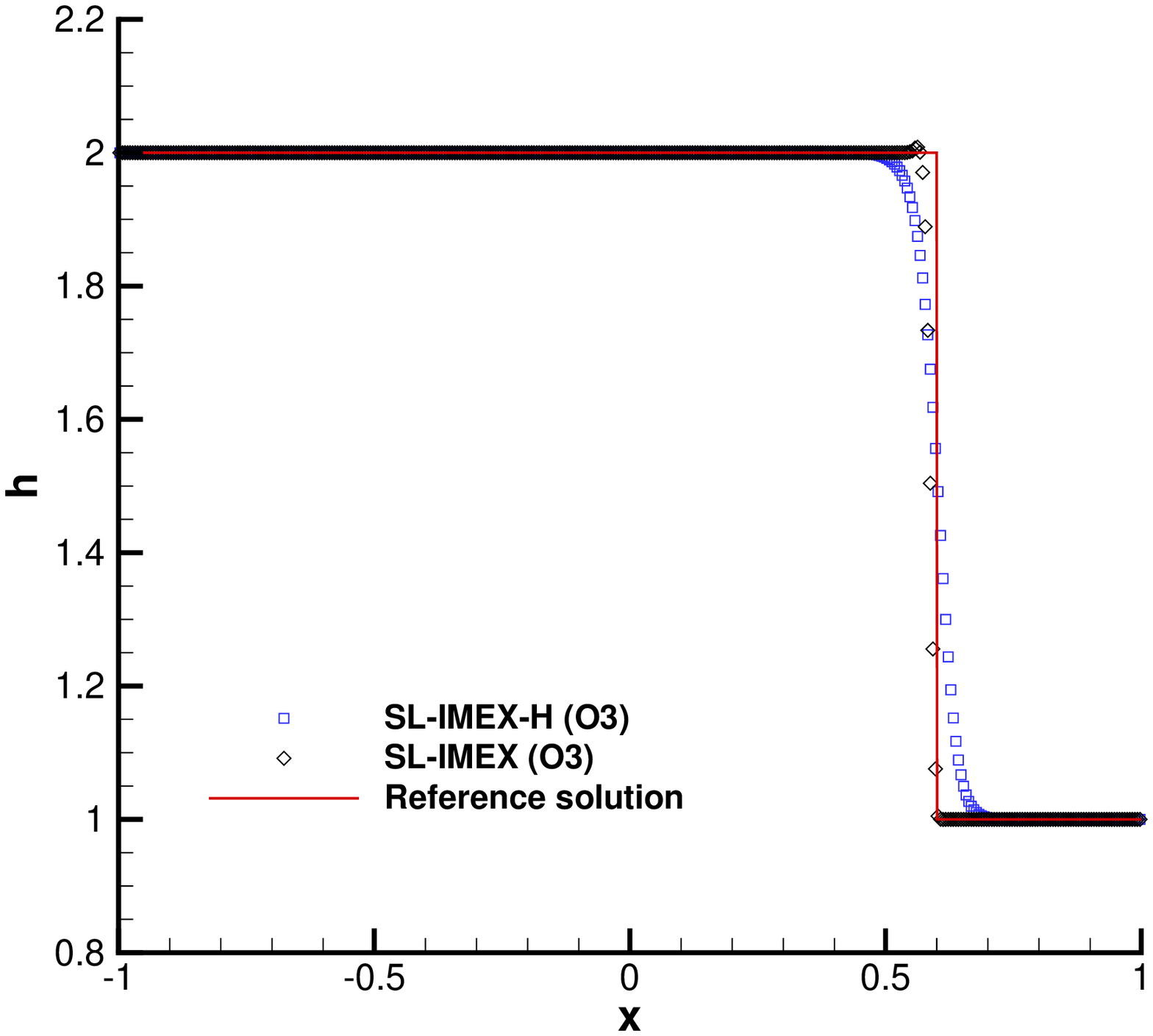}  &
				\includegraphics[width=0.47\textwidth]{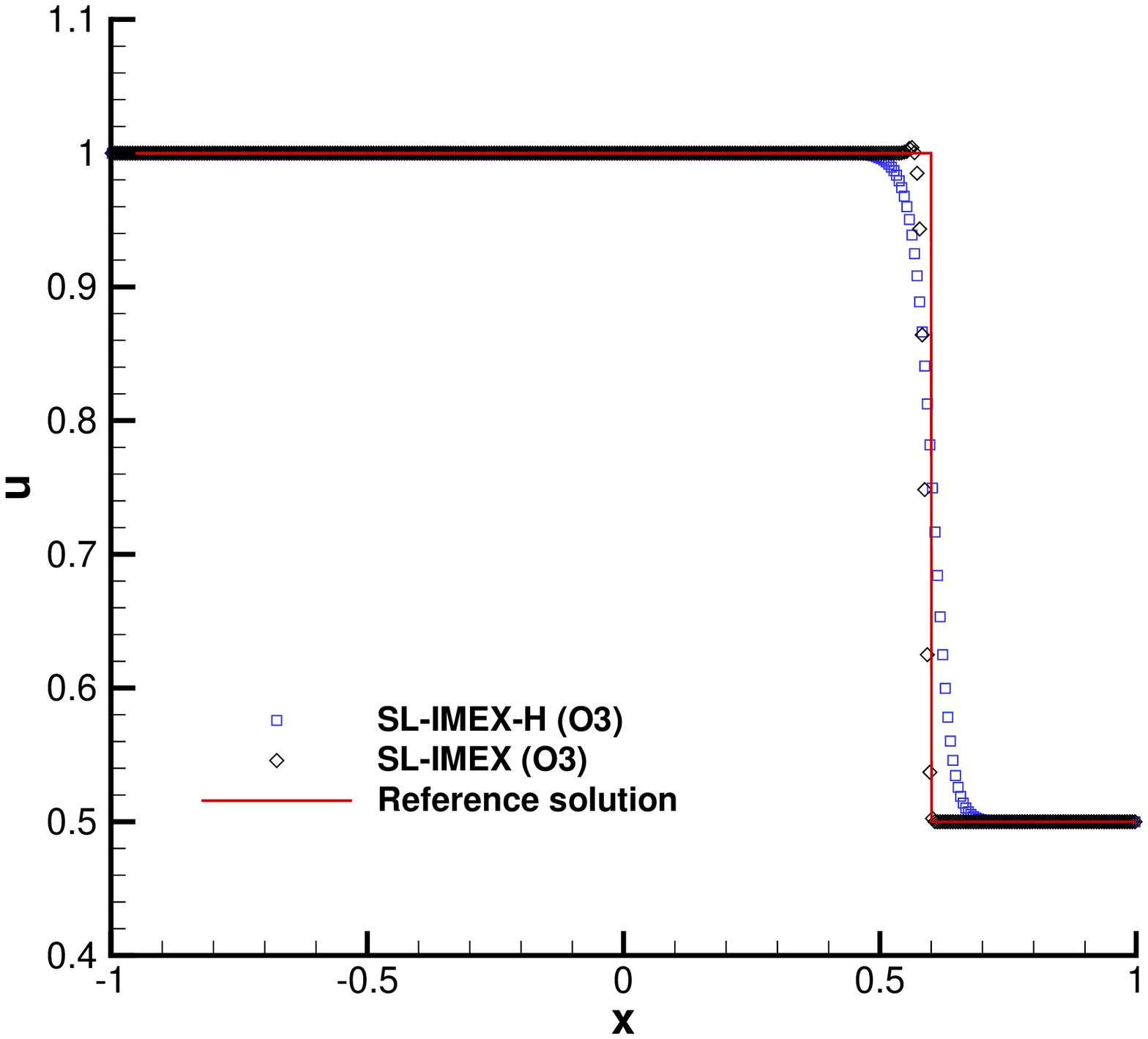} \\			
			\end{tabular}
			\caption{Riemann problems with third order SL-IMEX scheme with $N_x=400$ for the relaxation system of SWE with $\varepsilon=10^{-14}$. Top: rarefaction solution (B1 test) at time $t=0.3$. Bottom: shock solution (B2 test) at time $t=0.04$.  Comparison against the reference solution for water height (left) and velocity (right).}
			\label{fig.Burgers-RP}
		\end{center}
	\end{figure}

	Riemann problem B1 involves a rarefaction wave, while the test case B2 is concerned with a shock wave. Conservation in the numerical scheme is then crucial in order to correctly capture the location of the discontinuities, thus both B1 and B2 represent a validation of both SL-IMEX-H and SL-IMEX schemes in the stiff limit of the relaxation system \eqref{sw1ra}-\eqref{sw1rb}. The computational domain is defined as $\Omega=[-1;1]$ and counts a total number of control volumes $N_x=400$. The CFL number is set to $\textnormal{CFL}=4$ and a comparison between the third order numerical results and the exact solution is depicted in Figure \ref{fig.Burgers-RP}. Even in this case, both SL-IMEX-H and SL-IMEX schemes are able to capture the correct wave speeds while running the simulation with a time step which is four times larger than the one of an explicit scheme. Furthermore, the results are proven to be consistent with the limit model, i.e. the inviscid Burgers equation. Notice that also the velocity relaxes to the correct asymptotic limit, that is $u \to h/2$, for all the test cases shown in this section. The SL-IMEX scheme is clearly less dissipative than the SL-IMEX-H version, which is visible from the resolution of the tail and head of the rarefaction wave in test B1, as well as in the sharp resolution of the shock front for problem B2.
	
\begin{figure}[!htbp]
		\begin{center}
			\begin{tabular}{cc} 
				\includegraphics[width=0.47\textwidth]{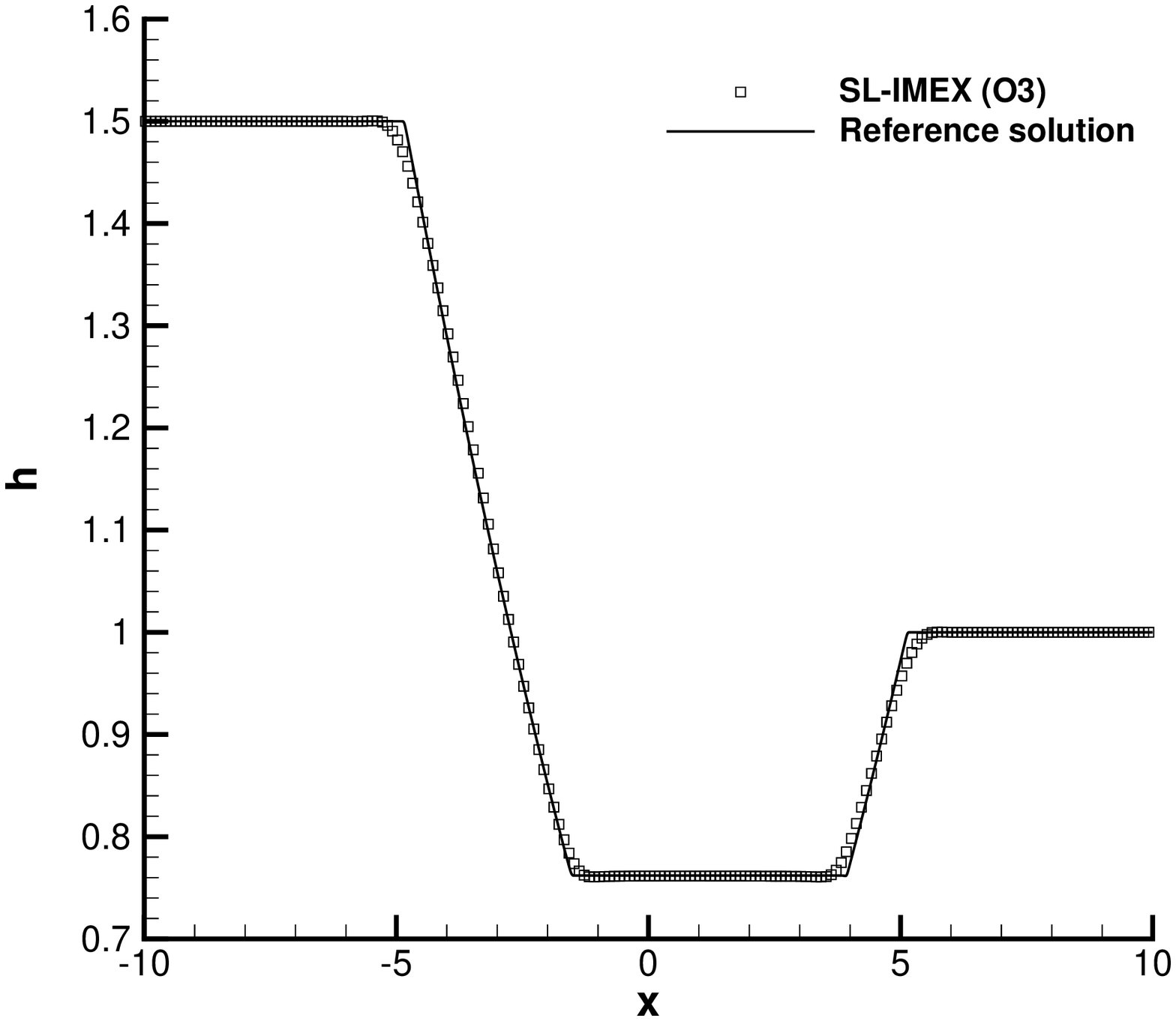}  &
				\includegraphics[width=0.47\textwidth]{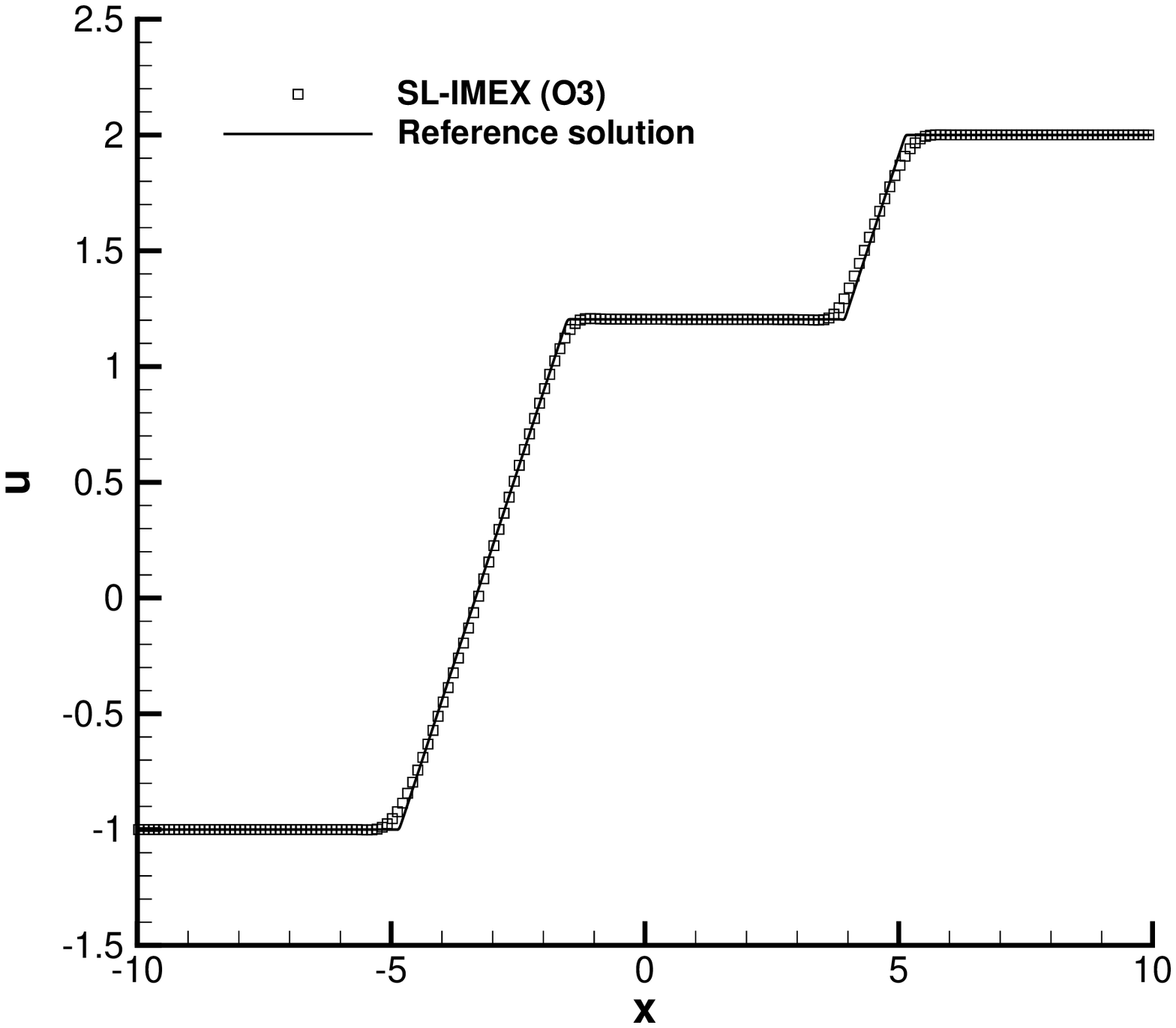} \\
				\includegraphics[width=0.47\textwidth]{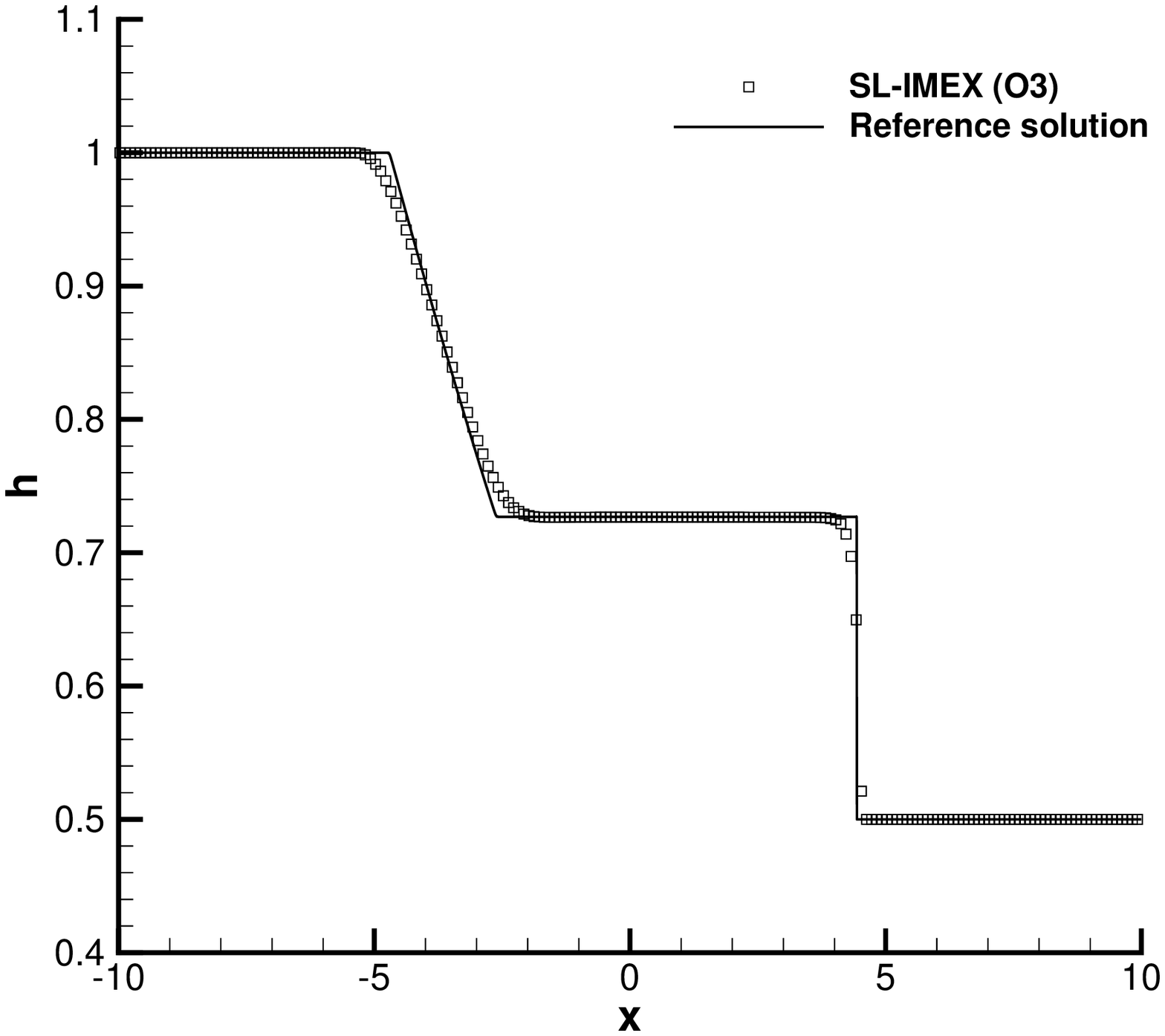}  &
				\includegraphics[width=0.47\textwidth]{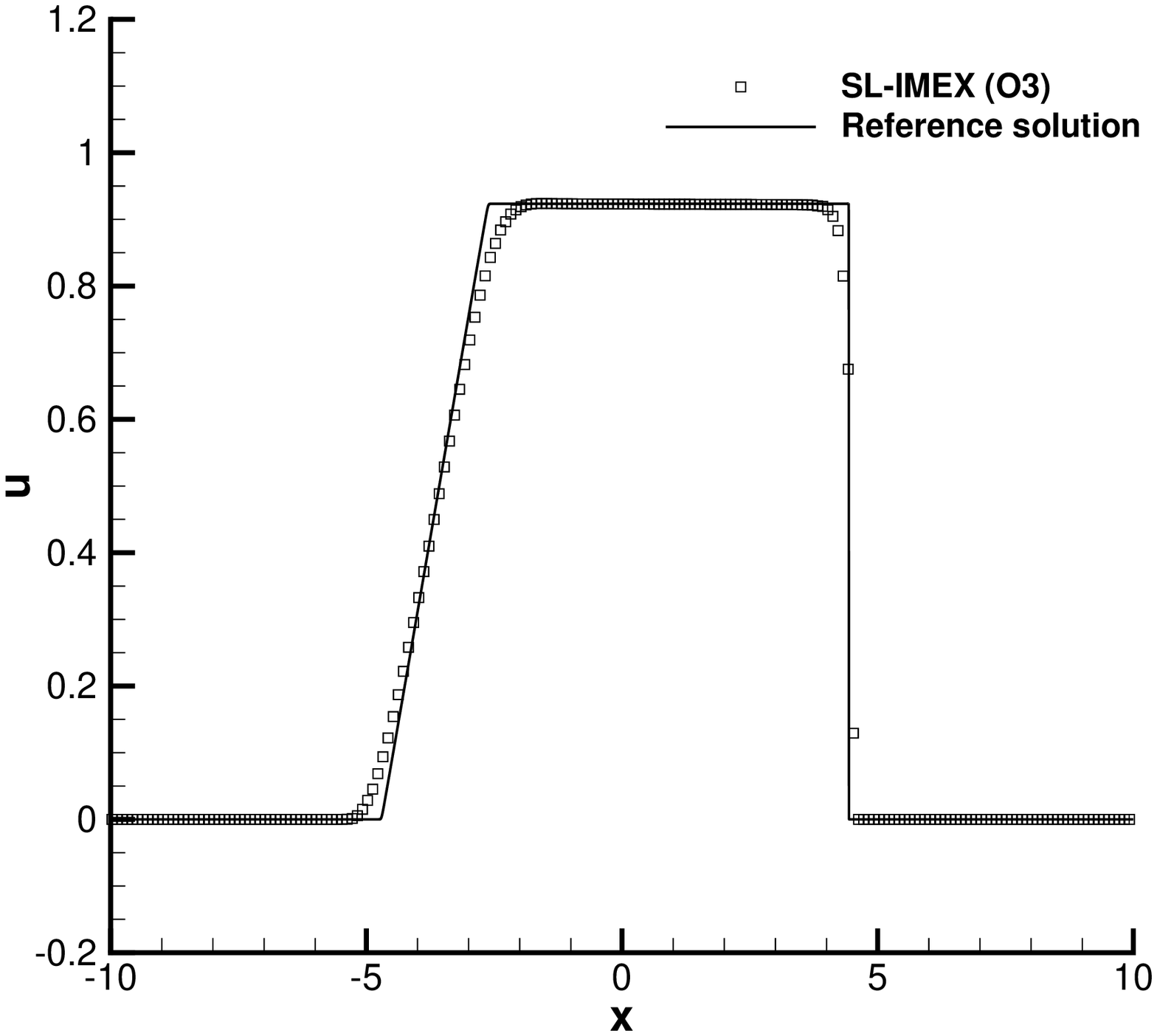} \\			
			\end{tabular}
			\caption{Riemann problems RP1 (top) and RP2 (bottom) with third order SL-IMEX scheme with $N_x=400$ for the relaxation system of SWE with $\varepsilon=10^{-14}$. Comparison against the reference solution for water height (left) and velocity (right).}
			\label{fig.RP_fullSL}
		\end{center}
	\end{figure}

\begin{figure}[!htbp]
	\begin{center}
		\begin{tabular}{cc} 
			\includegraphics[width=0.47\textwidth]{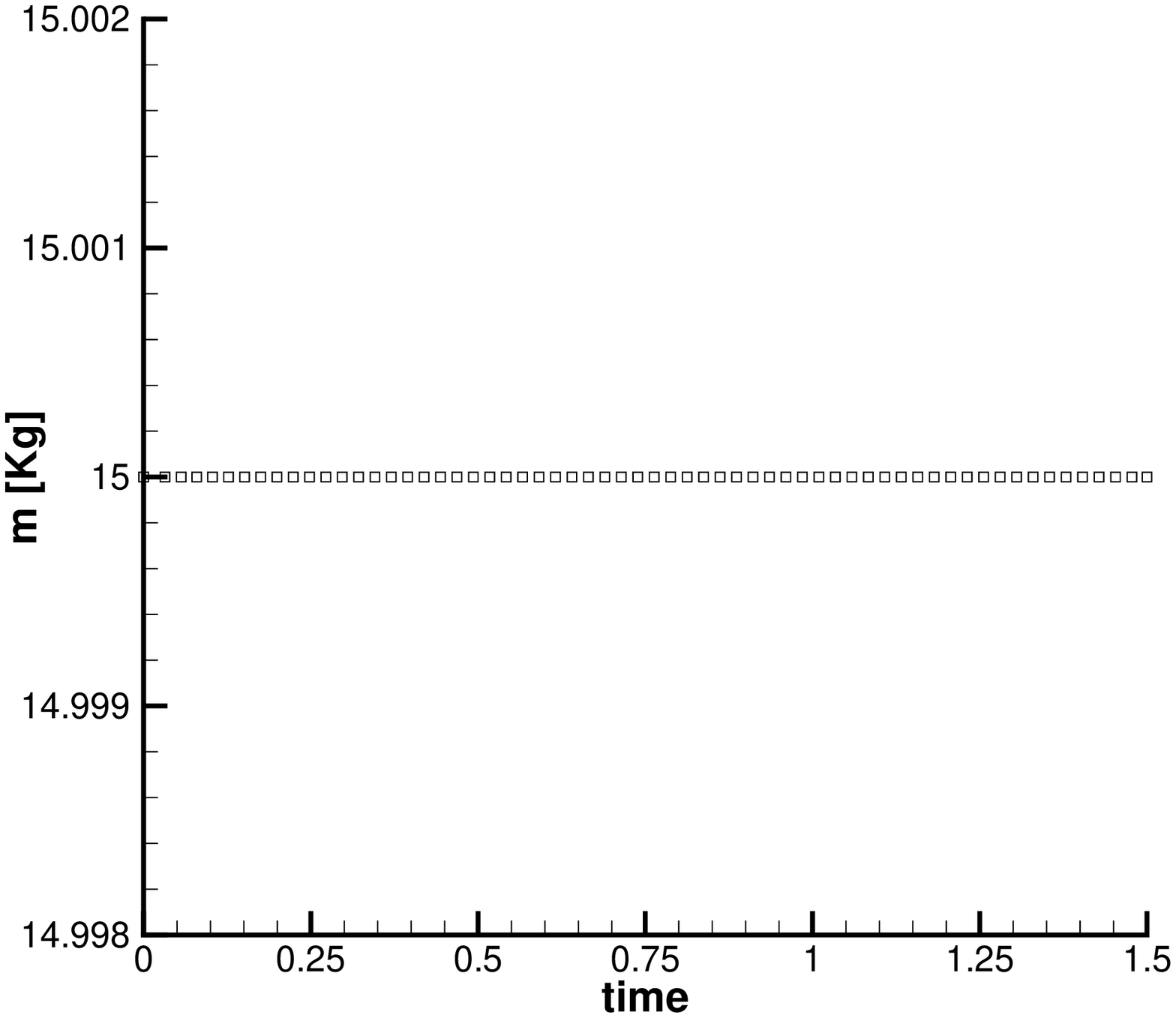}  &
			\includegraphics[width=0.47\textwidth]{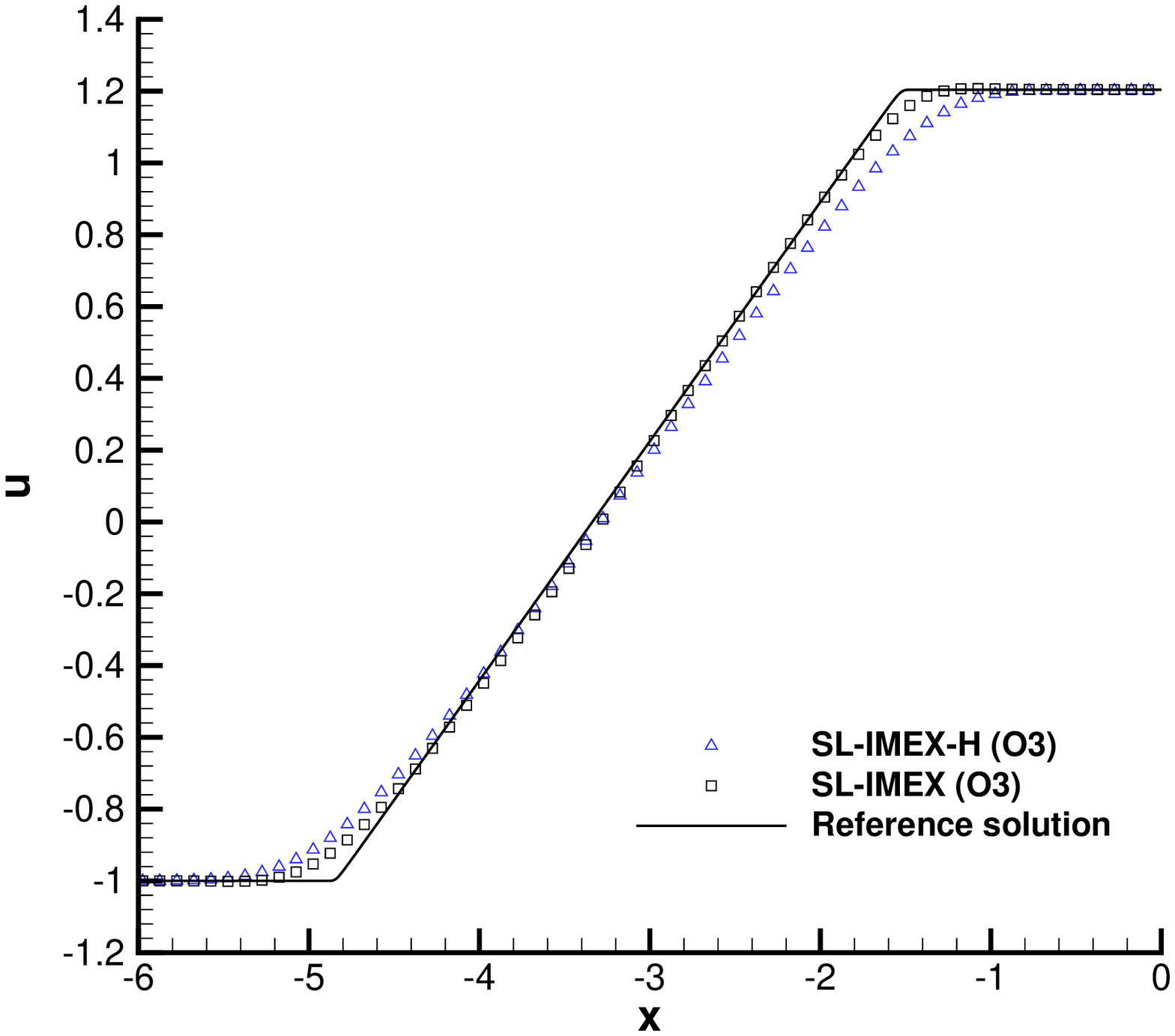} \\
		\end{tabular}
		\caption{Left: time evolution of the total mass $m$ for RP2 with third order SL-IMEX scheme. Right: comparison against the reference solution for RP1 with third order SL-IMEX (black squares) and SL-IMEX-H (blue triangles) scheme with $N_x=400$ for the relaxation system of SWE with $\varepsilon=10^{-14}$. Zoom on the velocity profile across the tail and head of the left rarefaction wave.}
		\label{fig.zoom_fullSL}
	\end{center}
\end{figure}

	Finally, the SL-IMEX scheme \eqref{sw2raSL}-\eqref{sw2rbSL} is used to run the Riemann problems RP1 and RP2 in the case $\epsilon \to \infty$. The initial condition is reported in Table \ref{tab:init} and the third order results are depicted in Figure \ref{fig.RP_fullSL}, where an excellent agreement with the reference solution can be appreciated. Moreover, mass is fully conserved and even for these test cases the SL-IMEX schemes are less dissipative than the SL-IMEX-H methods, as shown in Figure \ref{fig.zoom_fullSL} with a zoom on the tail and head of the rarefaction wave for RP1.
	
	\section{Conclusions} \label{sec.concl}
	In this work a novel class of semi-Lagrangian schemes with IMEX Runge-Kutta time stepping has been derived and discussed. First, the method is presented for a simple advection-diffusion equation, where the nonlinear convective terms have been discretized explicitly using a high order semi-Lagrangian technique, while the diffusion part of the equation has been treated implicitly. The resulting scheme is therefore unconditionally stable and a suite of test cases has been shown to demonstrate the achievement of the formal order of accuracy. Second, the SL-IMEX methods have been extended to ensure conservation of the transported quantities in the case of hyperbolic systems of conservation laws by means of a novel technique based on the integration of the governing PDE onto the space-time control volume generated by the motion of the grid points along the characteristics. As a prototype model the shallow water system has been used.
	
	The space discretization relies on high order CWENO reconstruction operators, while a Cauchy-Kovalevskaya procedure allows the semi-Lagrangian scheme to reach high order of accuracy in time as well. The conservative version of the schemes has then be applied to test problems involving shocks, rarefaction waves and contact discontinuities. The conservation property exhibited by the schemes is crucial to correctly capture and locate shock waves and plateau of the solution. The IMEX Runge-Kutta method is used for the time stepping, with explicit convective contribution and implicit pressure terms, so that the CFL-type stability condition can be completely relaxed. Finally, the asymptotic preserving (AP) property of the schemes has been studied as well considering the shallow water systems with relaxation terms, for which two different AP schemes are derived, both in terms of a semi-Lagrangian IMEX discretization. Numerical experiments demonstrate that the novel SL-IMEX methods provide a consistent discretization  in the stiff regime of the limit model given by the inviscid Burgers equation.
	
	Future research will concern the extension of the novel class of SL-IMEX schemes to multiple space dimensions, including the usage of unstructured meshes. Complex hyperbolic PDE systems like the magnetohydrodynamics equations (MHD) or the GPR unified model for continuum mechanics \cite{PeshRom2014} are also planned to be investigated. Finally, a space-time predictor strategy along the line of \cite{PNPM,ADERFSE} might also be useful to overcome the Cauchy-Kovalevskaya procedure, thus constituting another potential future topic of research.

	\begin{acknowledgements}
		WB and LP would like to thank the Italian Ministry of Instruction, University and Research (MIUR) to support this research with funds coming from PRIN Project 2017 (No. 2017KKJP4X entitled “Innovative numerical methods for evolutionary partial differential equations and applications”). 
		Part of this work was also supported by the projects "Advanced Computations and Experiments for anisotropic particle transport in turbulent flows" (ACE) (PRIN Project 2017 No. 2017RSH3JY) and CRC project HM: Hydropeaking mitigation, financed by the Free University of Bozen-Bolzano.
	\end{acknowledgements}
	
	\newpage
	\appendix
	\section{IMEX schemes} \label{app.IMEX}
	The Butcher tableau \eqref{eqn.butcher} for the IMEX schemes used in this work are reported hereafter. They have been derived in \cite{PR_IMEX,PR_IMEXHO} and each IMEX scheme is described with a triplet $(s,\tilde{s},p)$ which characterizes the number $s$ of stages of the implicit method, the number $\tilde{s}$ of stages of the explicit method and the order $p$ of the resulting scheme. The acronym SA stands for Stiffly Accurate, which is a crucial feature in the stiff limit of the governing PDE in order to be consistent with the limit model at the discrete level \cite{Dellacherie1,BosRus2019}. Strong Stability Preserving (SSP) methods are preferred when dealing with shock waves and other discontinuities.
	
	\begin{itemize}
		\item SP(1,1,1)
		
		\begin{equation}
			\begin{array}{c|c}
				0 & 0 \\ \hline & 1
			\end{array} \qquad
			\begin{array}{c|c}
				1 & 1 \\ \hline & 1
			\end{array}
			\label{eqn.IMEX1}
		\end{equation}
		
		\item SA SSP(3,3,2)
		
		\begin{equation}
			\begin{array}{c|ccc}
				0 & 0 & 0 & 0 \\ 1/2 & 1/2 & 0 & 0 \\ 1 & 1/2 & 1/2 & 0 \\ \hline & 1/3 & 1/3 & 1/3
			\end{array} \qquad
			\begin{array}{c|ccc}
				1/4 & 1/4 & 0 & 0 \\ 1/4 & 0 & 1/4 & 0 \\ 1 & 1/3 & 1/3 & 1/3 \\ \hline & 1/3 & 1/3 & 1/3
			\end{array}
			\label{eqn.IMEX2}
		\end{equation}
		
		\item SSP3(4,3,3)
		
		\begin{equation}
			\begin{array}{c|cccc}
				0 & 0 & 0 & 0 & 0 \\ 0 & 0 & 0 & 0 & 0 \\ 1 & 0 & 1 & 0 & 0 \\ 1/2 & 0 & 1/4 & 1/4 & 0 \\ \hline  & 0 & 1/6 & 1/6 & 2/3
			\end{array} \qquad
			\begin{array}{c|cccc}
				\alpha & \alpha & 0 & 0 & 0 \\ 0 & -\alpha & \alpha & 0 & 0 \\ 1 & 0 & 1-\alpha & \alpha & 0 \\ 1/2 & \delta & \eta & 1/2-\delta-\eta-\alpha & \alpha \\ \hline  & 0 & 1/6 & 1/6 & 2/3
			\end{array}
			\label{eqn.IMEX3}
		\end{equation}
		$\alpha=0.241694$, $\delta=0.060424$, $\eta=0.129153$
		
	\end{itemize} 
	\noindent
	The first order scheme \eqref{eqn.IMEX1} corresponds to the implicit Euler method and is stiffly accurate and stability preserving. Both properties are also exhibited by the second order scheme \eqref{eqn.IMEX2}, while the third order IMEX RK method \eqref{eqn.IMEX3} is not stiffly accurate.

	\section*{Declarations}
	
	\paragraph{Funding.} This research was partially funded by the Italian Ministry of Instruction, University and Research (MIUR) with PRIN Project 2017 No. 2017KKJP4X and PRIN Project 2017 No. 2017RSH3JY. Funding was also granted by the Free University of Bozen-Bolzano under the CRC project HM: "Hydropeaking mitigation".
	
	\paragraph{Conflicts of interest.} The authors declare that they have no conflict of interest.
	
	\paragraph{Availability of data and material.} Data and material are available upon reasonable request addressed to the corresponding author.
	
	\paragraph{Code availability.} The code is written in Matlab programming language and is available upon reasonable request addressed to the corresponding author.
	
	\paragraph{Ethics approval.} Not applicable.
	
	\paragraph{Consent to participate.} Not applicable.
		
	\paragraph{Consent for publication.} Not applicable.
	
	\section*{Data availability}
	The datasets generated during the current study are available from the corresponding author on reasonable request.

	\bibliography{biblio}   
	\bibliographystyle{plain}
	
\end{document}